%% file: arxiv_main.tex
\documentclass{article}

\usepackage{amsthm}
\usepackage{amsmath}

\usepackage{hyperref}
\input{package.tex}

\input{commands.tex}

\input{style.tex}

\newtheorem{remark}{Remark}

\newtheorem{lemma}{Lemma}

\usepackage{PRIMEarxiv}

\usepackage[utf8]{inputenc} 



\title{Network Inpainting via Optimal Transport
}
\author{
  Enrico Facca\\
  Department of Mathematics\\
  University of Bergen\\
  Bergen, Norway \\
  \texttt{enrico.facca@uib.no} \\
  \And
  Jan Martin Nordbotten\\
  Department of Mathematics\\
  University of Bergen\\
  Bergen, Norway\\
  \texttt{jan.nordbotten@uib.no}\\
  \And
  Erik Andreas Hanson\\
  Department of Computer science, Electrical engineering and Mathematical sciences\\
  Western Norway University of Applied Sciences\\
  Bergen, Norway\\
  \texttt{erik.andreas.hanson@hvl.no} \\
}

\begin{document}
\maketitle

\begin{abstract}
\input{abstract}
\end{abstract}

\newcommand{\sep}{\and}
\keywords{
  Inpainting  \sep
  Optimal transport \sep
  Branched transport \sep
  Image processing \sep
  Inverse problem \sep
}

\input{paper}

\appendix

\input{appendix}

\section*{Acknowledgments}
E.F. was founded by the European Union, under the Marie-Curie Action, project number 104109101.

\bibliographystyle{siam}  
\bibliography{biblio}

\end{document}

%% file: package.tex
\usepackage{graphicx}
\usepackage{array}

\usepackage{amssymb}

\usepackage[export]{adjustbox}
\usepackage{cleveref}
\usepackage{placeins}

\usepackage{xinttools}
\usepackage{grffile} 
\usepackage{caption}
\usepackage{subcaption}
\usepackage{color}
\usepackage{tikz}
\usepackage{multirow}

%% file: commands.tex
\newfont{\NUMBERS}{msbm8 scaled\magstep1}

\newcommand{\REAL}{\mbox{\NUMBERS R}}

\usepackage{ifthen}

\newcommand{\Vect}[2][]
{
  \ifthenelse{\equal{#1}{}}
  {\boldsymbol{#2}}
  {{#2}_{#1}}
}
\newcommand{\Matr}[2][]
{
  \ifthenelse{\equal{#1}{}}
  {\boldsymbol{#2}}
  {{#2}_{#1}}
}

\newcommand{\dx}{\, dx}

\DeclareMathOperator{\Div}{\operatorname{\text{div}}}
\DeclareMathOperator{\Grad}{\nabla}

\newcommand{\Dt}[1]{\partial_{t}{#1}}

\newcommand{\Obs}[1]{{#1}_{\text{obs}}}
\newcommand{\Rec}[1]{{#1}_{\text{rec}}}
\newcommand{\True}[1]{#1_{\text{true}}}
\newcommand{\Opt}[1]{{#1}_{\text{opt}}}

\newcommand{\ImgTd}{\mathcal{I}}

\newcommand{\Discr}{\mathcal{D}}
\newcommand{\Img}{I}
\newcommand{\Reg}{\ensuremath{\mathcal{R}}}
\newcommand{\confidence}{W}

\newcommand{\Inpaint}{\mathcal{J}}
\newcommand{\InpaintReg}{\mathcal{R}}

\newcommand{\Lyap}{\Reg}

\newcommand{\Wmass}{\ensuremath{\mathcal{M}_{\gamma}}}
\newcommand{\Ene}{\ensuremath{\mathcal{E}}}

\newcommand{\Dim}{d}
\newcommand{\Domain}{\ensuremath{\Omega}}

\newcommand{\tstep}{k}

\newcommand{\Deltat}[1][]{
  \ifthenelse{\equal{#1}{}}
  {\Delta t_{\tstep}}
  {\Delta t_{#1}}
  }

\newcommand{\MeshPar}{h}
\newcommand{\Cell}[1][]{\ifthenelse{\equal{#1}{}}{T}{T_{#1}}}
\newcommand{\Tsymb}{\mathcal{T}}
\newcommand{\Triang}[1][]
{
  \ifthenelse{\equal{#1}{}}
  {\Tsymb}
  {\Tsymb_{#1}}
}
\newcommand{\TriangH}{\Triang[\MeshPar]}

\newcommand{\Lspacechar}{L}
\newcommand{\Lspace}[2][]{
  \ifthenelse{\equal{#1}{}}
  {\Lspacechar^{#2}}
  {\Lspacechar_{#1}^{#2}}
}

\newcommand{\Tdens}{\mu}

\newcommand{\TdensIni}{\Tdens_0}
\newcommand{\TdensH}{\Tdens_{\MeshPar}}

\newcommand{\RecTdensH}{\Tdens_{\text{rec},h}}
\newcommand{\PouTdens}{\Tdens_{\text{Pou}}}
\newcommand{\OptTdensH}{\Tdens_{\text{opt},h}}
\newcommand{\OptTdens}{\Opt{\Tdens}}
\newcommand{\Pot}{u}
\newcommand{\PotH}{\Pot_{\MeshPar}}

\newcommand{\Vel}{v}

\newcommand{\Dirac}[1]{\delta_{#1}}

\newcommand{\Of}[1]{(#1)}

\newcommand{\Graph}{\mathcal{G}}
\newcommand{\Graphs}{\mathbf{G}}
\newcommand{\Vertices}{\mathcal{V}}
\newcommand{\Edges}{\mathcal{E}}

\newcommand{\Pbranch}{\alpha}

\newcommand{\Forcing}{f}

\newcommand{\Source}{\Forcing^{+}}
\newcommand{\Sink}{\Forcing^{-}}

\newcommand{\OT}{OT}

\newcommand{\BT}{BT}

\newcommand{\Lower}{\Tdens_{\min}}
\newcommand{\Mask}{D}

\newcommand{\DGzero}{DG^{0}}

\newcommand{\Scaling}{\alpha}
\newcommand{\PM}{PM}
\newcommand{\pmu}{\rho}
\newcommand{\PouseilleConstant}{\kappa}
\newcommand{\confmask}{\confidence_{\text{mask}}}
\newcommand{\TdensLift}{\Tdens_{+}}
\newcommand{\NIOT}{NIOT}

\newcommand{\tmpframe}[1]{\fbox{#1}}

\usepackage{changes}

%% file: style.tex
\newtheorem{Problem}{Problem}

\crefname{Problem}{Problem}{Problems}

\newcommand{\citet}{\cite}
\newcommand{\citep}{\cite}

%% file: abstract.tex
 In this work, we present a novel tool for reconstructing networks from corrupted images. The reconstructed network is the result of a minimization problem that has a misfit term with respect to the observed data, and a physics-based regularizing term coming from the theory of optimal transport. Through a range of numerical tests, we demonstrate that our suggested approach can effectively rebuild the primary features of damaged networks, even when artifacts are present.

%% file: paper.tex
\section{Introduction}
\label{sec:inpainting}
The problem of reconstructing images from corrupted data is a fundamental problem in science, arising in a wide range of problems, from medical imaging \cite{ben2021deep} to art restoration~\cite{guillemot2013image}. This reconstruction process is known in the literature of image processing as \emph{inpainting} \cite{ii2000,schonlieb2015partial,elharrouss2020image}. 

Inpainting problems may become particularly challenging when the image we are interested in represents a network and we want to recover certain its characteristics. For example, given the corrupted image shown in~\cref{fig:network-inpainting-obs}, we may want to restore the connected network shown in~\cref{fig:network-inpainting-true}. This requirement is fundamental if, for example, we want to simulate a fluid passing through the network, like the blood vessels in a medical image or the network of rivers in a satellite image. Networks impressed in images may be corrupted by loss of data, noise, and/or the presence of artifacts,
such as the example shown in~\cref{fig:network-inpainting-obs} or in many real-world problems \cite{kaufman1998two,update:2019,rs12172737}. If not properly reconstructed, the data cannot be used in numerical simulations, may require long manual preprocessing, or lead to wrong predictions.

\def \fraction {0.32}
\def \subfraction {0.99}
\begin{figure}
  \centering
  \begin{subfigure}[t]{\fraction\textwidth}
     \tmpframe{\includegraphics[width=\subfraction\columnwidth]{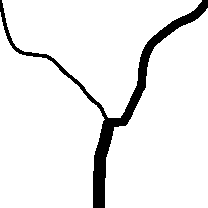}}
     \caption{Original $\True{\Img}$}
     \label{fig:network-inpainting-true}
  \end{subfigure}
  \qquad
  \begin{subfigure}[t]{\fraction\textwidth}
    \centering
    \tmpframe{\includegraphics[trim={0.9cm 2.7cm 0.9cm 2.7cm},clip,width=\subfraction\columnwidth]{{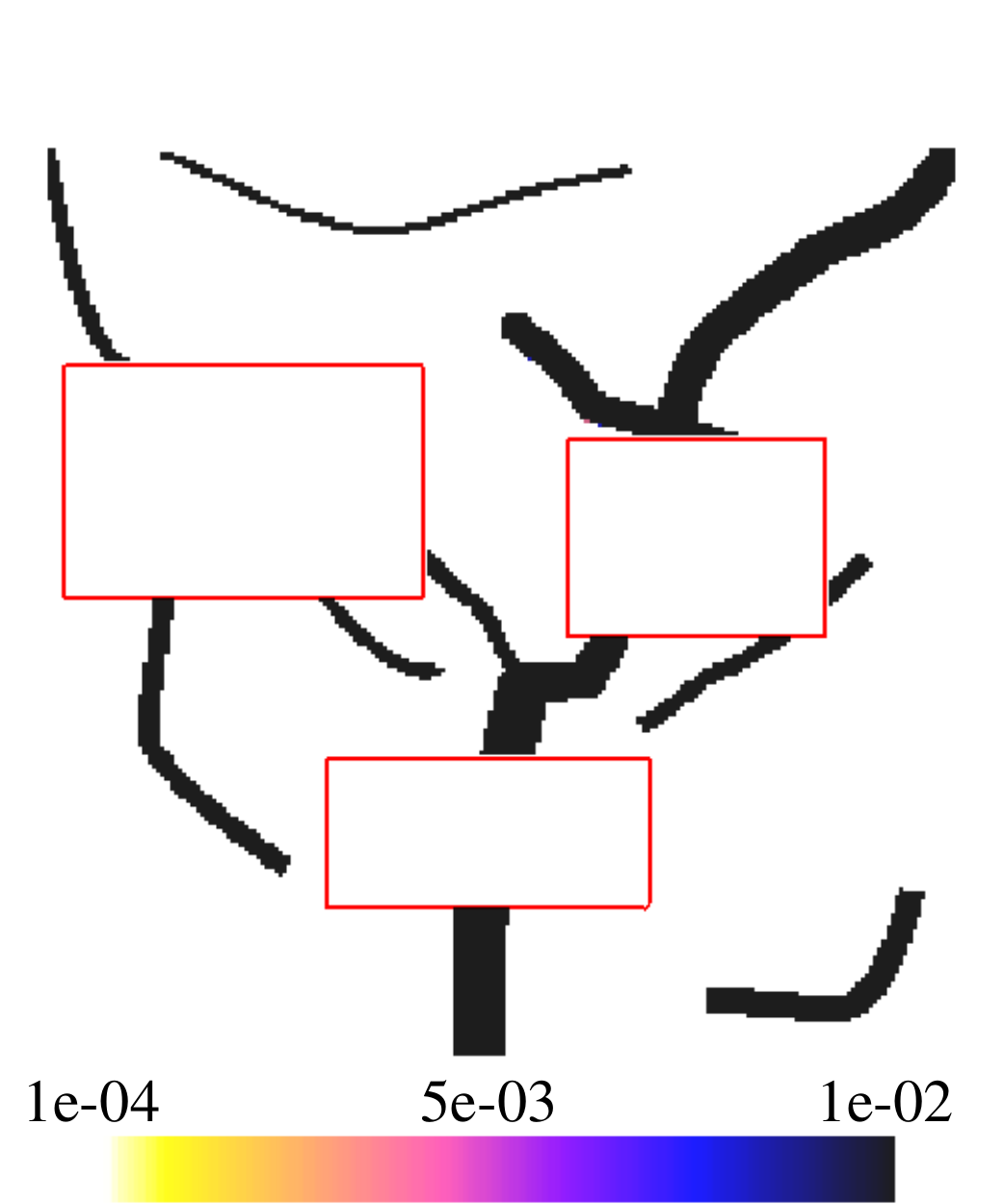}}}
    \caption{Observed $\Obs{\Img}$, corrupted by loss of data and by artifacts.}
    \label{fig:network-inpainting-obs}
  \end{subfigure}
  \\
     \begin{subfigure}[t]{\fraction\textwidth}
     \centering
     \tmpframe{\includegraphics[width=\subfraction\columnwidth]{{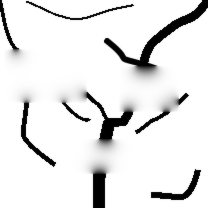}}}
     \caption{Harmonic $\Rec{\Img}$}
     \label{fig:network-inpainting-harmonic}
     \end{subfigure}
       \begin{subfigure}[t]{\fraction\textwidth}
     \centering
\tmpframe{\includegraphics[width=\subfraction\columnwidth]{{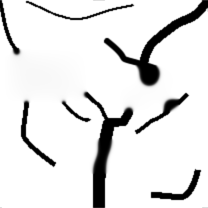}}}
  \caption{Cahn-Hilliard $\Rec{\Img}$ 
  }
  \label{fig:network-inpainting-ch}
  \end{subfigure}
  \begin{subfigure}[t]{\fraction\textwidth}
    \centering
\tmpframe{\includegraphics[width=\subfraction\columnwidth]{{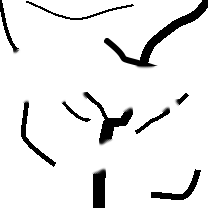}}}
  \caption{Transport $\Rec{\Img}$ 
  }
  \label{fig:network-inpainting-transport}
  \end{subfigure}
  \\
  \begin{subfigure}[t]{\fraction\textwidth}
  \centering
\tmpframe{\includegraphics[width=\subfraction\columnwidth]{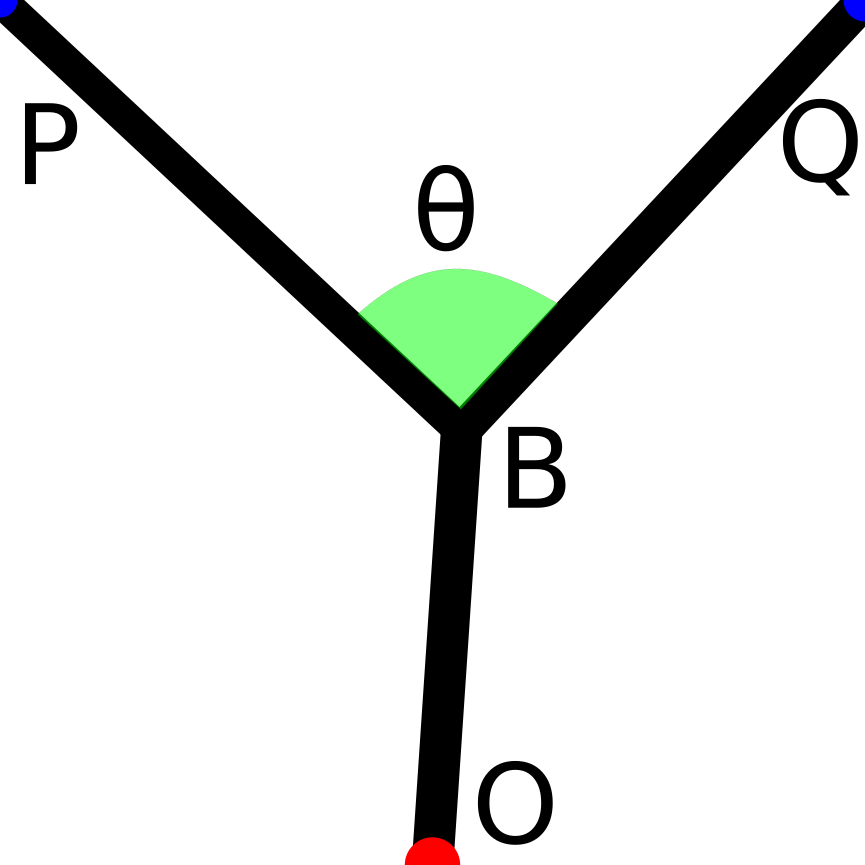}}
  \caption{Optimal transport network solving~\cref{prob:gilbert} for $\Source=\Dirac{O}$ (red) and $\Sink=1/3\Dirac{P}+2/3\Dirac{Q}$ (blue) for $\alpha=0.5$.}
  \label{fig:network-inpainting-btp}
  \end{subfigure}
  \qquad
  \begin{subfigure}[t]{\fraction\textwidth}
     \centering
    \tmpframe{\includegraphics[trim={0.9cm 2.7cm 0.9cm 2.7cm},clip,width=\subfraction\columnwidth]{{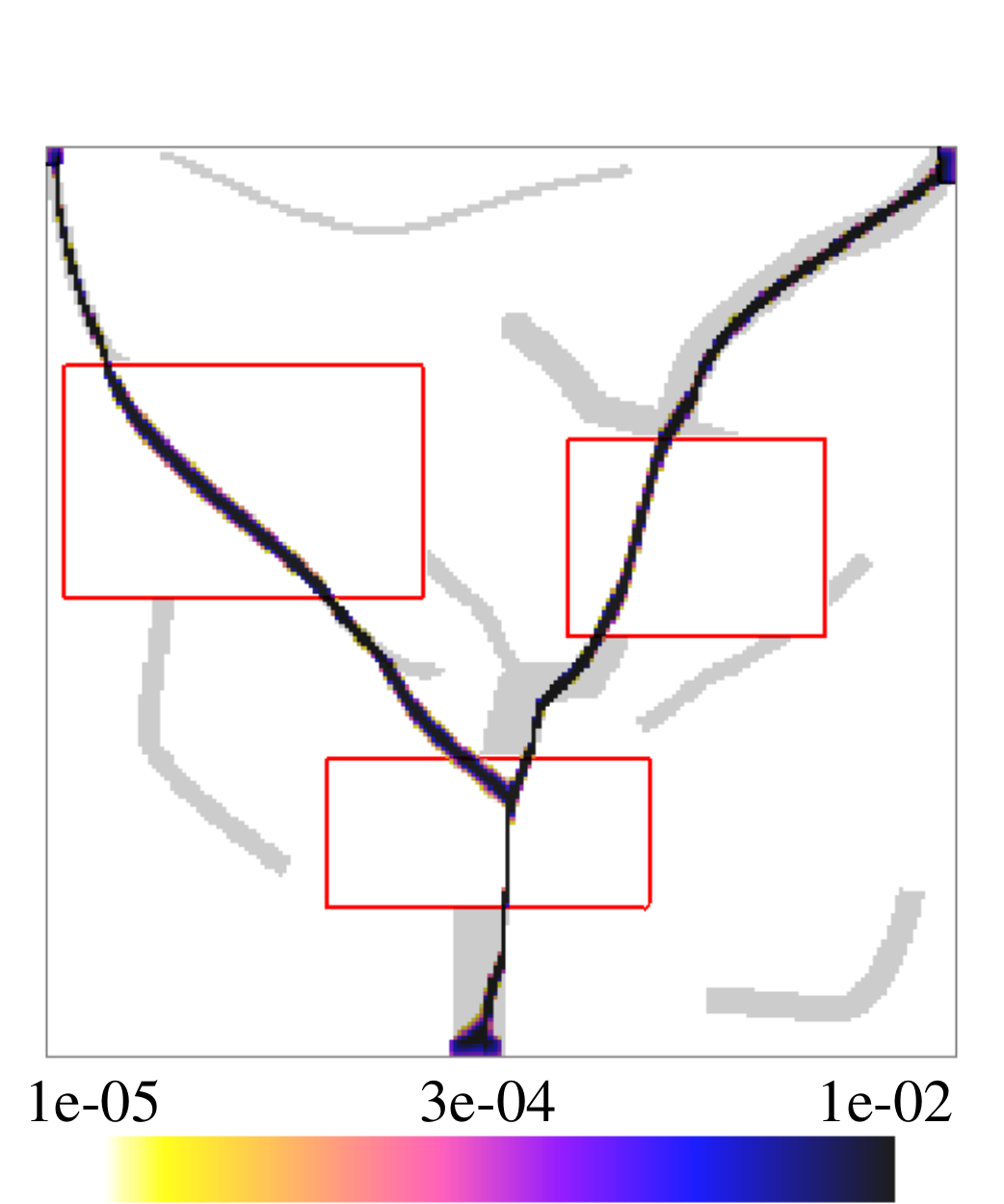}}}
    \caption{\NIOT{} $\Rec{\Img}$ ($\Obs{\Img}$ in light gray)}
  \label{fig:network-inpainting-rec}
  \end{subfigure}
  \caption{Schematic representation of the network reconstruction problem. In the upper panels, we report an image $\True{\Img}$ and a corrupted version $\Obs{\Img}$. In the central panels, we report the images reconstructed using three inpainting algorithms implemented in~\cite{ParSch2016}, the harmonic, Cahn-Hilliard, and transport inpainting methods. In~\cref{fig:network-inpainting-btp} we show the optimal network transporting $\Source$ to $\Sink$. \Cref{fig:network-inpainting-rec} reports a result obtained by the method presented in this paper taken from~\cref{fig:artifacts}.} 
  \label{fig:network-inpainting}
\end{figure}

The approaches used to solve inpainting problems can be divided into two main groups: machine learning-based methods and traditional methods. The former has shown versatile capabilities in recent years (see \cite{YI2019101552,MOCCIA201871,Jam2021} for relatively recent reviews on the topic). However, these methods rely on the availability of a large amount of high-quality data, while in many applications a transparent algorithmic approach is desired.

Among traditional inpainting methods, the most commonly used is to recover missing information by solving PDEs~\cite{nitzberg1993filtering,mm98,shen2002mathematical}. 
In this setting, an image is seen as a function $\Img:\Domain \subset \REAL^{\Dim}\to [0,\infty[$, in most problems $\Dim=2,3$. In the Y-shaped network in~\cref{fig:network-inpainting-true}, the domain $\Domain$ is a 2d square and the image/function $\True{\Img}$ is a binary map that indicates whether a point belongs to a network or not. 

Traditional methods include the so-called variational inpainting methods, where a reconstructed image $\Rec{\Img}$ is obtained by solving, for a fixed $\lambda>0$, the following minimization problem 
\begin{equation}
  \label{eq:inpainting-variational}
  \min_{\Img} \Inpaint(\Img) = 
  \InpaintReg(\Img) 
  + \lambda
  \Discr (\Img,\Obs{\Img}),
\end{equation}
where different discrepancy functionals $\Discr$ and regularization functionals $\InpaintReg$ can be used according to sought images. The parameter $\lambda>0$ determines the balance between the two terms. Standard discrepancy measures are weighted the $L^2$ and $H^{-1}$ norms. The most commonly used regularization functionals are the TV norm~\cite{rudin1992nonlinear,shen2002mathematical}, the Cahn-Hilliard functional~\cite{Bertozzi07,Burger09}, or the Mumford-Shah-Euler functional~\cite{esedoglu2002digital}. 

However, when the image describes a network, this type of regularization functional may not preserve the network structure. In~\cref{fig:network-inpainting-harmonic,fig:network-inpainting-ch,fig:network-inpainting-transport} we show how three inpainting algorithms (described in~\cite{schonlieb2015partial} and implemented in~\cite{ParSch2016}) fail in reconstructing the main topological features of the network in $\True{\Img}$. This is because they capture only local information and do not express any physical concepts. If we only consider natural networks such as blood vessels, rivers, and plant roots, these methods ignore our knowledge that we are dealing with networks whose main purpose is to transport nutrients and/or water, possibly also minimizing the energy spent for such transport. 

The root of this problem is the lack of a mathematical framework capable of describing a network in a convenient form for inpainting algorithms. The main idea of this paper is to address these issues thanks to recent advances in the optimal transport theory described in~\cite{Facca2021}, including a physics-based regularization term $\Reg$ in~\cref{eq:inpainting-variational}. This functional is written with respect to a nonnegative function $\Tdens$ that can be seen as a conductivity. Hence, we rewrite the minimization problem in~\cref{eq:inpainting-variational} in terms of this variable, solving the following problem 
\begin{equation}
\label{eq:inpainting-niot-intro}
  \min_{\Tdens\geq 0} 
  \Inpaint(\Tdens):=\Reg(\Tdens) + \lambda \Discr(\ImgTd(\Tdens), \Obs{\Img}),
\end{equation}
where $\ImgTd$ denotes a suitable map that transform a conductivity $\Tdens$ into an image $\Img$. The reconstructed image/network will simply be $\Rec{\Img}=\ImgTd(\Rec{\Tdens})$.
The variational problem in~\cref{eq:inpainting-niot-intro} represents the core of the methodology proposed in this paper, which we call Network Inpainting via Optimal Transport (NIOT).

The manuscript is organized as follows. In~\cref{sec:otp-btp} we give a brief overview of the optimal transport problem and, in particular, of the branched transport problem that studies the ramified structures we are interested in. We also present the functional $\Reg$ used in~\cref{eq:inpainting-niot-intro} and how it allows us to include branched transport ideas in the framework of image inpainting problems. This integration and the details of its numerical implementation are described in~\cref{sec:niot}. In~\cref{sec:experiments} we present a series of numerical experiments that present the capabilities and limitations of the proposed approach.

\section{Optimal and branched transport problems}
\label{sec:otp-btp}
The Optimal Transport (\OT{}) problem deals with finding the most efficient way to move a nonnegative measure $\Source$ into another measure of equal mass $\Sink$, assuming that the transport cost of moving one unit of mass is given by a function $c:\Omega\times \Omega \to \REAL$, where $\Omega\subset \REAL^{\Dim}$ is the domain of measures $\Source$ and $\Sink$. 

The typical transport cost is defined as $c(x,y)=|x-y|^{p}$, with $x,y\in \Domain$ and $p\geq 1$, for which a rich literature exists (see~\cite{Santambrogio:2015,Peyre2019} for complete overviews of theoretical and numerical advances on the topic). The fundamental solution to \OT{} problems with this type of cost is that the mass travels along straight lines \cite{villani2009optimal}. For example, in the problem reported in~\cref{fig:network-inpainting-btp}, the \OT{} solution follows a V-shaped path connecting the point $O$ to $P$ and $Q$, sending $1/3$ of the mass to $P$ and $2/3$ to $Q$. 

However, one may want to model transport phenomena in which it is more convenient to move the mass together, following a path such as the Y-shaped network reported in~\cref{fig:network-inpainting-btp}. The mass at the point $O$ first moves to the point $B$ and then splits along two different branches. The subarea of \OT{} that studies this type of branching structures as minimal energy solutions is called the Branched Transport (\BT{}) problem (see also~\cite{bernot2008optimal} for a general overview of the topic). 

The first \BT{} formulation can be traced back to the work of Gilbert in~\cite{Gilbert:1967}, but in this paper we use the formulation described in~\cite{xia:2003}. Before presenting it, it is useful to state the following simple lemma on the mass balance equation on graphs.
\begin{lemma}
\label{lem:div-graph}
Consider a graph $\Graph=(\Vertices,\Edges)$, with nodes $\Vertices$ and edges $\Edges$ that contains no loops (this type of graph is said to be acyclic). Fix an orientation on the edges $\Edges$ and let $\Vect{\Forcing}\in \REAL^{|\Vertices|}$ whose entries sum up to zero on each connected components of $\Graph$. Then there exists a unique flux $\Vel\in \REAL^{|\Edges|}$ that satisfies the following mass balance equation 
\begin{equation}
\label{eq:div-graph-pure}
\sum_{e\in \delta^{+}(k)} \Vect[e]{\Vel} - \sum_{e\in \delta^{-}(k)} \Vect[e]{\Vel} = f_k  \quad \forall k=1,\ldots,|\Vertices|,
\end{equation}
where $\delta^{+}(k)$ and $\delta^{-}(k)$ denote the set of edge entering and existing the $k$-th node (according to the fixed orientation).
\end{lemma}

The proof of this lemma can be found in~\cite{bapat2010graphs}. It is based on the facts that on each connected component of acyclic graphs the number of vertices is one plus the number of edges. We can now present the \BT{} problem as follows.
\begin{Problem}
\label{prob:gilbert}
Let us denote with $\Dirac{x}$ a Dirac measure centered at $x\in \Domain$.
Consider two measures
\begin{equation}
    \label{eq:forcing-gilbert}
\Source=\sum_{i=1}^n\Source_i\Dirac{x_i},
\quad
\Sink=\sum_{j=1}^m\Sink_j\Dirac{y_j},
\end{equation}
where:
\begin{itemize}
    \item $n,m$ are two positive integers;
    \item $\Source_i,\Sink_j$ are strictly positive real numbers that sum to zero;
    \item $x_i,y_j$ are distinct points in $\Omega$.
\end{itemize}
Consider the set $\Graphs$ of those graphs $\Graph=(\Vertices,\Edges)$ such that:
\begin{itemize}
  \item the graph $\Graph$ is a union of acyclic graphs;
  \item the set $\Vertices$ contains all nodes $x_i$ and $y_j$, and possibly extra nodes;
  \item to each edge $e\in \Edges$ we associate a length $\Vect[e]{\ell}$ given by the Euclidean distance between the connected nodes. We denote by $\Vect{\ell}\Of{\Graph}$ the vector in $\REAL^{|\Edges|}$ with all these lengths.
\end{itemize}

Consider the following vector $\Vect{\Forcing}\in \REAL^{|\Vertices|}$ with entries $\Vect[k]{\Forcing}$ given by
\begin{equation}
\label{eq:div-graph}
\Vect[k]{\Forcing}
  =
  \left\{
    \begin{aligned}
    &\Source_i &\text{if $k$-th node}=x_i,\\
    &-\Sink_j &\text{if $k$-th node}=y_j,
    \\
    & 0 &\text{otherwise},
    \end{aligned}
  \right.
\end{equation}
for $k=1,\ldots,|\Vertices|$. If a graph $\Graph\in \Graphs$ and falls under assumptions~\cref{lem:div-graph}, then there exists a unique flux that satisfies the mass balance equation~\cref{eq:div-graph-pure} with the forcing term $\Vect{\Forcing}$, which we denote by $\Vect{\Vel}\Of{\Graph}$.

Now, fix an exponent $0<\Pbranch<1$. We want to find the optimal graph $\Opt{\Graph}\in \Graphs$ that minimizes the energy 
\begin{equation}
    \label{eq:gilbert-energy}
    E^{\alpha}(\Graph):=\sum_{e\in \Edges}|\Vect[e]{\Vel}\Of{\Graph}|^\alpha\Vect[e]{\ell}\Of{\Graph}.
\end{equation}
\end{Problem}

The presence of the concave and sub-additive term $|\cdot|^{\alpha}$ favors graphs with fewer edges but a higher flux rate, making the ramified structures optimal.~\Cref{prob:gilbert} admits a solution since there is a finite number of possible combinations of graphs that meet its requirements~\cite[Prop. 2.1, 2.2]{xia2015motivations}. However, the number of these combinations increases exponentially with the total number of nodes $x_i,y_j$, making the \BT{} problem NP hard.

An explicit description of the optimal graph can be found only in simple cases, such as the problem in~\cref{fig:network-inpainting-btp}, where there is an analytic formula for the branching point $B$ given in~\cite{bernot2008optimal}. This formula also shows that the branching angle $\theta$ between the edges $PB$ and $QB$ increases as the exponent $\alpha$ passes from one to zero, reaching the maximum value of $120^{\circ}$. More generally, the minimal branching angle of the network is bounded from below by a quantity that decreases with the exponent $\alpha$~\cite[Proposition 2.1]{xia2011boundary}.

Extensions of~\cref{prob:gilbert} to general measures $\Source,\Sink$ have been proposed in recent years~\cite{maddalena2003variational,brasco2011benamou,xia:2003}, with the latter being the most closely related to~\cref{prob:gilbert}. 

The theory of \BT{} provides a unified mathematical framework for various studies that model natural networks as optimal transport solutions. Examples include rivers~\cite{rigon1993optimal,Rodriguez-Rinaldo:2001}, blood vessels~\cite{Murray1926,banavar1999size}, and plant roots~\cite{McCulloh2003}. However, two fundamental issues limit the use of this theory in applied problems and, more specifically, in inpainting problems considered in this paper.

The first issue is merely computational. Given the NP hard nature of~\cref{prob:gilbert} there is no hope of having an efficient numerical method to solve the \BT{} problem, neither for $\Source$ and $\Sink$ given in~\cref{eq:forcing-gilbert} nor for general measures. All methods found in the literature work only for small numbers of nodes~\cite{xia-numerics,lippmann2022theory}, converge toward local minima~\cite{oudet2011modica,Facca2021}), or require additional assumptions on the transported measures $\Source$ and $\Sink$~\cite{Dirks2022}. 

The second issue is the incompatibility between~\cref{prob:gilbert} and the variational inpainting problem in~\cref{eq:inpainting-variational}. The first is written in terms of graph, while in the second, the information on the observed image represents some form of intensity variation in space and is typically stored in pixels/voxels. There is no simple way to integrate these two problems. In the following section, we will introduce an optimization problem closely related to the \BT{} problem, where this integration can be done more naturally.

\subsection{A reformulation of the \BT{}}
\label{sec:dmk}
In~\cite{Facca2021}, the following problem has been proposed.
\begin{Problem}
  \label{prob:dmk}
  Consider $\Source,\Sink$ being two nonnegative densities in $L^2(\Domain)$ having the same mass. Fix $\Lower>0$ and $\gamma\in]0,1]$. Find an optimal ``conductivity'' $\OptTdens:\Domain\to [0,+\infty[$ solving
  \begin{equation}  
    \label{eq:lyap}
    \min_{\Tdens\geq 0}
    \Lyap(\Tdens) :=
    \underbrace{
    \int_{\Domain}\frac{(\Tdens+\Lower)|\Grad\Pot\Of{\Tdens}|^2}{2}\dx
    }_{
    =:\Ene(\Tdens)
    } 
    +
    \underbrace{
    \int_{\Omega}
    \frac{\Tdens^{\gamma}}{2\gamma}\dx
    }_{
=:\Wmass(\Tdens)
    },
\end{equation}
where, for any given $\Tdens:\Omega\to[0,\infty[$, the function $\Pot\Of{\Tdens}$ denotes the solution $\Pot:\Domain\to \REAL$ of the equation 
\begin{equation}
  \label{eq:poisson}
    \left\{
    \begin{aligned}
    &-\Div((\Tdens+\Lower)\Grad \Pot)=\Source-\Sink,\\
    &(\Tdens+\Lower)\Grad \Pot \cdot n_{\partial \Domain}=0.
    \end{aligned}
    \right.
\end{equation}
\end{Problem}
Here $\Lower>0$ is a small scalar that ensures well-posedness of the PDE~\cref{eq:poisson}. In this paper, we fix $\Lower=1e-8$, while in~\cite{Facca2021} $\Lower$ was set to zero. Similar problem formulations can be found also in~\cite{Haskovec-et-al:2015,Albi2017}. \Cref{prob:dmk} with $\gamma=1$ and $\Lower=0$ is equivalent with the \OT{} problem with a cost equal to the Euclidean distance, as shown in~\cite{facca2020numerical}.

Intuitively, the support of a function $\Tdens$ in~\cref{eq:lyap} describes a ``pathway'' connecting $\Source$ to $\Sink$. The divergence constraint ensures its linkage. An optimal solution $\OptTdens$ of~\cref{prob:dmk} provides the best trade-off between energy $\Ene(\Tdens)$, which is the Joule energy dissipated via potential flow, and cost $\Wmass(\Tdens)$, which represents the nonlinear cost of building the ``pathway'' $\Tdens$.  If $0<\gamma<1$,  this cost is concave and encourages the concentration of the support of $\Tdens$ in small channels with high values of $\Tdens$. Increasing $\gamma$ from zero to one, this effect concentration effect is excepted to weaken, passing from networks having few channels, wide branching angles, and high conductivities to networks with more channels, narrower branching angles, and reduced conductivities. 

These intuitions are supported by a series of numerical experiments presented in~\cite{Facca2021} in which a discrete approximation of the minimization problem in~\cref{eq:lyap} is obtained using finite elements, solution denoted hereafter by $\OptTdensH$. In~\cref{fig:dmk} we show the approximate $\OptTdensH$ obtained for the measures $\Source$ and $\Sink$ in~\cref{fig:network-inpainting-btp} using different values of $\gamma$ and refining an initial grid $\Triang_{h}$ used to discretizated the problem (see~\cref{sec:algorithm} for the details of the numerical solution of~\cref{prob:dmk}). 

{
\def \fraction {0.225}
\def \extrah {0.0em}
  \begin{figure}
    \centering
    \begin{tabular}{|c|@{\hspace{\extrah}}c@{\hspace{\extrah}}|@{\hspace{\extrah}}c@{\hspace{\extrah}}|@{\hspace{\extrah}}c@{\hspace{\extrah}}|c}
      \cline{1-4}
      & $\Triang_{h}$ & $\Triang_{h/2}$ &  $\Triang_{h/4}$ &
      \\
      \cline{1-4}
      \raisebox{-.5\normalbaselineskip}[0pt][0pt]{$\gamma=0.8$}
      &
      \tmpframe{\adjustbox{valign=m,vspace=0pt}{\includegraphics[trim={4.3cm 0.9cm 4.3cm 0.9cm},clip,width=\fraction\columnwidth,valign=t]{{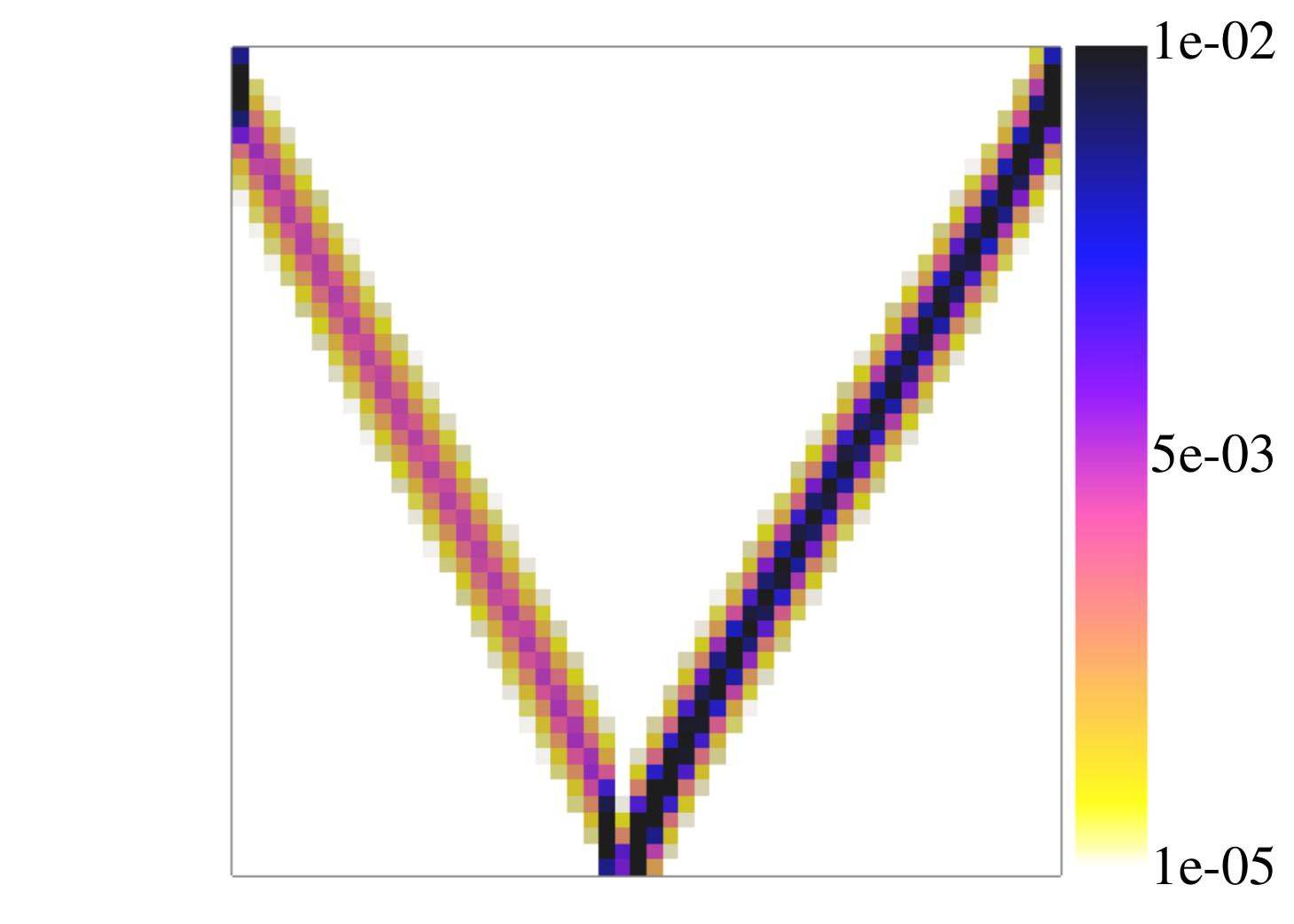}}}}
      &
      \tmpframe{\adjustbox{valign=m,vspace=0pt}{\includegraphics[trim={4.3cm 0.9cm 4.3cm 0.9cm},clip,width=\fraction\columnwidth,valign=t]{{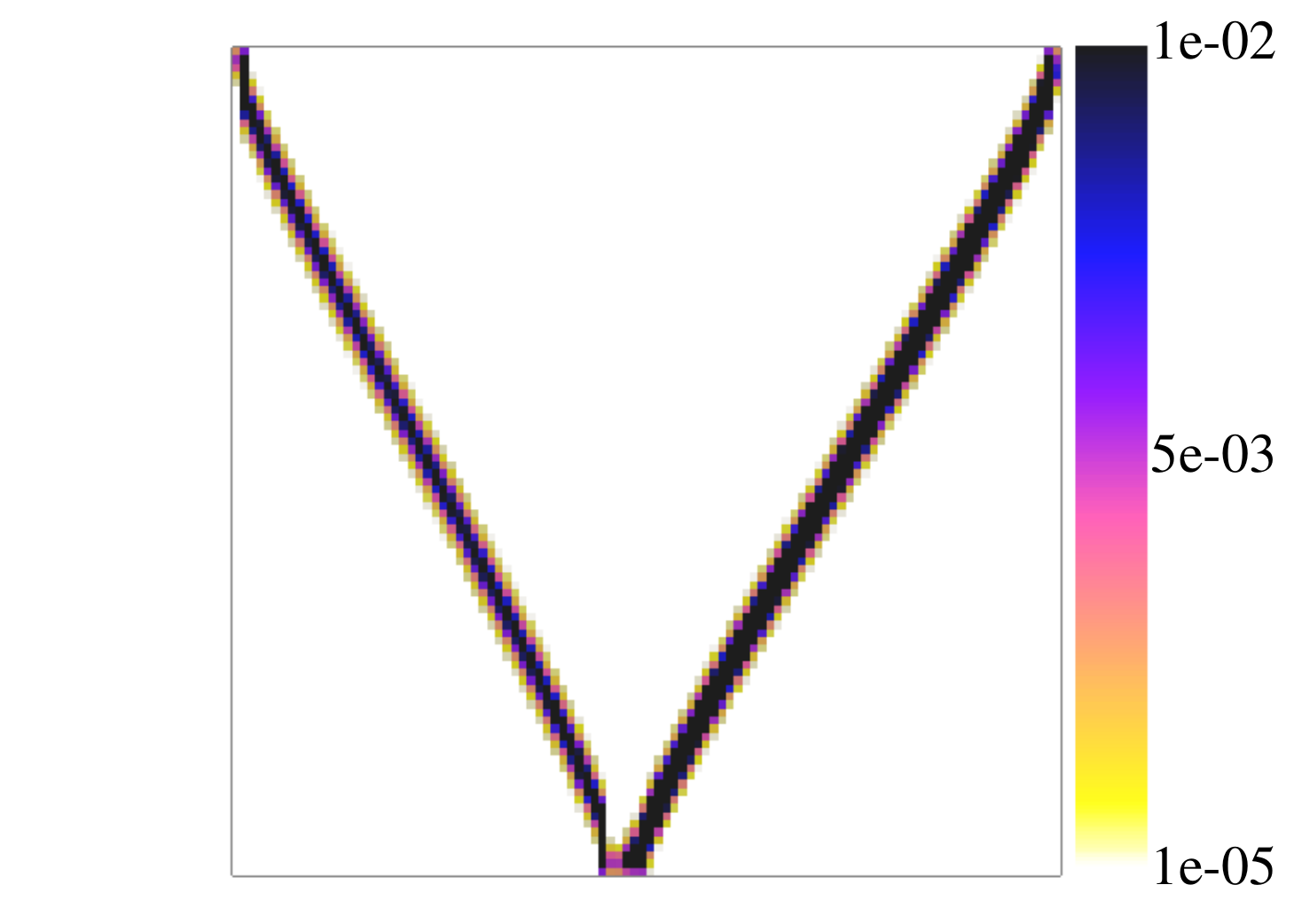}}}}
      &
      \tmpframe{\adjustbox{valign=m,vspace=0pt}{\includegraphics[trim={4.3cm 0.9cm 4.3cm 0.9cm},clip,width=\fraction\columnwidth,valign=t]{{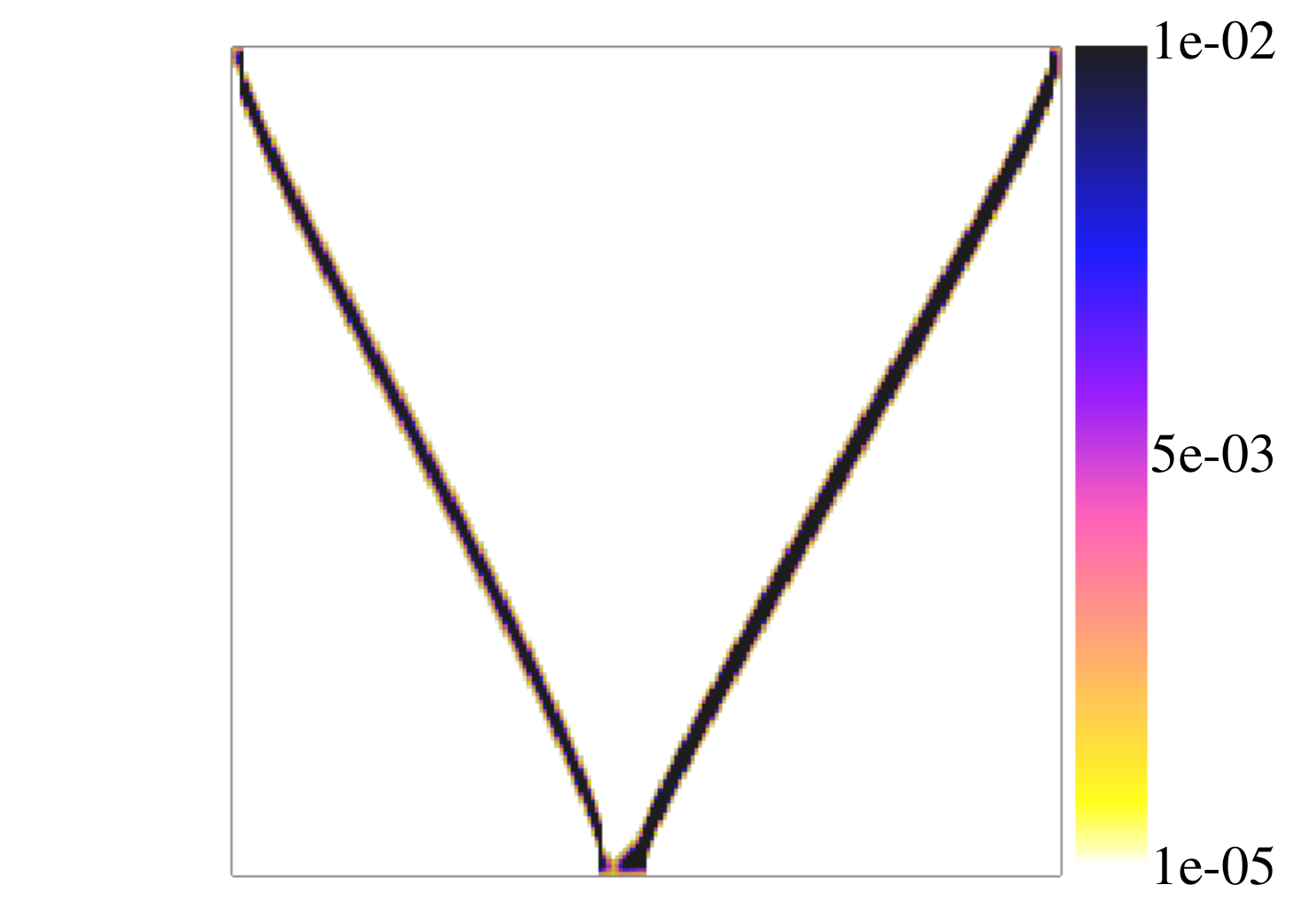}}}}
      &
      \multirow{3}{*}{ \adjustbox{valign=m,vspace=0pt}{\includegraphics[trim={20.3cm 0cm 0.3cm 0cm},clip,width=0.08\columnwidth,valign=t]{{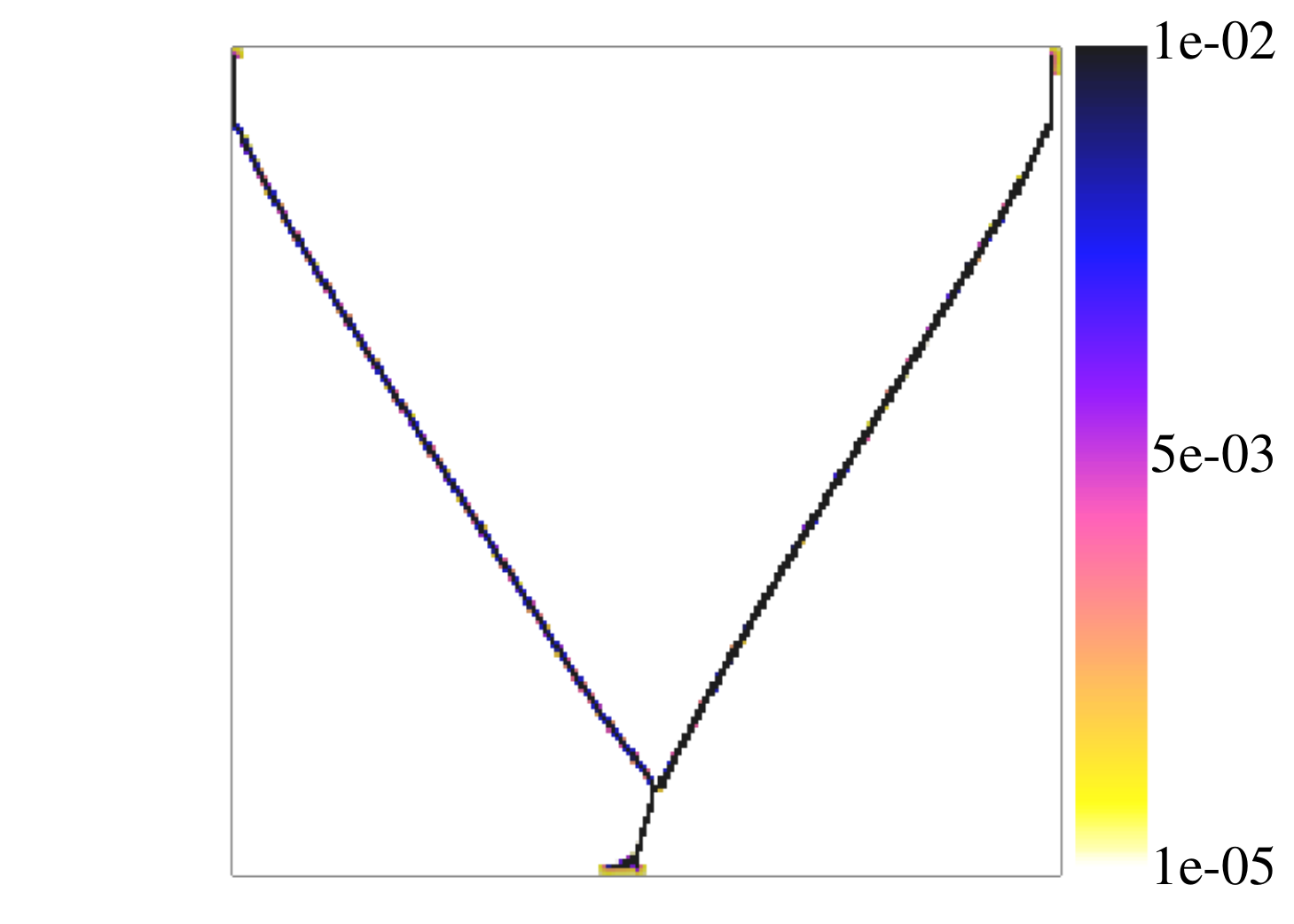}}}}
      \\
      \cline{1-3}
      \raisebox{-.5\normalbaselineskip}[0pt][0pt]{$\gamma=0.5$}
      &
      \tmpframe{\adjustbox{valign=m,vspace=0pt}{\includegraphics[trim={4.3cm 0.9cm 4.3cm 0.9cm},clip,width=\fraction\columnwidth,valign=t]{{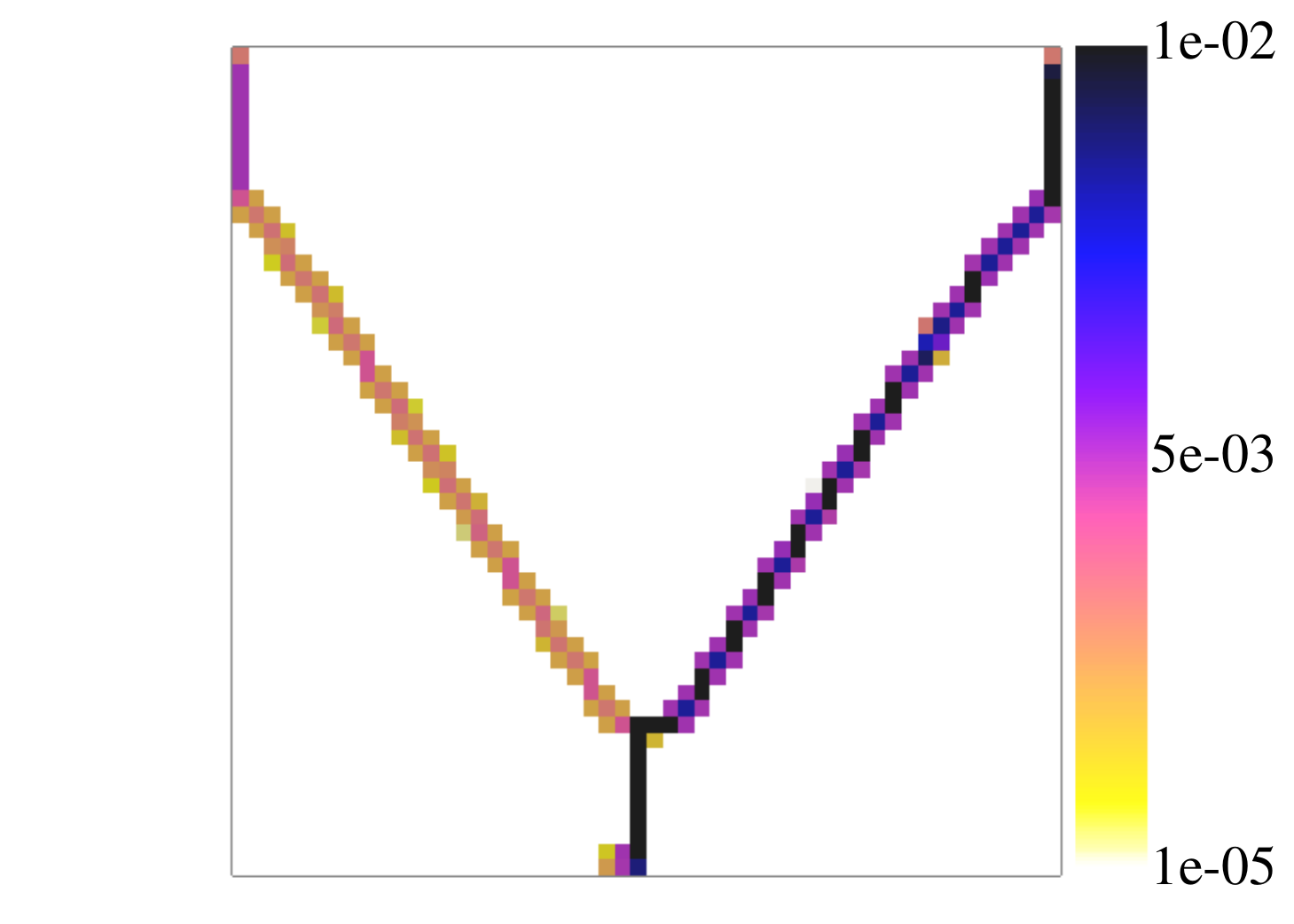}}}}
      &
      \tmpframe{\adjustbox{valign=m,vspace=0pt}{\includegraphics[trim={4.3cm 0.9cm 4.3cm 0.9cm},clip,width=\fraction\columnwidth,valign=t]{{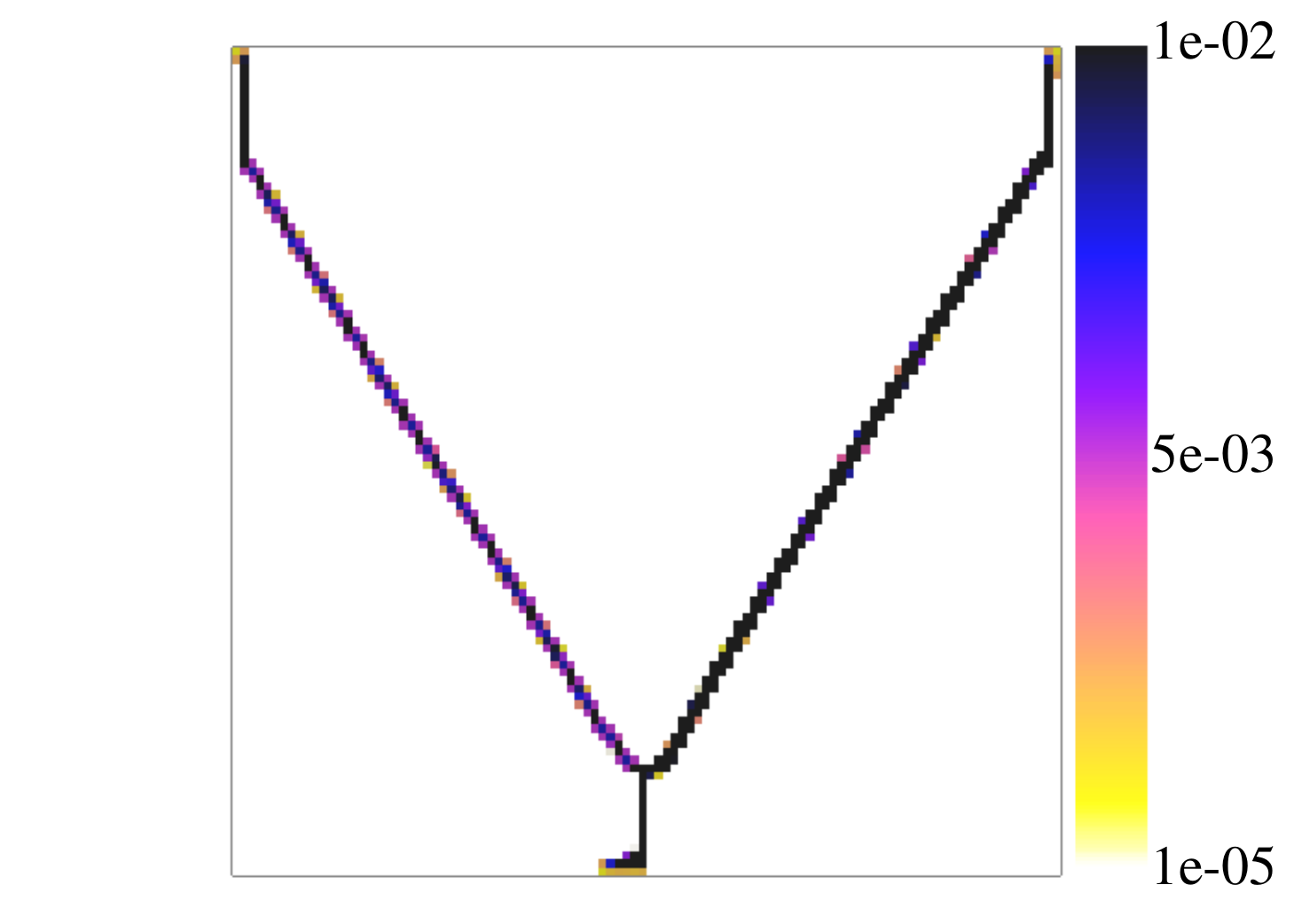}}}}
      &
      \tmpframe{\adjustbox{valign=m,vspace=0pt}{\includegraphics[trim={4.3cm 0.9cm 4.3cm 0.9cm},clip,width=\fraction\columnwidth,valign=t]{{imgs/y_net/mask_medium/matrix_nref2_femDG0DG0_gamma5.0e-01_wd0.0e+00_wr0.0e+00_netnetwork_ini0.0e+00_confMASK_mu2iidentity_scaling1.0e+01_methodte_tdens.pdf}}}}
      &
      \\
      \cline{1-4}
      \raisebox{-.5\normalbaselineskip}[0pt][0pt]{$\gamma=0.2$}
      &
      \tmpframe{\adjustbox{valign=m,vspace=0pt}{\includegraphics[trim={4.3cm 0.9cm 4.3cm 0.9cm},clip,width=\fraction\columnwidth,valign=t]{{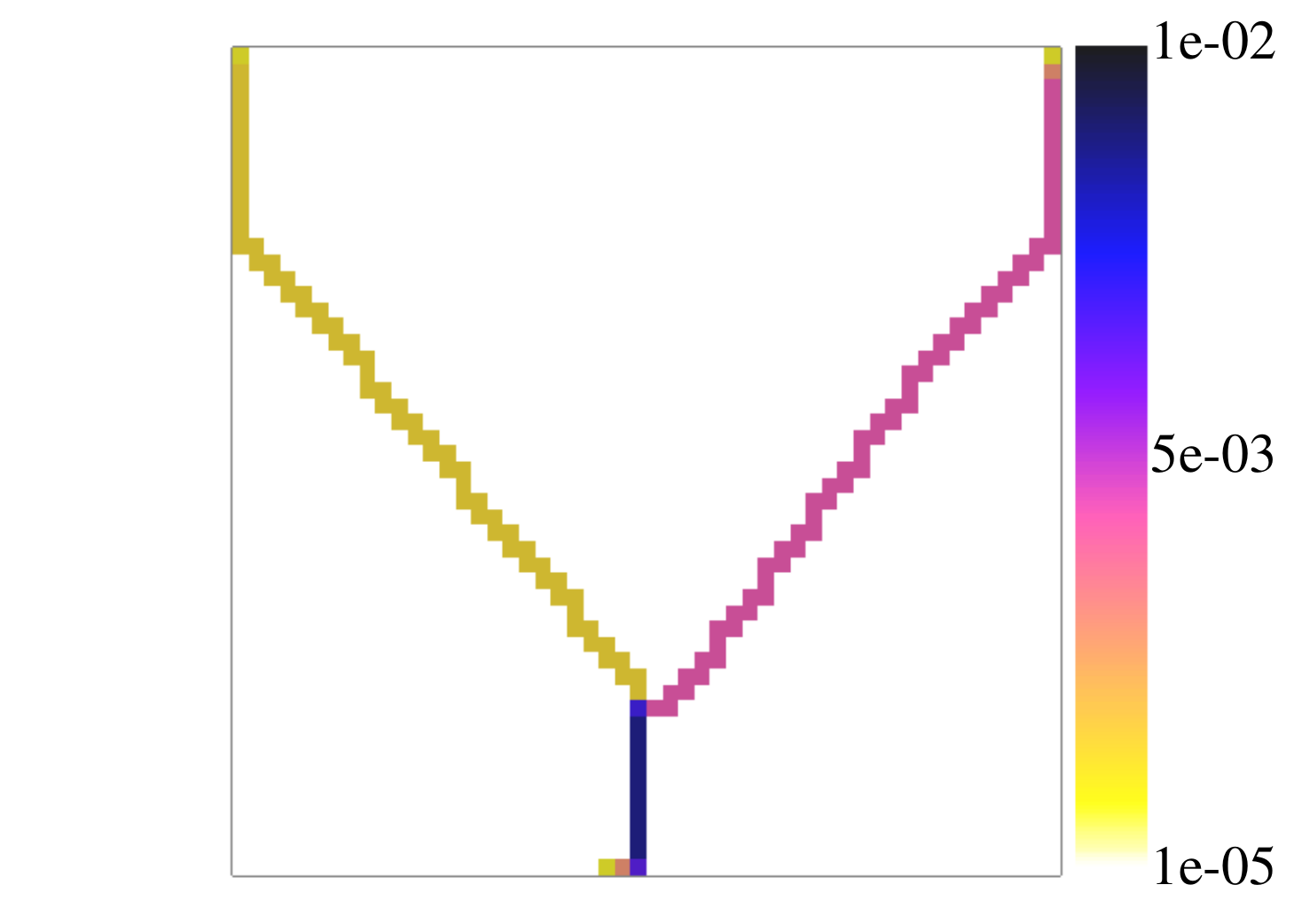}}}}
      &
      \tmpframe{\adjustbox{valign=m,vspace=0pt}{\includegraphics[trim={4.3cm 0.9cm 4.3cm 0.9cm},clip,width=\fraction\columnwidth,valign=t]{{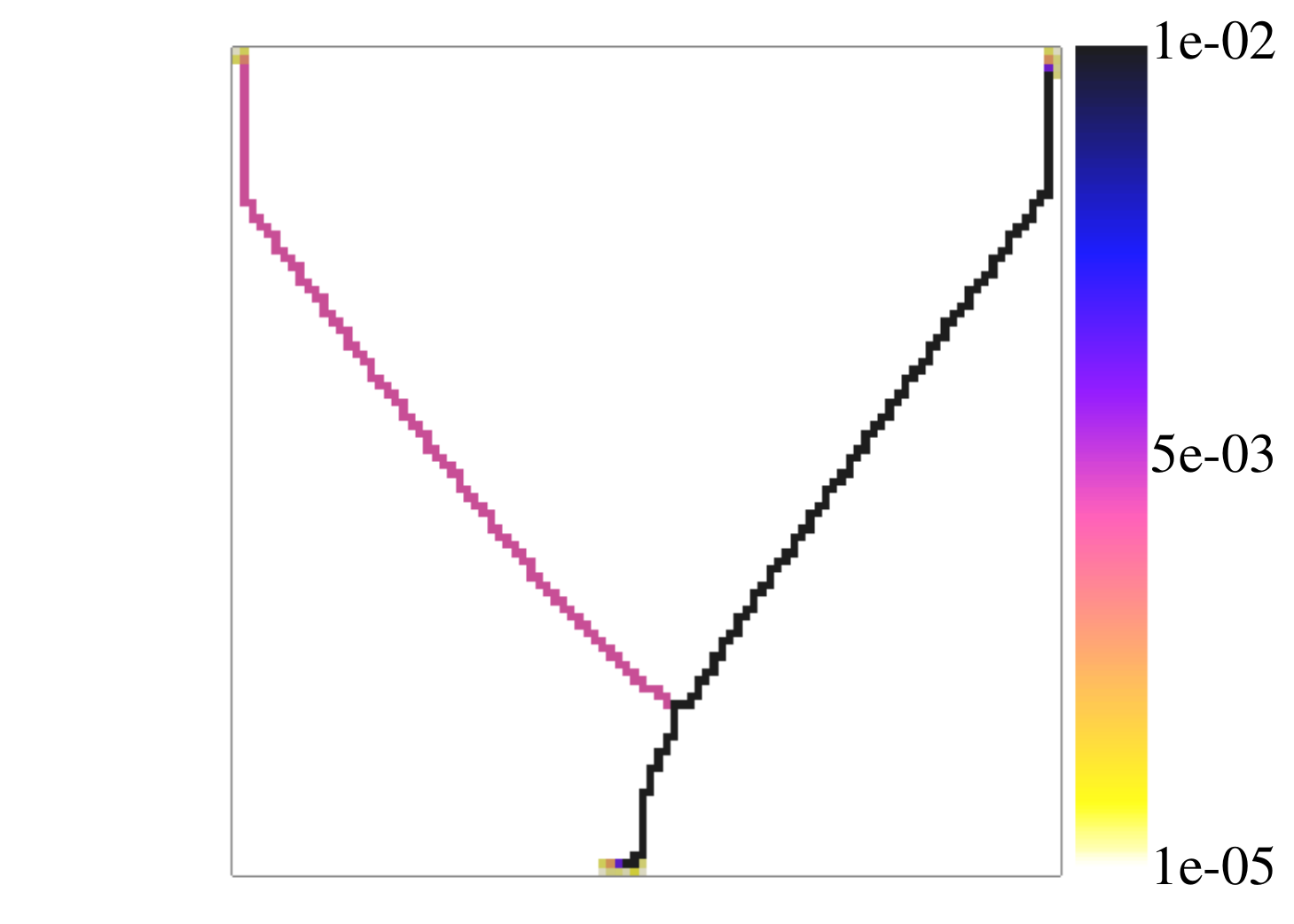}}}}
      &
      \tmpframe{\adjustbox{valign=m,vspace=0pt}{\includegraphics[trim={4.3cm 0.9cm 4.3cm 0.9cm},clip,width=\fraction\columnwidth,valign=t]{{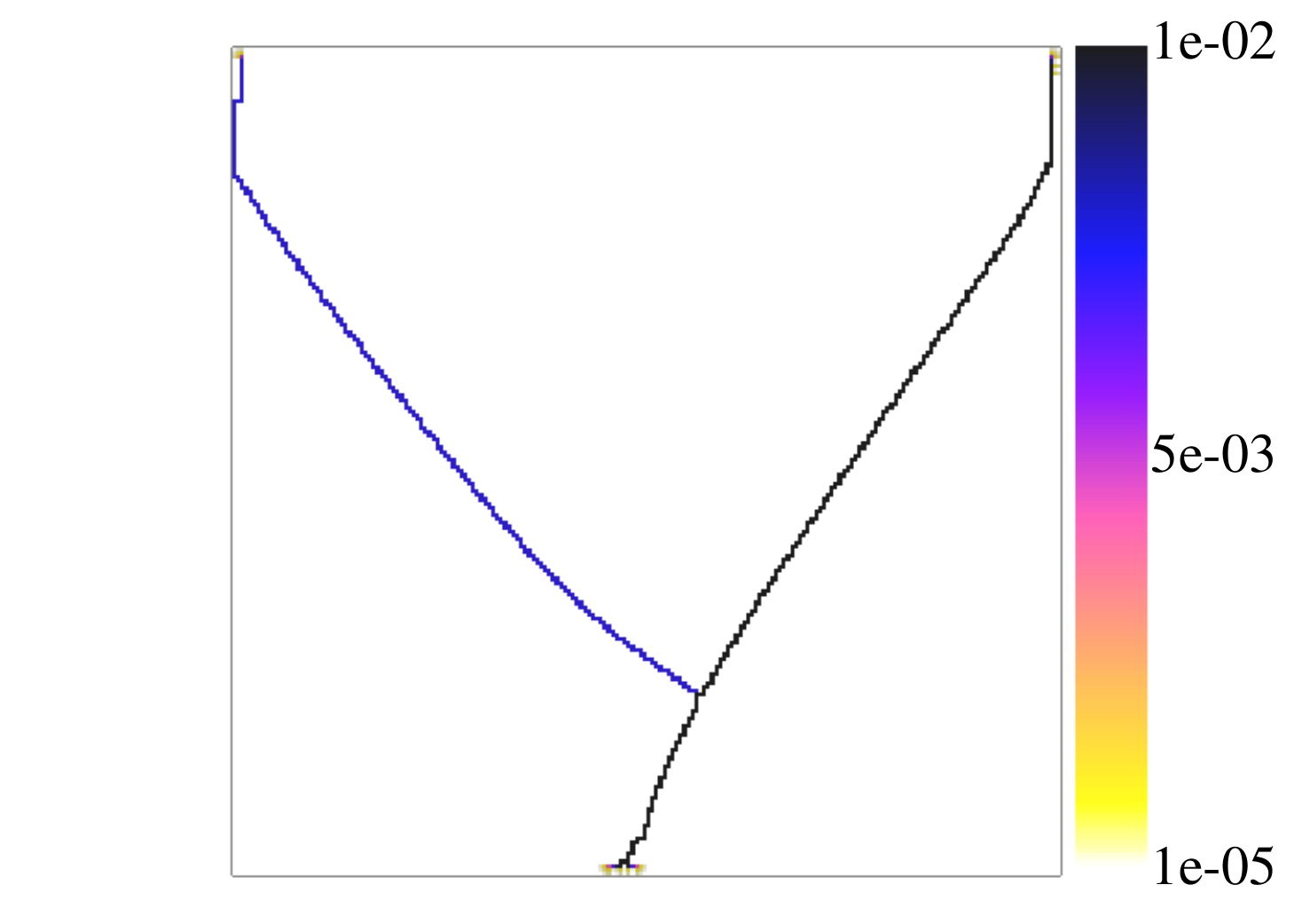}}}}
      &
      \\
      \cline{1-4}
    \end{tabular}
    \caption{Spatial distribution of $\DGzero$-approximation of the $\OptTdensH$ solving~\cref{prob:dmk} for the $\Source$ and $\Sink$ shown in~\cref{fig:network-inpainting-btp} using $\gamma=0.8,0.5,0.2$ (top to bottom) and refining the initial $52\times 52$ grid $\Triang_{h}$(left to right)}
    \label{fig:dmk}
  \end{figure}
}

We see that for $\gamma=0.8$ a V-shaped network is obtained, while for $\gamma=0.5$ and $\gamma=0.2$ a Y-shaped network is obtained. We can also see that $\OptTdensH$ depends on the mesh $\TriangH$ used to discretize the problem. In fact, its support tends to concentrate on narrow ``channels'' composed of a sequence of single cells in the mesh $\TriangH$, indicating convergence towards a solution having a support composed of a union of 1d lines. Similar results were shown in~\cite[Fig. 5-6-7]{Facca2021}. Due to these phenomena, it is not clear how to recover the infinite-dimensional problem in~\cref{eq:lyap} from its discrete counterpart. Therefore, there is no formal connection between~\cref{prob:dmk} and classical \BT{} formulations. 

However, this matter has a slight impact on the variational inpainting problem that we are focusing on. In such problems, the grid is given by the format used to store $\Obs{\Img}$, and we do not concern ourselves with the issue of cells with size going to zero. More crucially, the mesh appears to have a minor influence on the overall network topology that is formed by the support of $\OptTdensH$, as shown in~\cref{fig:dmk} or in~\cite[Fig. 8]{Facca2021}. For a given grid, \cref{prob:dmk} offers the substantial advantage of defining an optimal network as the support of the function $\OptTdensH$. This facilitates a seamless incorporation of \BT{} principles into inpainting problems, where data are typically stored as pixels/voxels.

Another issue related to~\cref{prob:dmk}, either in its continuous or discrete form, arises from the fact that the functional in~\cref{eq:lyap} is nonconvex. In this setting, there is no assurance that the gradient descent-based technique employed in~\cite{Facca2021} will discover a global minimum. As a result, any potential minimizer that we identify is merely local and depends on the initial data. Therefore, when we mention optimal conductivity $\OptTdensH$, we will specifically refer to the local minimum achieved using the gradient descent method and will detail the initial data utilized.

\section{Network inpainting via optimal transport}
\label{sec:niot}

Two additional building blocks are required to define the NIOT problem in~\cref{eq:inpainting-niot-intro}.
\begin{itemize}
    \item A map $\ImgTd$ that transforms a nonnegative conductivity $\Tdens$ into a nonnegative image $\Img$, so that $\Img=\ImgTd(\Tdens)$.
    \item A discrepancy term suitable for the network inpainting problem.
\end{itemize}
The forthcoming two sections are devoted to defining the desired traits and building these two components.

The fundamental concept is that by minimizing the functional $\Reg(\Tdens)$, we ensure simultaneously that the support of $\Tdens$ connects $\Source$ and $\Sink$, as a result of the divergence constraint in~\cref{eq:poisson}, and that it approximates a 1d structure. 
The allocation of $\Source$ and $\Sink$ is our indirect method to integrate our knowledge that the observed network is a structure that transports resources into the reconstruction problem. On the other hand, by introducing the discrepancy term, we force the reconstructed network $\Rec{\Img}$ to closely match the observed data, encompassing other physics-based information not encapsulated in~\cref{prob:dmk}. This suggested method can be viewed as a model-driven inverse problem, similar to those described in~\cite{arridge_maass_oktem_schonlieb_2019}.

\subsection{Discrepancy term $\Discr$}
\label{sec:discrepancy}
The discrepancy term considered in this paper is the following weighted $L^2$-distance 
\begin{equation}
    \label{eq:discrepancy}
\Discr(\Img, \Obs{\Img})
  :=
  \frac{1}{2}
  \int_{\Domain} 
   (\Img - \Obs{\Img})^2 \confidence \dx,
\end{equation}
where the weight $\confidence:\Domain \to [0,1]$ measures our confidence to the observed data in each portion of the domain. In this paper, we consider two types of weights.

The first weight $\confidence$ can be defined only if there exists a region $\Mask\subset \Domain$, hereafter denoted as mask, where we know that the network is corrupted. For example, in~\cref{fig:network-inpainting-obs} the mask $\Mask$ is the region within the red rectangles. Some inpainting algorithms (for example~\cite{nitzberg1993filtering,mm98,shen2002mathematical}) work only if this mask can be identified, inpainting the missed area only. Under these assumptions, it is rather natural to define the confidence weight as
\begin{equation}
\label{eq:confidence-mask}
    \confidence_{\text{mask}}(x)
    =
    \left\{
    \begin{aligned}
    &1 \text{ if } x\in \Omega\setminus \Mask \\
    &0 \text{ if } x\in \Mask
    \end{aligned}
    \right.,
\end{equation}
in order to localize the reconstruction of the network only in the corrupted regions. However, in many inpainting problems, we may not have a priori information on where such regions are located and how the data are corrupted. In this case, a more agnostic discrepancy weight should be preferred. Our second candidate is simply $\confidence=1$.

\subsection{Mapping $\Tdens$ into $\Img$}
\label{sec:map}
The transformation $\ImgTd$ has the scope of connecting two objects that may be different, the conductivity $\Tdens$ in the \BT{} problem, which is an model-based variable, into an image $\Img$, which is instead an observation-based variable. This transformation must send nonnegative conductivity $\Tdens$ to nonnegative images $\Img$. Furthermore, it must be differentiable for all $\Tdens\geq 0$. This is necessary to compute the gradient of the functional $\Inpaint$ in~\cref{eq:inpainting-niot-intro}. We do not require the map to be invertible, but it may be useful to infer information on $\Tdens$ from observations.

Any map $\ImgTd$ will necessarily include some multiplicative factor $\Scaling>0$ such that $\ImgTd(\Tdens)$ and $\Obs{\Img}$ within the discrepancy term in~\cref{eq:inpainting-niot-intro} have the same order of magnitude. For example, when $\lambda=0$, if the terms $\Source$ and $\Sink$ are both multiplied by a constant $c$, it is easy to see that the optimal conductivity $\OptTdensH$ will be multiplied by a factor $c^{\frac{2}{\gamma+1}}$. This multiplicative factor can also be deduced from the physical information encoded in the images. 

Moreover, the map $\ImgTd$ should be designed such that $\ImgTd(\Tdens)$ falls within the range of values in which data $\Obs{\Img}$ are stored. For example, in the case of standard grayscale images, it must be ensured that $\ImgTd(\Tdens)\in[0,1]$. 

In this paper, two maps $\ImgTd$ are considered. The first map is the identity map, that is, $\ImgTd(\Tdens)=\Scaling\Tdens$. This map is extremely simple and is suitable for inpainting problems where the support of $\Tdens$ is expected to be similar to the support of $\Obs{\Img}$, and the observed data $\Obs{\Img}$ are proportional to the conductivity $\Tdens$.
Its main advantage is that it is linear and the gradient of the functional $\Inpaint$ in~\cref{eq:inpainting-niot-intro} can be computed without any additional cost.

The second map examined in this paper is designed to handle the case where some fluid is passing through the observed network. Under mild assumptions, this flow can be modeled using the Hagen-Poiseuille approximation, which states that the velocity of a fluid flowing through a pipe is proportional to the gradient of a pressure multiplied by a conductivity factor. For tubular pipes, this conductivity factor is expressed as a power law of the radius $r>0$ of the channel, given by
\begin{equation}
\label{eq:pouseille-law}
  \PouTdens = \PouseilleConstant r^p  
\end{equation}
where $\PouseilleConstant$ is a constant factor, while $p=3$ in $\REAL^2$, while $p=4$ in $\REAL^3$. 
Note that in the \BT{} formulation in~\cref{eq:poisson} we implicitly assumed the Hagen-Poiseuille law to hold. 

A map $\ImgTd$ adapted to this type of networks should transform a conductivity $\Tdens$ concentrated on a union of 1d curves into an image $\Img$ whose support satisfies the relation in~\cref{eq:pouseille-law}. This implies mapping 1d channels with higher conductivities into tubular neighborhood with wider support, making the map necessarily nonlinear. 

In the variational setting of~\cref{eq:inpainting-niot-intro}, we cannot take a simple tubular expansion of a 1d channel, possibly tuned according to $\Tdens$, since this transformation is not differentiable. Moreover, it would be complicated to define this transformation, where the support of $\Tdens$ is not exactly supported on a 1 pixel wide curve. 

We can define a map having the desired characteristics using the porous media (\PM{}) equation~\cite{Vazquez2007,fasano2006problems}, which can be formulated as follows. Given $m\geq1$, $T>0$, and a nonnegative initial condition $\pmu_0$, we seek a function $\pmu:[0,T]\times \Omega\to [0,+\infty]$ that satisfies the following equation
\begin{equation}
    \label{eq:porous-media}
    \left\{
    \begin{aligned}
    &\Dt{\pmu}-\Div(\pmu^{m-1}\Grad \pmu)=0,\\
    &\pmu(0)=\pmu_0,\\
    &\pmu^{m-1}\Grad \pmu \cdot n_{\partial \Domain}=0.
    \end{aligned}
    \right.
\end{equation}
Denoting by $\pmu(t,m,\pmu_0)$ the solution of the PDE~\cref{eq:porous-media} at time $t$ with initial condition $\pmu_0$, we can define a map $\ImgTd$ (hereafter referred as the \PM{} map) as follows
\begin{equation}
  \label{eq:pm-map}
  \ImgTd(\Tdens):=\Scaling \pmu(t^*,m,\Tdens),
\end{equation}
with $t^*$ being a suitable time. In~\cref{sec:pm} we explain how to deduce the exponent $m$ and $t^*$ from the exponent $p$ and the constant $\PouseilleConstant$ in~\cref{eq:pouseille-law}, with the formulas given in~\cref{eq:mt}.

Note also that for $m=1$, the equation in~\cref{eq:porous-media} reduces to the heat equation, and thus the maps would essential act as a Gaussian filter. This map is suitable for the inpainting problem, with the advantage of being easier to compute. However, it has the disadvantage of infinite propagation of the support and, being linear, does not provide the ``higher'' to ``wider'' mapping needed to approximate the law in~\cref{eq:pouseille-law}. For this reason, we did not consider it in this paper.

\subsection{Discretization and optimization}
\label{sec:algorithm}
The spatial discretization of~\cref{eq:inpainting-niot-intro} is performed using a finite-volume/finite-difference scheme based on a Cartesian grid $\TriangH$ of the domain $\Domain$. This discretization is particularly convenient in inpainting problems, where we mostly deal with images stored in pixels/voxels. We now look for the pair $(\PotH,\TdensH)\in\DGzero(\TriangH)\times\DGzero(\TriangH)$, where $\DGzero(\TriangH)$ is the space of piecewise functions on $\TriangH$.

The gradient of the discretized functional $\Inpaint$ in~\cref{eq:inpainting-niot-intro} with respect to $\TdensH$ is given by 
\begin{equation}
  \label{eq:inpainting-gradient-discrete}
  \Grad_{\TdensH}{\Inpaint}(\TdensH)=
  \frac{
  -|\Grad \Pot(\TdensH)|^2+\TdensH^{\gamma-1}
  }{2}
  + \lambda(\ImgTd(\TdensH)-\Obs{\Img}) \confidence(x) \ImgTd'(\TdensH)
\end{equation}
where $\PotH\Of{\TdensH}\in \DGzero(\TriangH)$ is the numerical solution of the elliptic equation in~\cref{eq:poisson}. In this process, a map from the cell-centered values of $\TdensH$ to the faces is needed. We used the harmonic average for this purpose to ensure that the flux is properly conserved.

When the \PM{} map is used, there is an additional cost to compute the term involving the Jacobian $\ImgTd'(\TdensH)$ coming from the discrepancy term in~\cref{eq:inpainting-gradient-discrete}. To handle this, we employ adjoint-based methods \cite{M2AN_1986__20_3_371_0}. In our case, these methods require solving the \PM{} equation in~\cref{eq:porous-media} to obtain the term $\ImgTd(\TdensH)$. We use the backward Euler method combined with the Newton method for this purpose. Additionally, we need to solve a linear system with the Jacobian of the \PM{} map $\ImgTd'(\TdensH)$ as the matrix. To improve accuracy, we divide the time $t^*$ in~\cref{eq:pm-map} into five subintervals, gradually increasing the initial time step. Consequently, in the calculation of the gradient in~\cref{eq:inpainting-gradient-discrete}, we solve five equations of porous media and five linear systems associated with the adjoint-based method.

The optimization part follows the approach used in~\cite{Facca2021}. 
Starting from an initial data $\TdensH^0=\TdensIni$, we update $\TdensH$ as follows
 \begin{align}
  \label{eq:niot-dyn}
    & \TdensH^{k+1} = \TdensH^{k} - \Delta t^{k}(\TdensH^k)^{\frac{2}{\gamma}} \Grad_{\TdensH}{\Inpaint}(\TdensH^k),\quad k=0,\ldots,k_{\text{max}},
\end{align}
with $\Delta t^{k}>0$ being a sequence of time step adaptively to ensure the stability of the scheme.  The direction of update in~\cref{eq:niot-dyn} is (minus) the gradient of the functional $\Inpaint$ scaled by the factor $\TdensH^{\frac{2}{\gamma}}$. This scaling procedure ensures that $\TdensH$ remains positive and mimics a scheme known in the optimization community as the mirror descent approach~\cite{nemirovsky1983wiley}. We continue updating $\TdensH$ until the $L^1$-norm of the discrete right-hand side of~\cref{eq:niot-dyn} is below $10^{-5}$.

The method we have introduced falls within the class of first-order methods, as it relies only on the gradient of $\Inpaint$. Approaches that incorporate second-order derivatives, such as a Newton-based scheme, appear to be not suitable for our problem. This is mainly due to the nonconvex nature of the functional $\Inpaint$. Furthermore, these methods require specific linear algebra solutions, which are challenging to modify to fit our circumstances, similar to those outlined in~\cite{FaccaBenzi2021}.

\subsection{Software and implementation details}
\label{sec:implementation}
The NIOT algorithm is implemented using Firedrake~\cite{rathgeber2016firedrake,FiredrakeUserManual}, a Python finite element library based on the Unified Form Language (UFL) introduced in~\cite{alnaes2014unified}. Using these libraries, we can write our problem directly from the functional $\Inpaint$ in~\cref{eq:lyap} and then use automatic differentiation to compute the gradient of the functional. The adjoint-based computations are performed using the dolfin-adjoint library~\cite{Mitusch2019}. Firedrake is also built on top of the PETSc library~\cite{petsc-user-ref,petsc-efficient}, which provides a wide range of efficient linear and nonlinear solvers.

We solve linear systems associated to elliptic PDE~\cref{eq:poisson}, and to the adjoint-based method using iterative methods preconditioned with the Hypre BoomerAMG algebraic multigrid solver described in~\cite{falgout2006design}. 

The source code used to produce the results in this paper is available at \url{https://github.com/enricofacca/niot}, together with the data and scripts to reproduce the results. A persistent identifier of the code used in this article is also available in the Zenodo repository \url{https://doi.org/10.5281/zenodo.11241387}.

\section{Numerical experiments}
\label{sec:experiments}
We now present a set of experiments designed for dual purpose. Our first objective is to address reconstruction problems in which data degradation is the result of data loss that causes network disconnection, and/or the existence of artifacts, as exemplified in~\cref{fig:network-inpainting-obs}. Our second goal is to investigate the interaction between the parameters and functions of the NIOT model, summarized in~\cref{tab:parameters}. With this goal in mind, we examine, in different setups, two test cases: a simple hand-drawn Y-shaped network and a more complex network extracted from a real-world geometry. 

\begin{table}[htbp]
  \centering
  \begin{tabular}{|c|c|c|}
    \hline
    \textbf{Category} & \textbf{Parameter} & \textbf{Description} \\
    \hline
    \multirow{2}{*}{Modeling} & $\Source$, $\Sink$ & Mass distribution densities in the \BT{} problem \\
    \cline{2-3}
    & $\gamma$ & Scalar determining branching angles \\
    \hline
    \multirow{3}{*}{Fitting} & $\ImgTd$, $\Scaling$ & Map from $\Tdens$ to $\Img$ \\
    \cline{2-3}
    & $\lambda$ & Scalar measuring global trust in data\\
    \cline{2-3} 
    & $\confidence$ & Function measuring local trust in data\\
    \hline
    Optimization &$\Tdens_0$ & Initial data in~\cref{eq:niot-dyn} \\
    \hline 
  \end{tabular}
  \caption{Parameters and functions controlling the NIOT model.}
  \label{tab:parameters}
\end{table}

All images reported in the following experiments should be seen in their high-resolution colored version, which is available directly within the electronic version of the paper. These images were generated using the PyVista library~\cite{sullivan2019pyvista}.

\subsection{Y-shaped hand drawing}
\label{sec:y-net}
The first test case presented in this paper is the reconstruction of the hand-drawn Y-shaped network in~\cref{fig:network-inpainting-true} after being corrupted using the mask made up of three rectangles shown in~\cref{fig:network-inpainting-obs}. In this test case, the artifacts shown in~\cref{fig:network-inpainting-obs} are removed.

First, we fix the resolution of the mesh used in our inpainting problem. We choose a mesh $\Triang_{h/4}$ with $208\times208$ pixels, used in the third column in~\cref{fig:dmk}. For this grid, the smallest branch is at least 4 pixels wide, while at coarser resolutions it is hard to see the effects of using the $\PM{}$ map used in the next experiments.

We also need to fix the exponent $\gamma$ that influences the branching angle of the optimal network. Looking at the networks obtained in~\cref{fig:dmk} we can see how the exponent $\gamma=0.8$ does not provide a Y-shaped network, which is instead obtained for $\gamma=0.2$ and $\gamma=0.5$. 
We decided to fix $\gamma=0.5$ for all experiments, because the network obtained for $\gamma=0.2$ is not that different from the one obtained for $\gamma=0.5$, while solving the \BT{} problem with smaller values of $\gamma$ is typically computationally more intense. This is because shorter time steps are required due to the factor $\Tdens^{\frac{2}{\gamma}}$ in~\cref{eq:niot-dyn}, and therefore more time steps to reach convergence.

\subsubsection{Images with lost data (case $\confidence=\confidence_{\text{mask}}$)}
\label{sec:netnetwork_confMASK}
We now turn our attention to the discrepancy term in~\cref{eq:inpainting-niot-intro}, which defines the fitting component of the NIOT model.

We consider two maps $\ImgTd$ described in~\cref{sec:map}, the identity map and the \PM{} map with $m=2$ and $t^*=1.0e-3$. The exponent $m$ is the result of the formula in~\cref{eq:mt} with $p=3$, while the time $t^*$ was experimentally found to guarantee that the support of $\Obs{\Img}$ and the reconstructed image $\ImgTd(\Tdens)$ have approximately the same supports.

The identity map is scaled by the factor $\Scaling=1.0e1$, while the \PM{} map is scaled by $\Scaling=5.0e1$. We choose these values by looking at the optimal $\OptTdensH$, reported in the rightmost column in~\cref{fig:dmk}, which correspond to the solution of~\cref{eq:inpainting-niot-intro} with $\lambda=0$. We tune the scaling parameter $\Scaling$ to have the same order of magnitude of the maximum values of the reconstructed image $\ImgTd(\Tdens)$ and the data term, which is a binary image. This calibration based on the maximum values presents some limitations. In fact, in full generality, we expect to have some mismatch between the observed data $\Obs{\Img}$ and the $\ImgTd(\Tdens)$ in the sub-branches of the network, where the conductivity is lower. Nonetheless, this method is rather straightforward and prevents the observed and reconstructed data from having values that vary by multiple orders of magnitude.

We fix the confidence term $\confidence=\confmask$ given in~\cref{eq:confidence-mask} and the initial data $\Tdens_0=1$ (the influence of the second choice is discussed in~\cref{sec:tdensini-obs}). In the upper panels of~\cref{fig:netnetwork-all} we report the result obtained for the different values of $\lambda$ and the two maps described above.

{
  \def \fraction {0.22}
  \def \extrah {0.01em}
  \def \labelh {0.5em}
  \begin{figure}
    \centering
  \begin{tabular}{|@{}c@{}|@{\hspace{0.25em}}c@{\hspace{\labelh}}|@{\hspace{\extrah}}c@{\hspace{\extrah}}|@{\hspace{\extrah}}c@{\hspace{\extrah}}|@{\hspace{\extrah}}c@{\hspace{\extrah}}|@{\hspace{\extrah}}c@{\hspace{\extrah}}|}
      \hline
      && $\lambda=1e-2$ & $\lambda=5e-2$  &  $\lambda=1e-1$ & $\lambda=1e0$ 
      \\
      \hline
      \multirow{2}{*}{
      \raisebox{-.5\normalbaselineskip}[0pt][0pt]{\rotatebox[origin=c]{90}{$\confidence=\confmask$, $\Tdens_0=1$}}
      }
      &
      \raisebox{-.5\normalbaselineskip}[0pt][0pt]{\rotatebox[origin=c]{90}{$\ImgTd(\Tdens)=1e1\Tdens$}}
      &
        \tmpframe{\adjustbox{valign=m,vspace=0pt}{\includegraphics[trim={0.9cm 2.7cm 0.9cm 2.7cm},clip,width=\fraction\columnwidth,valign=t]{{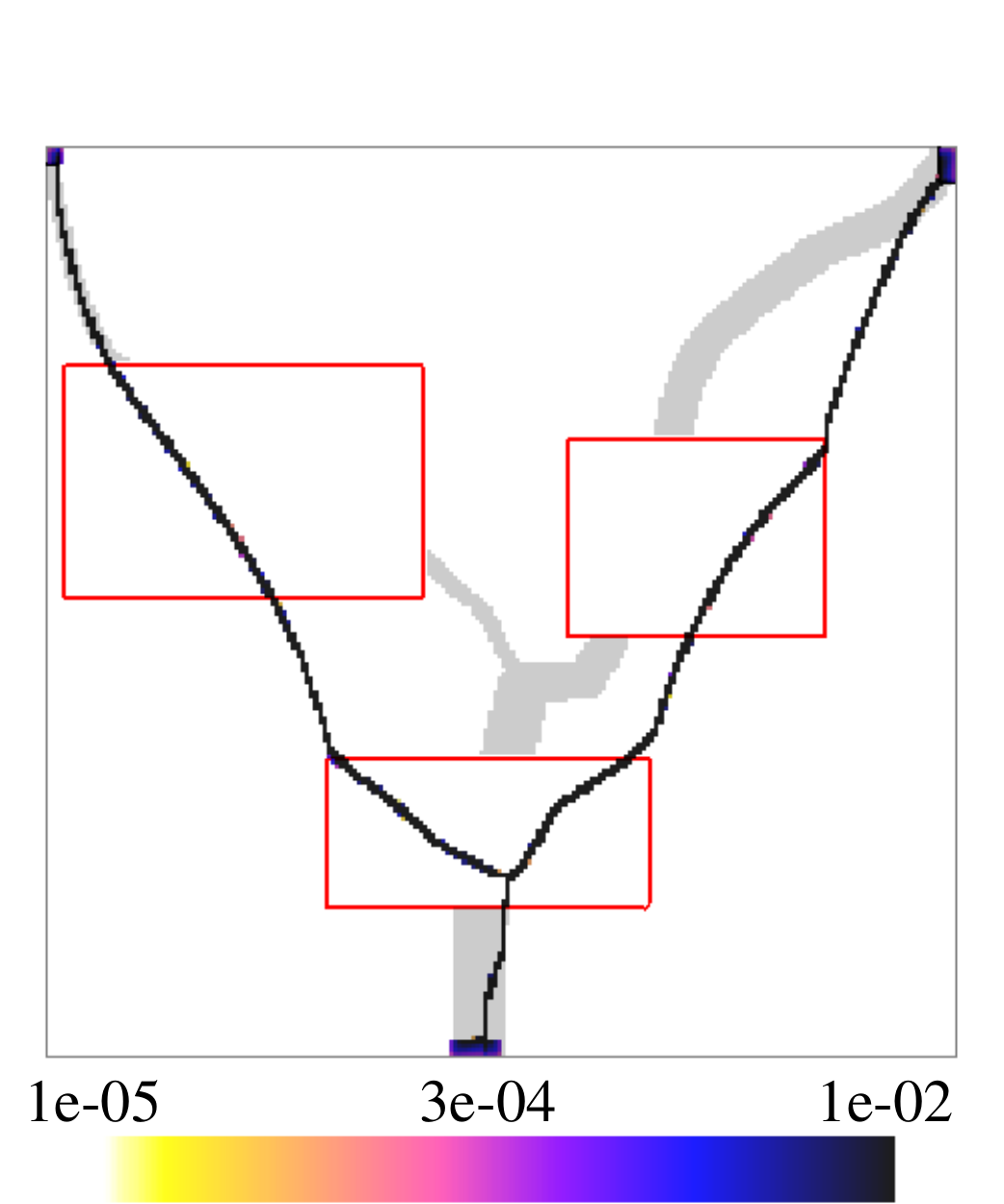}}}}
        &
        \tmpframe{\adjustbox{valign=m,vspace=0pt}{\includegraphics[trim={0.9cm 2.7cm 0.9cm 2.7cm},clip,width=\fraction\columnwidth,valign=t]{{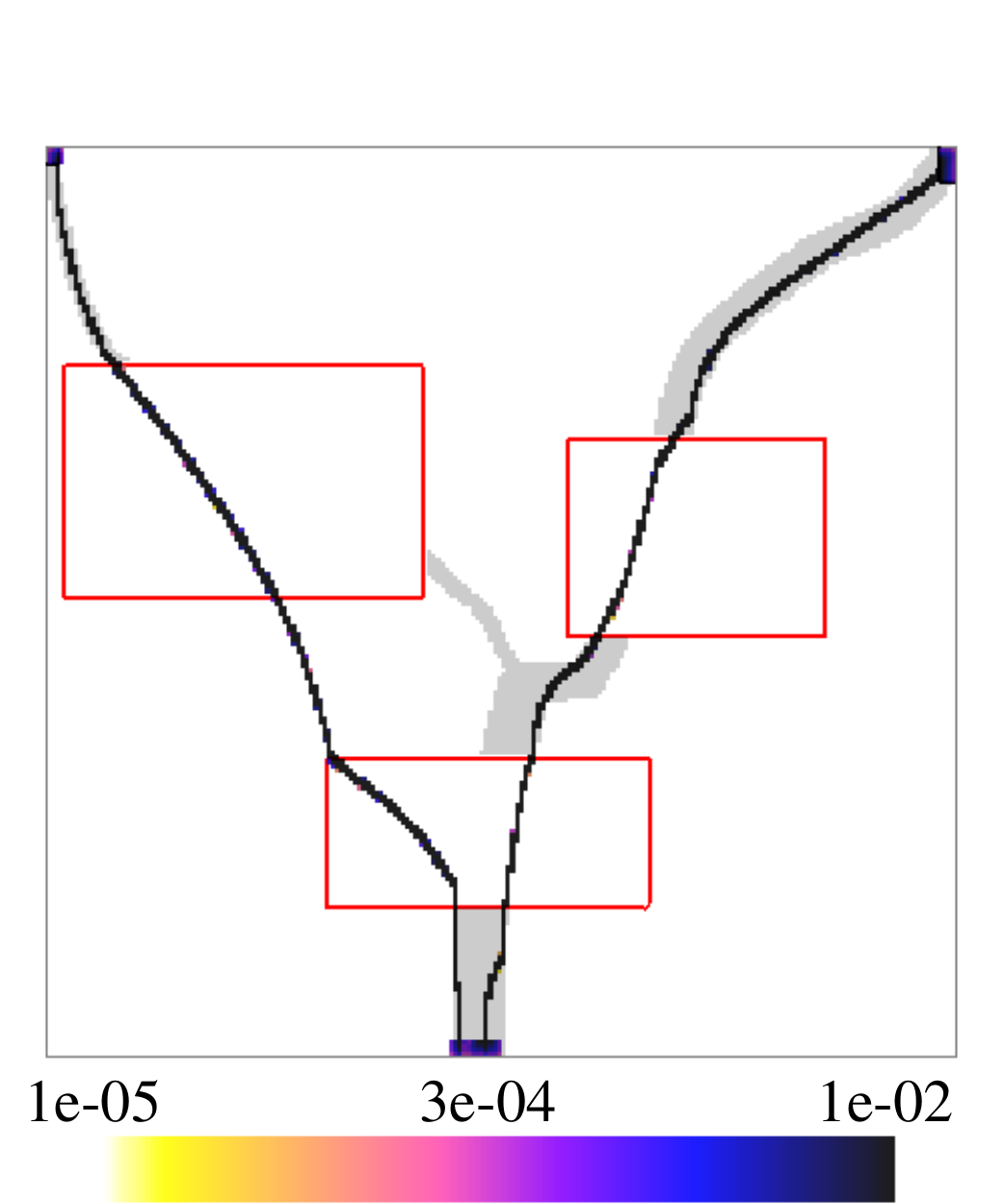}}}}
        &
        \tmpframe{\adjustbox{valign=m,vspace=0pt}{\includegraphics[trim={0.9cm 2.7cm 0.9cm 2.7cm},clip,width=\fraction\columnwidth,valign=t]{{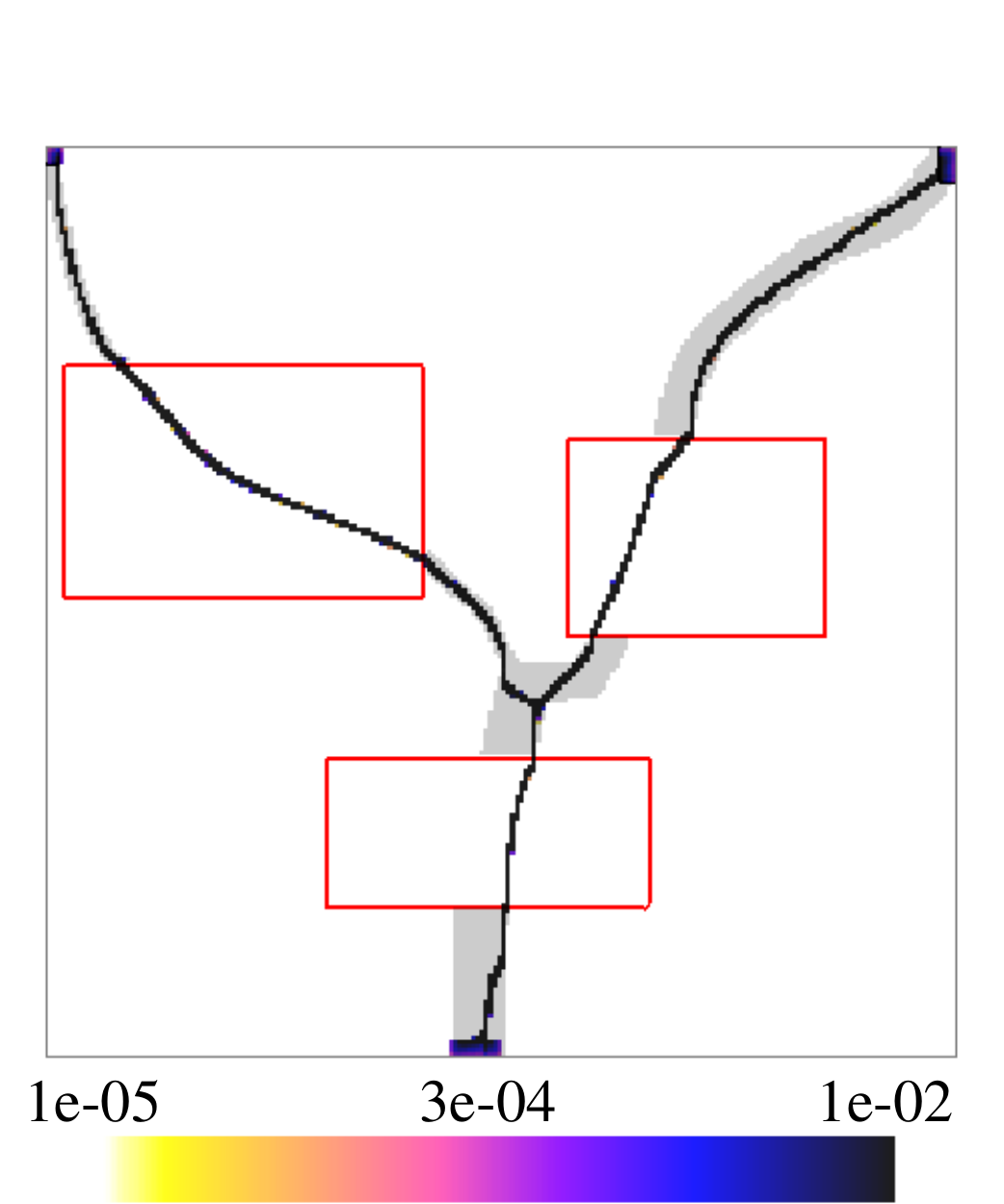}}}}
        &
        \tmpframe{\adjustbox{valign=m,vspace=0pt}{\includegraphics[trim={0.9cm 2.7cm 0.9cm 2.7cm},clip,width=\fraction\columnwidth,valign=t]{{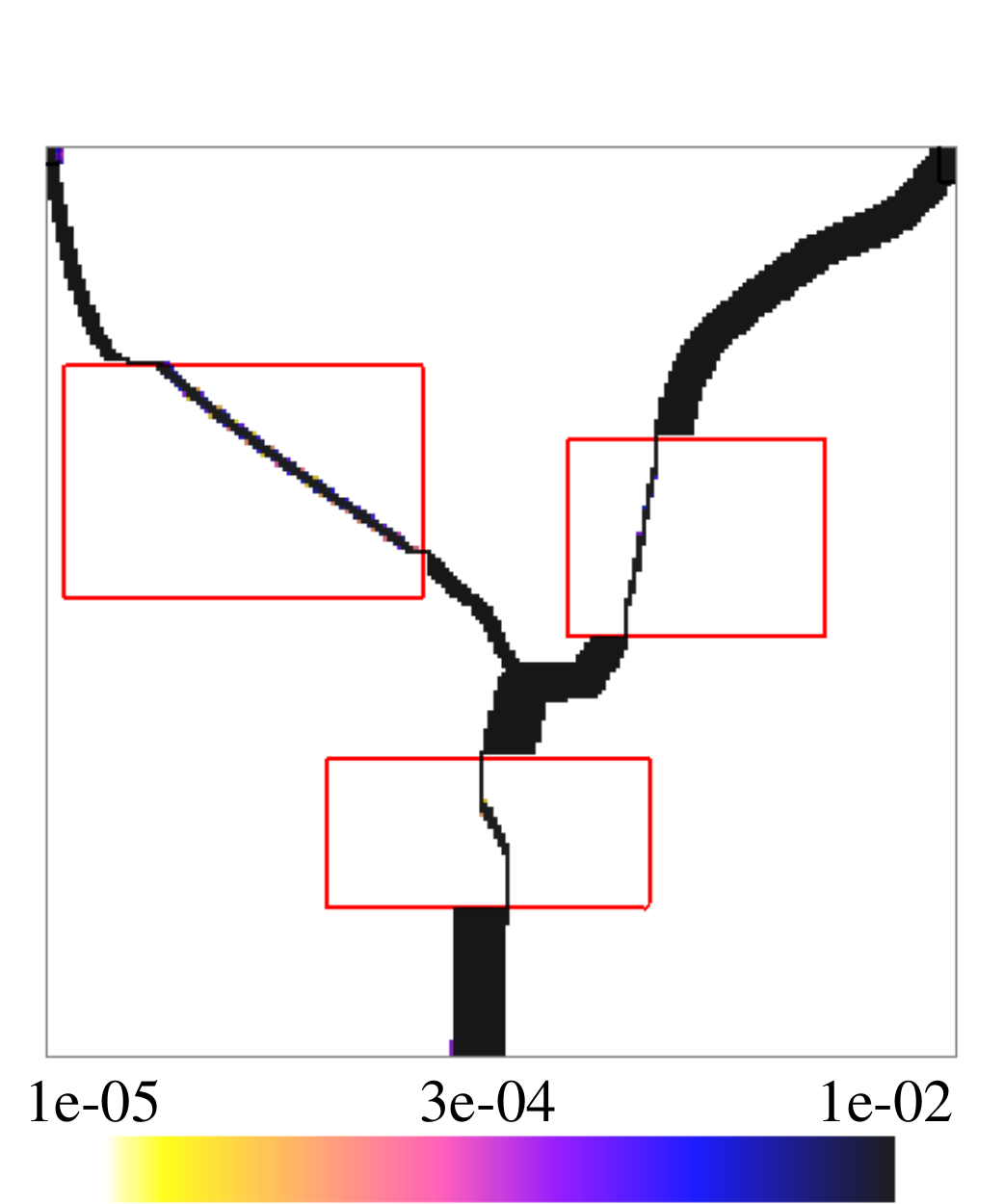}}}}
      \\
      \cline{2-5}
      &
      \raisebox{-.5\normalbaselineskip}[0pt][0pt]{\rotatebox[origin=c]{90}{\quad $\ImgTd(\Tdens)=50\PM(\Tdens)$}}
      &
      \tmpframe{\adjustbox{valign=m,vspace=0pt}{\includegraphics[trim={0.9cm 2.7cm 0.9cm 2.7cm},clip,width=\fraction\columnwidth,valign=t]{{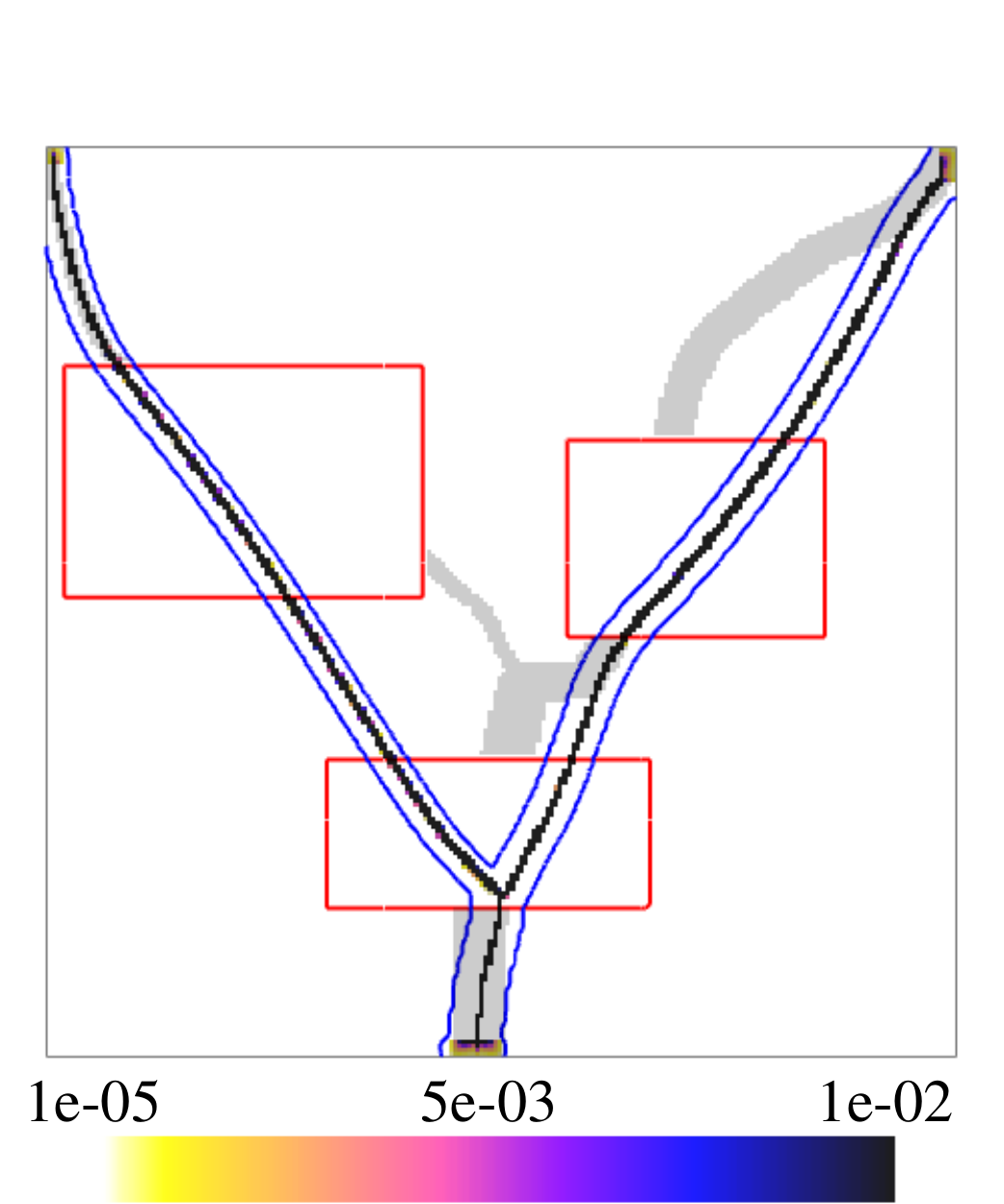}}}}
      &
      \tmpframe{\adjustbox{valign=m,vspace=0pt}{\includegraphics[trim={0.9cm 2.7cm 0.9cm 2.7cm},clip,width=\fraction\columnwidth,valign=t]{{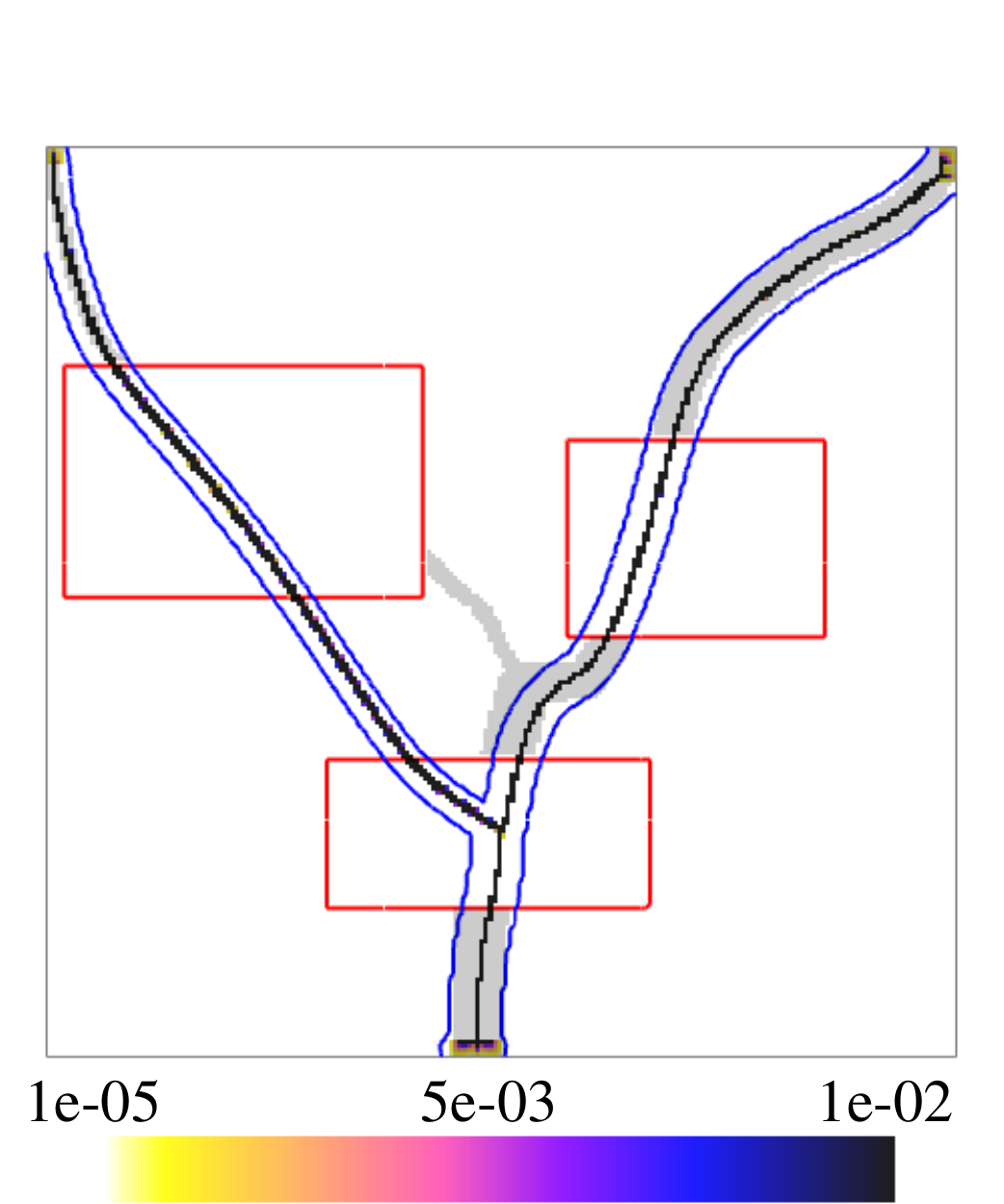}}}}
      &
      \tmpframe{\adjustbox{valign=m,vspace=0pt}{\includegraphics[trim={0.9cm 2.7cm 0.9cm 2.7cm},clip,width=\fraction\columnwidth,valign=t]{{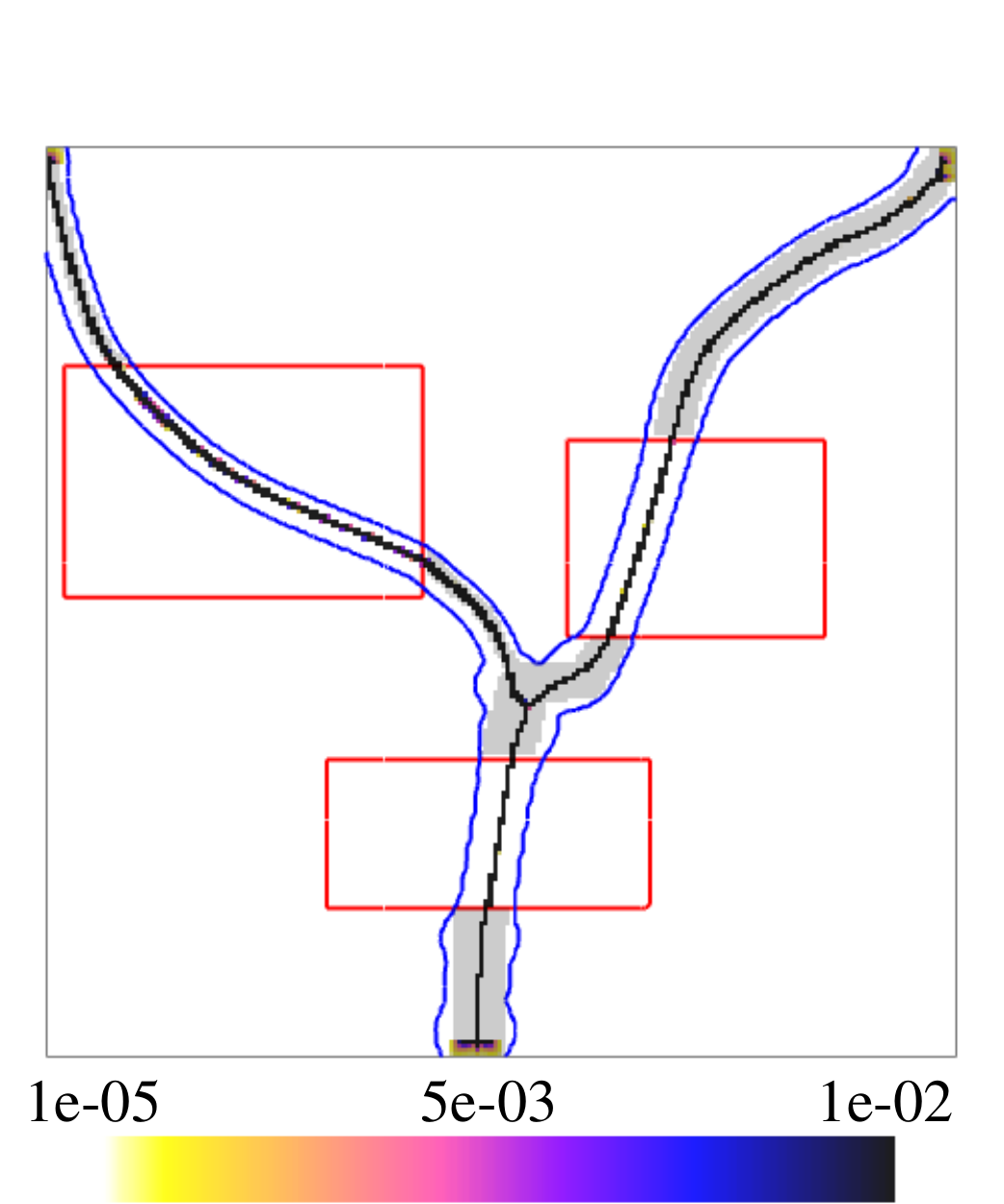}}}}
      &
      \tmpframe{\adjustbox{valign=m,vspace=0pt}{\includegraphics[trim={0.9cm 2.7cm 0.9cm 2.7cm},clip,width=\fraction\columnwidth,valign=t]{{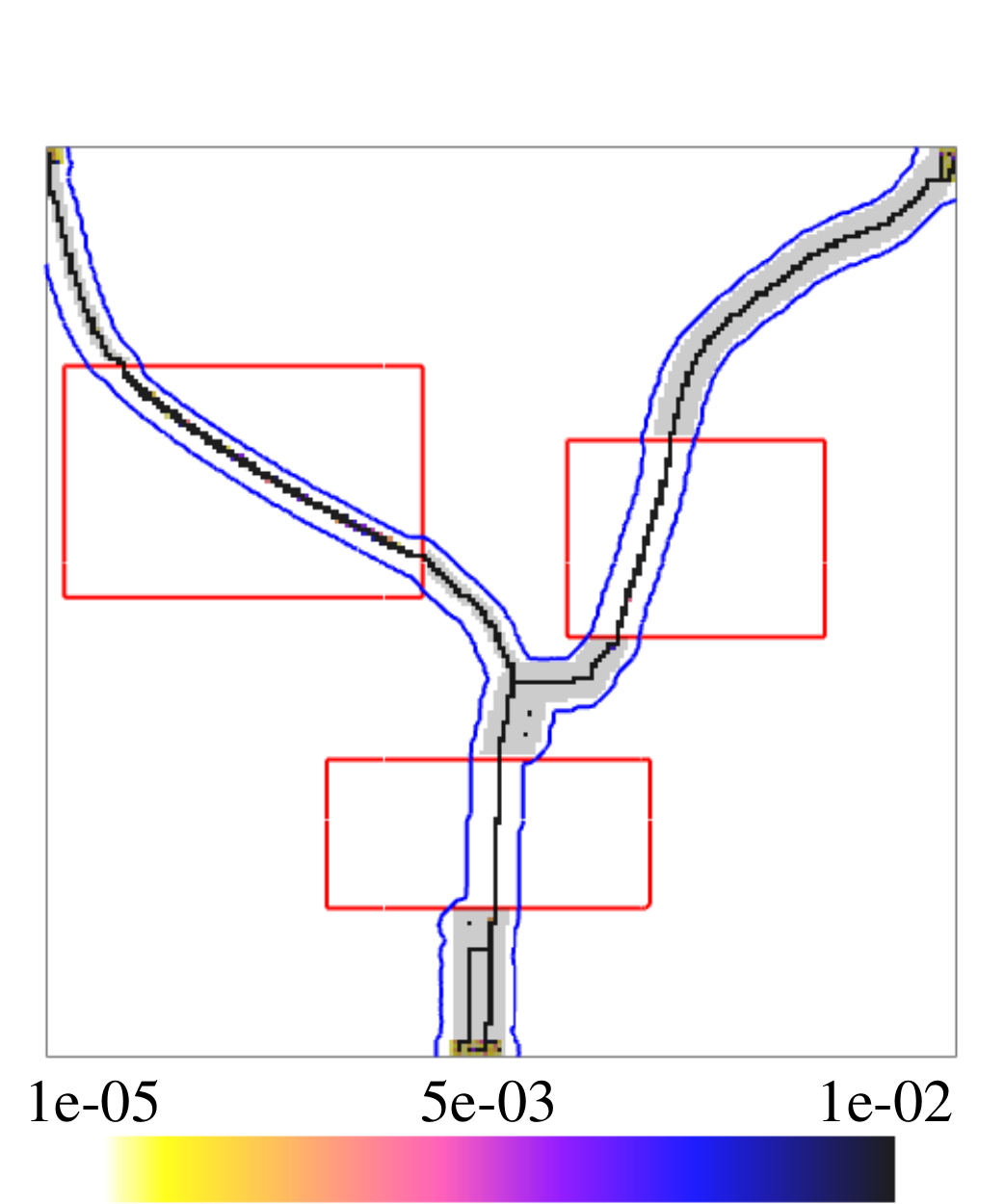}}}}
      \\
      \hline
      \multirow{2}{*}{
      \raisebox{-.5\normalbaselineskip}[0pt][0pt]{\rotatebox[origin=c]{90}{$\confidence=1$, $\Tdens_0=1$}}
      }
      &
      \raisebox{-.5\normalbaselineskip}[0pt][0pt]{\rotatebox[origin=c]{90}{\quad $\ImgTd(\Tdens)=10\Tdens$}}
      &
      \tmpframe{\adjustbox{valign=m,vspace=0pt}{\includegraphics[trim={0.9cm 2.7cm 0.9cm 2.7cm},clip,width=\fraction\columnwidth,valign=t]{{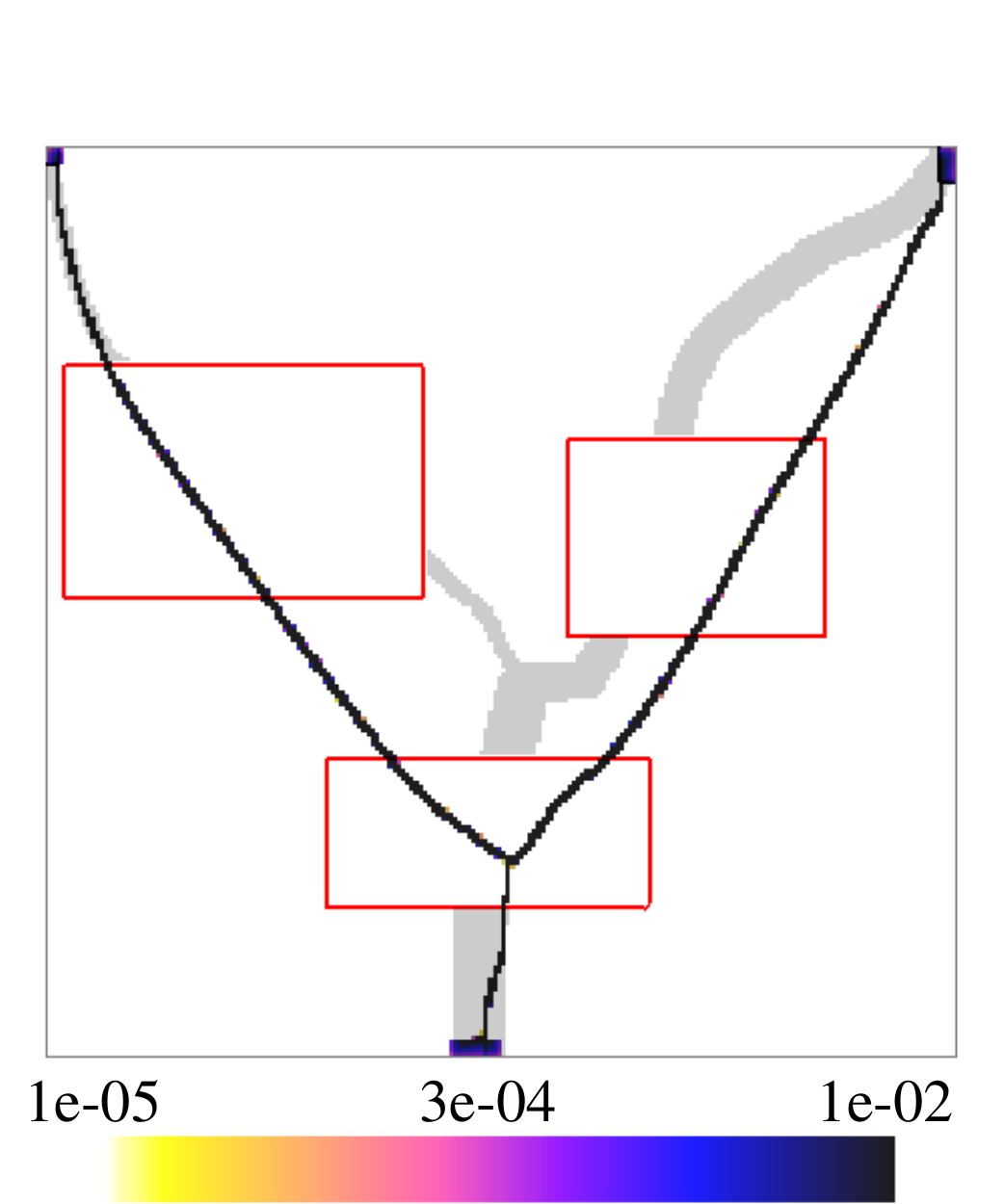}}}}
      &
      \tmpframe{\adjustbox{valign=m,vspace=0pt}{\includegraphics[trim={0.9cm 2.7cm 0.9cm 2.7cm},clip,width=\fraction\columnwidth,valign=t]{{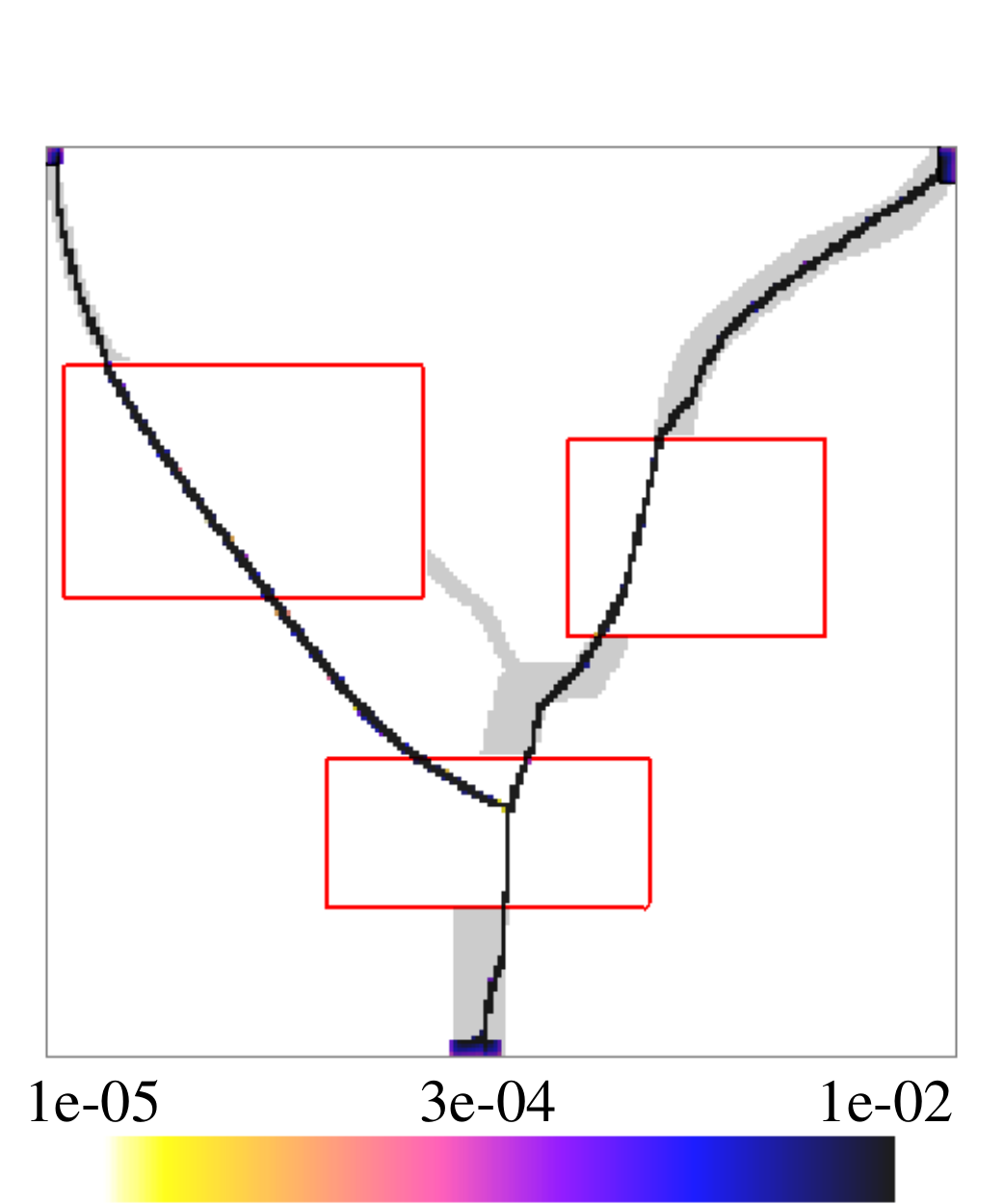}}}}
      &
      \tmpframe{\adjustbox{valign=m,vspace=0pt}{\includegraphics[trim={0.9cm 2.7cm 0.9cm 2.7cm},clip,width=\fraction\columnwidth,valign=t]{{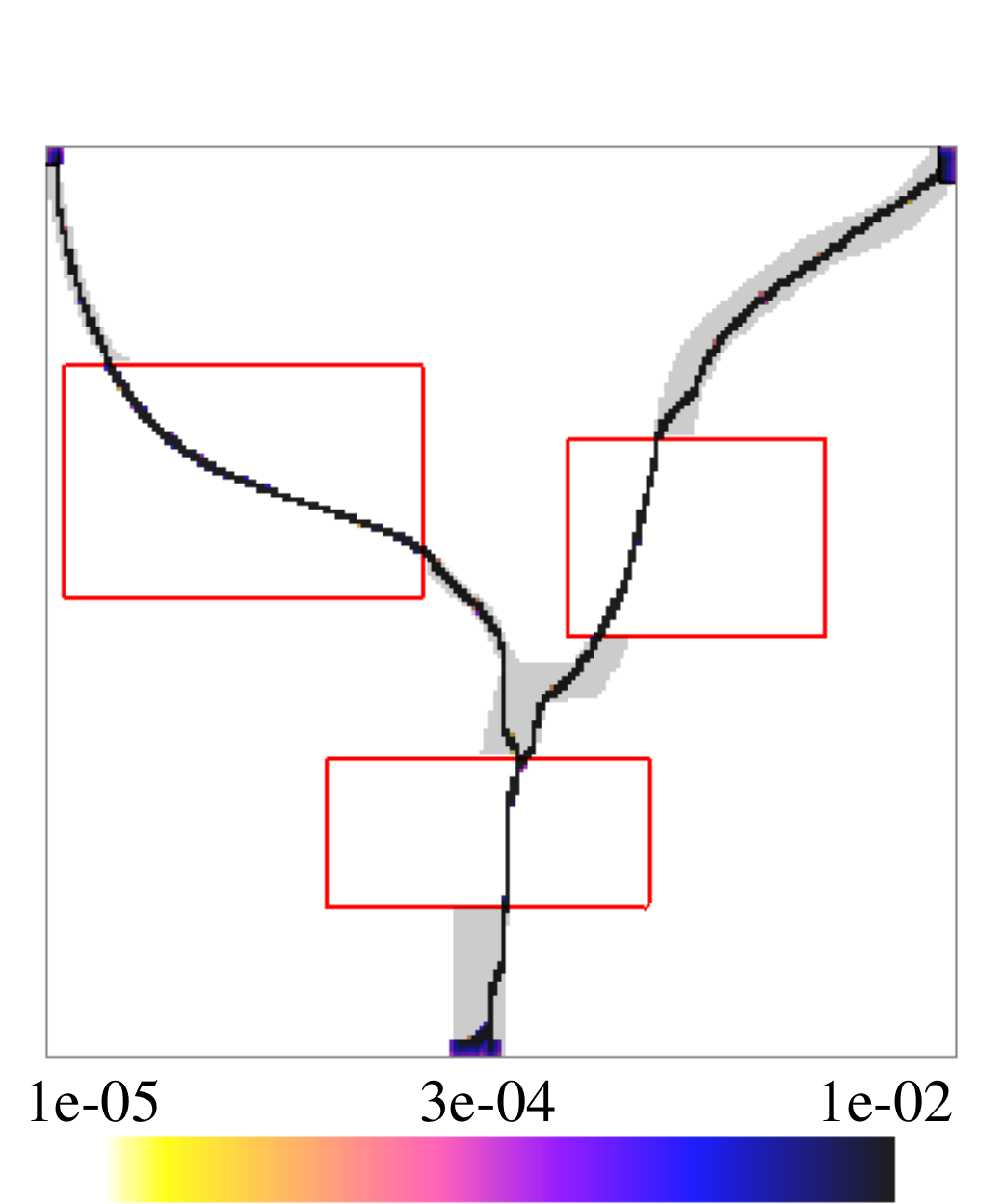}}}}
      &
      \tmpframe{\adjustbox{valign=m,vspace=0pt}{\includegraphics[trim={0.9cm 2.7cm 0.9cm 2.7cm},clip,width=\fraction\columnwidth,valign=t]{{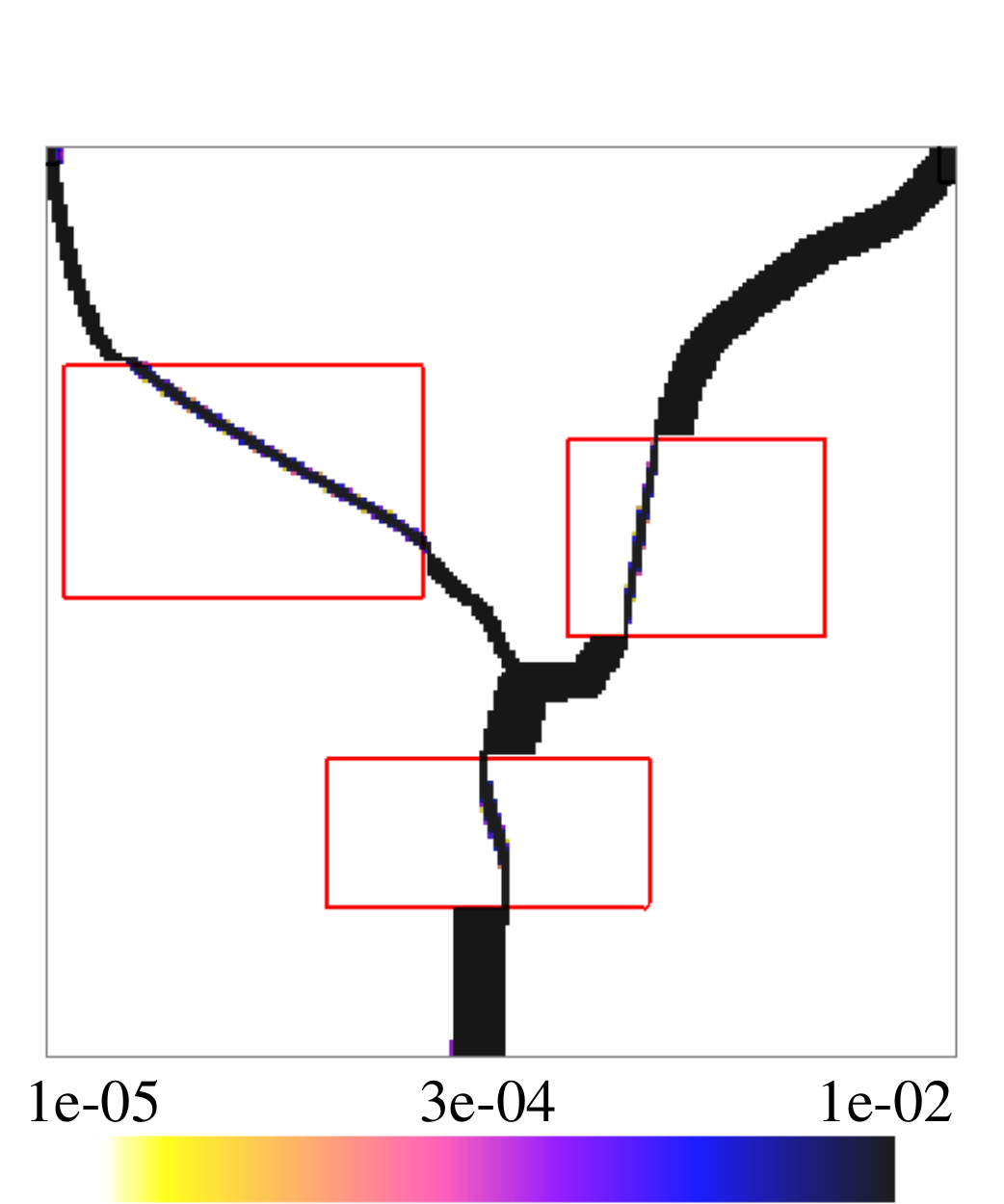}}}}
      \\
      \cline{2-5}
      &
      \raisebox{-.5\normalbaselineskip}[0pt][0pt]{\rotatebox[origin=c]{90}{\quad $\ImgTd(\Tdens)=50\PM(\Tdens)$}}
      &
      \tmpframe{\adjustbox{valign=m,vspace=0pt}{\includegraphics[trim={0.9cm 2.7cm 0.9cm 2.7cm},clip,width=\fraction\columnwidth,valign=t]{{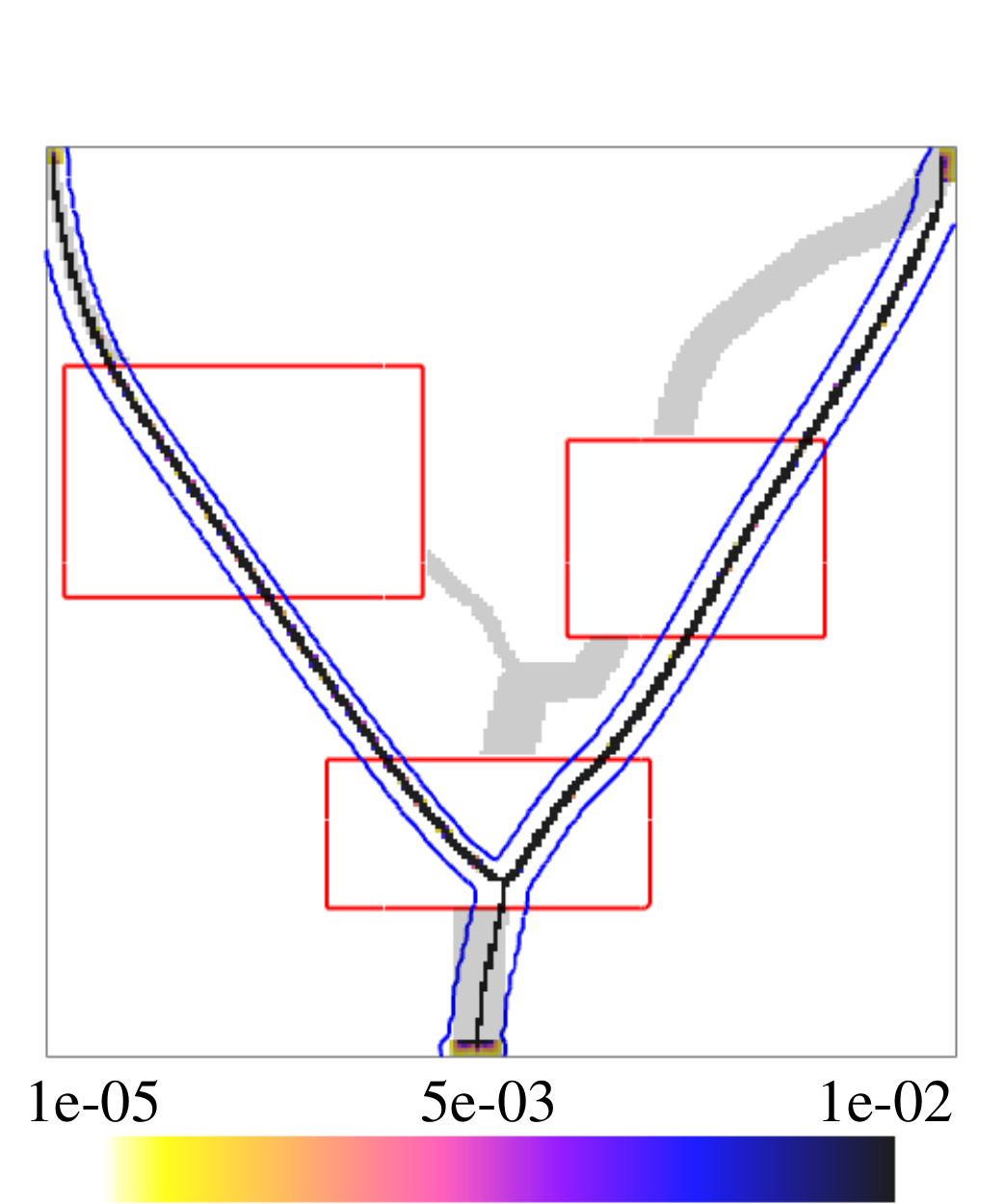}}}}
      &
      \tmpframe{\adjustbox{valign=m,vspace=0pt}{\includegraphics[trim={0.9cm 2.7cm 0.9cm 2.7cm},clip,width=\fraction\columnwidth,valign=t]{{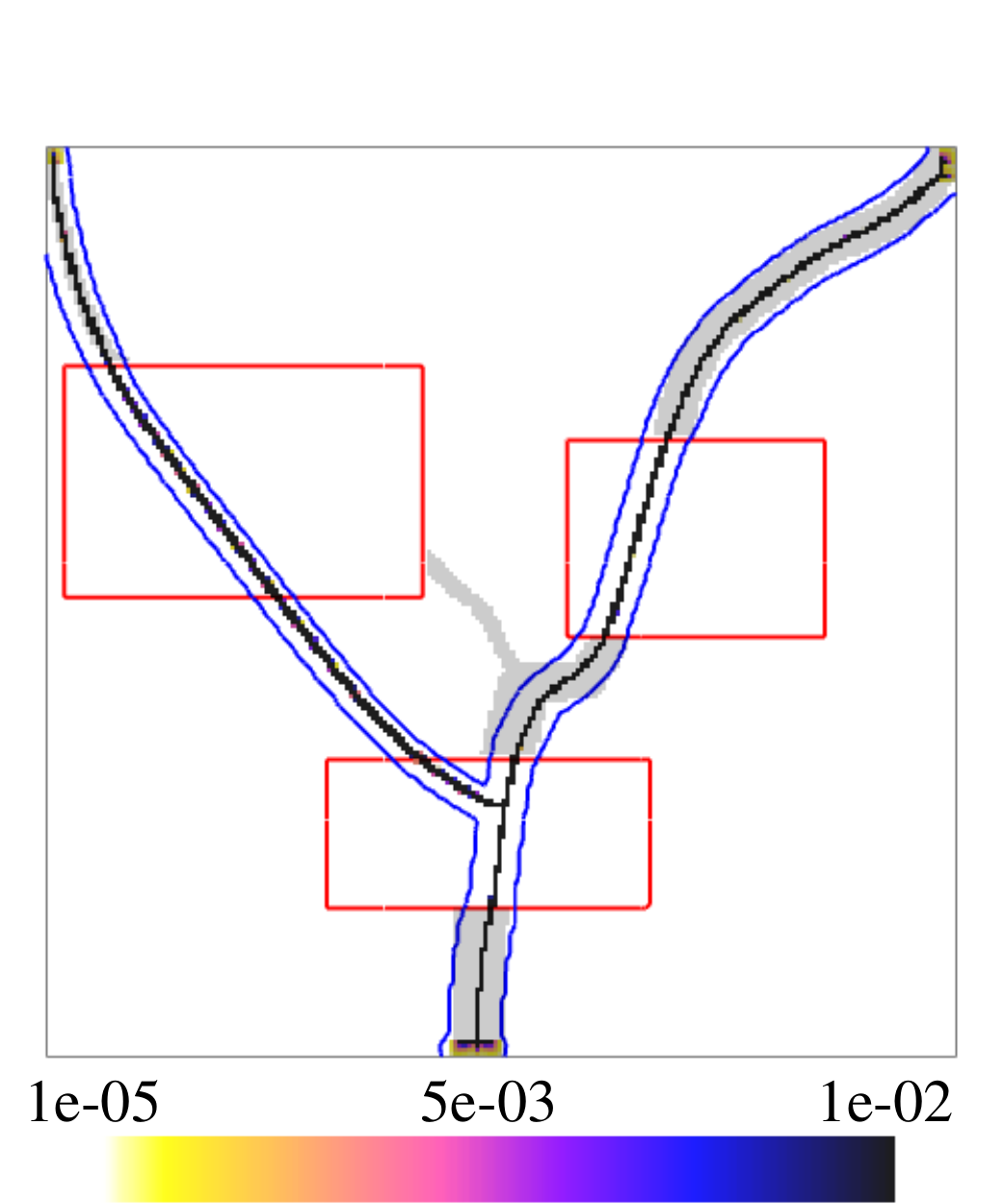}}}}
      &
      \tmpframe{\adjustbox{valign=m,vspace=0pt}{\includegraphics[trim={0.9cm 2.7cm 0.9cm 2.7cm},clip,width=\fraction\columnwidth,valign=t]{{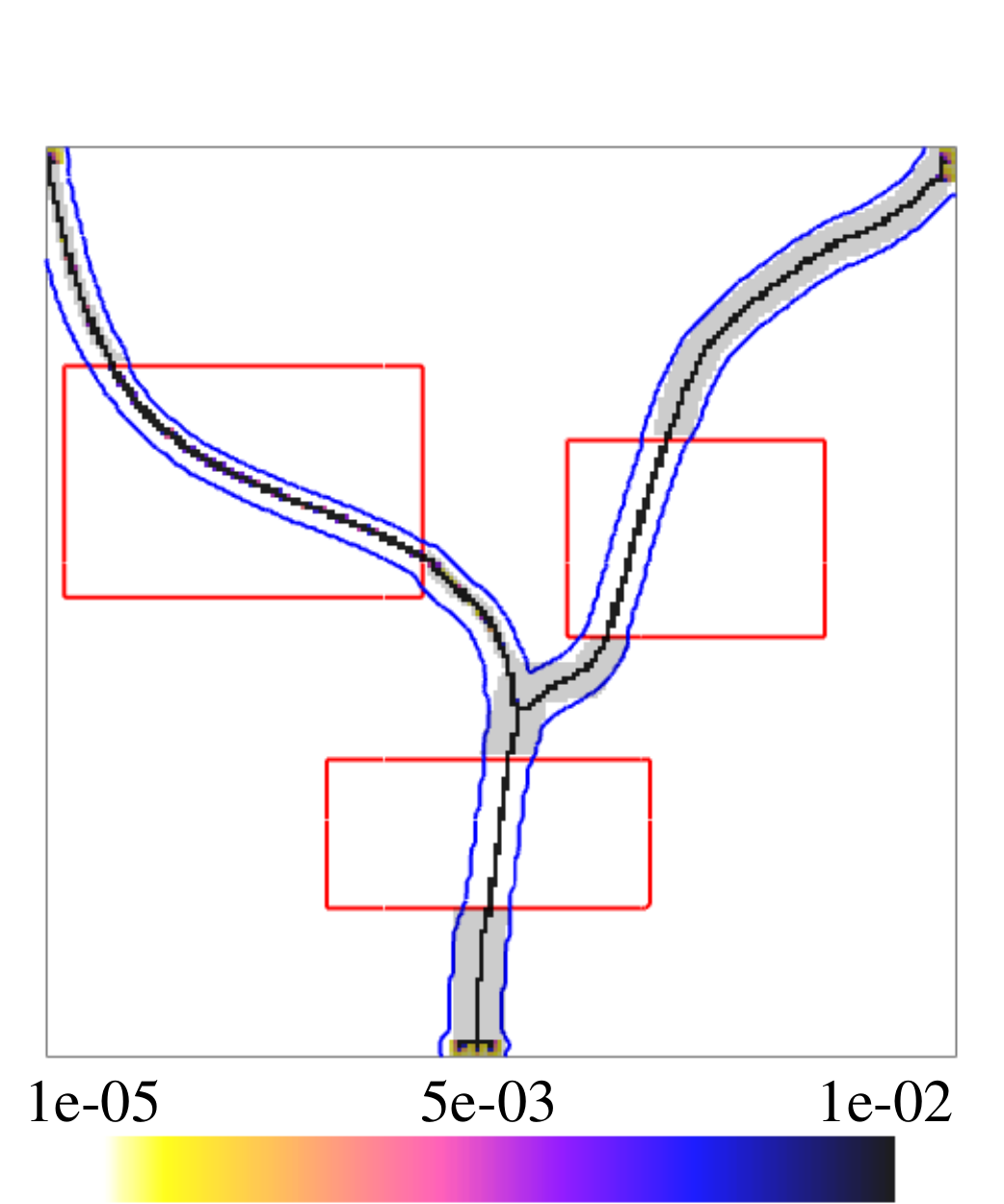}}}}
      &
      \tmpframe{\adjustbox{valign=m,vspace=0pt}{\includegraphics[trim={0.9cm 2.7cm 0.9cm 2.7cm},clip,width=\fraction\columnwidth,valign=t]{{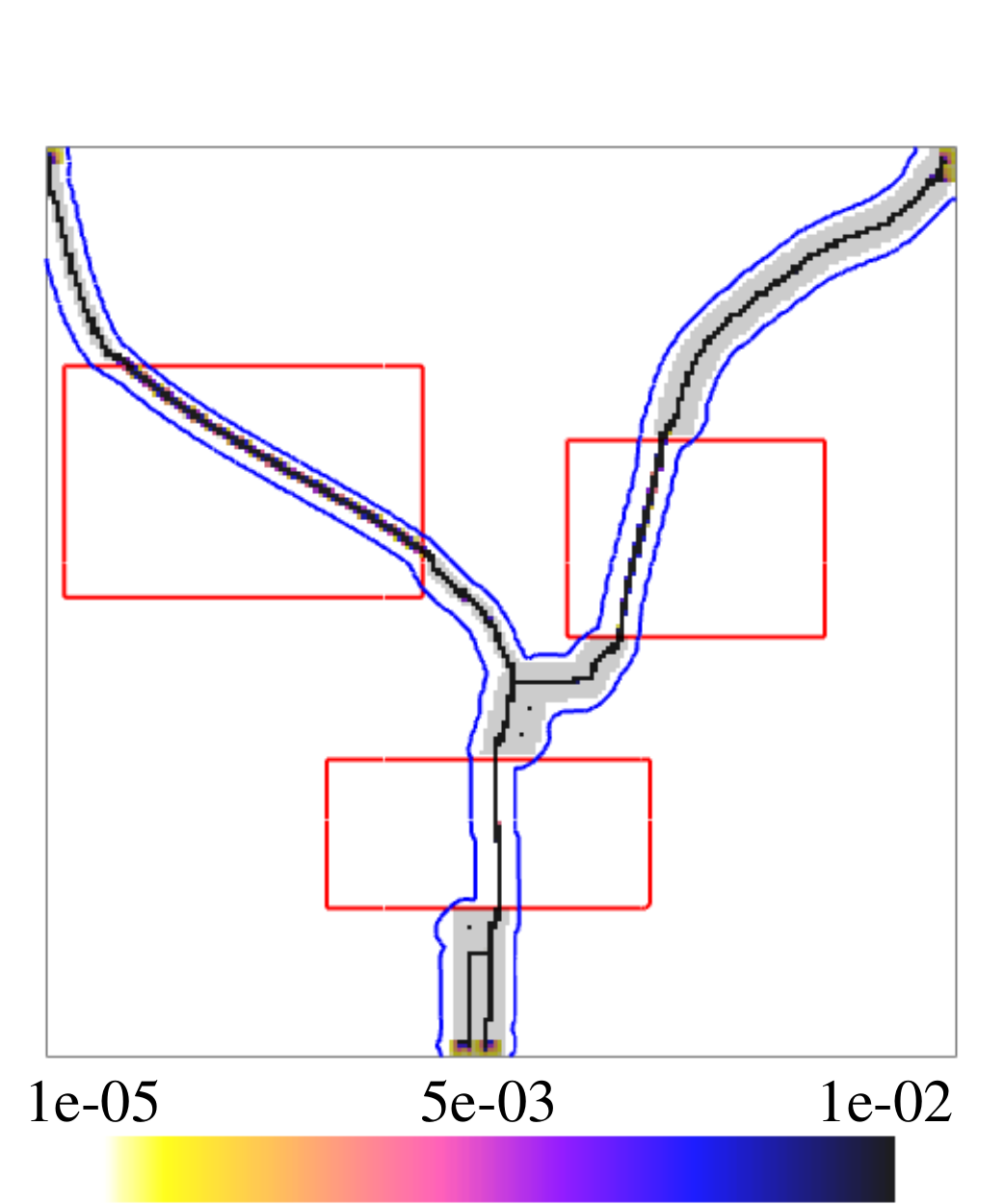}}}}
        \\
    \hline
    \multirow{2}{*}{
      \raisebox{-.5\normalbaselineskip}[0pt][0pt]{\rotatebox[origin=c]{90}{$\confidence=\confmask$, $\Tdens_0=\Obs{\Tdens}$}}
      }
    &
    \raisebox{-.5\normalbaselineskip}[0pt][0pt]{\rotatebox[origin=c]{90}{\quad $\ImgTd(\Tdens)=10\Tdens$}}
    &
      \tmpframe{\adjustbox{valign=m,vspace=0pt}{\includegraphics[trim={0.9cm 2.7cm 0.9cm 2.7cm},clip,width=\fraction\columnwidth,valign=t]{{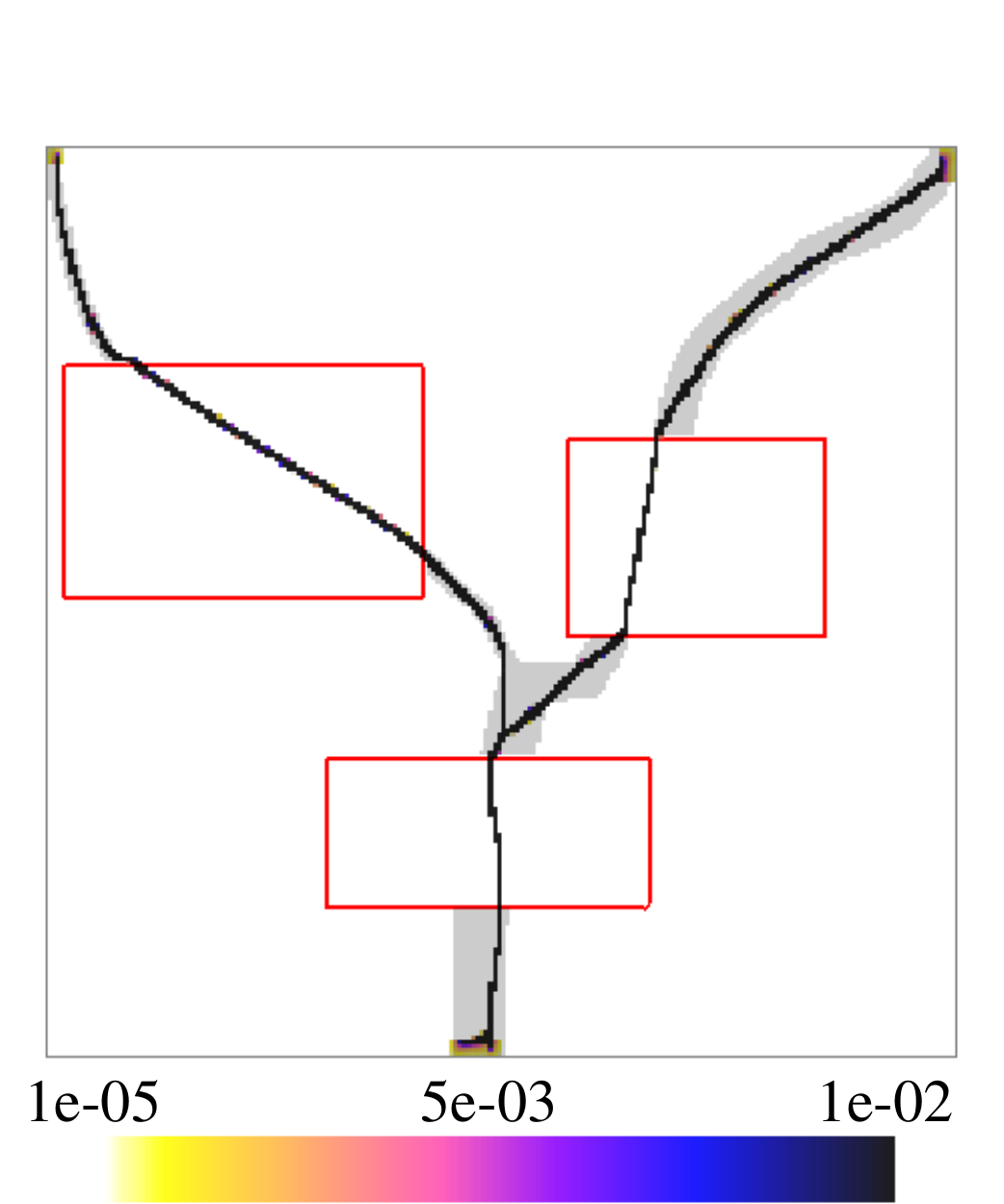}}}}
      &
      \tmpframe{\adjustbox{valign=m,vspace=0pt}{\includegraphics[trim={0.9cm 2.7cm 0.9cm 2.7cm},clip,width=\fraction\columnwidth,valign=t]{{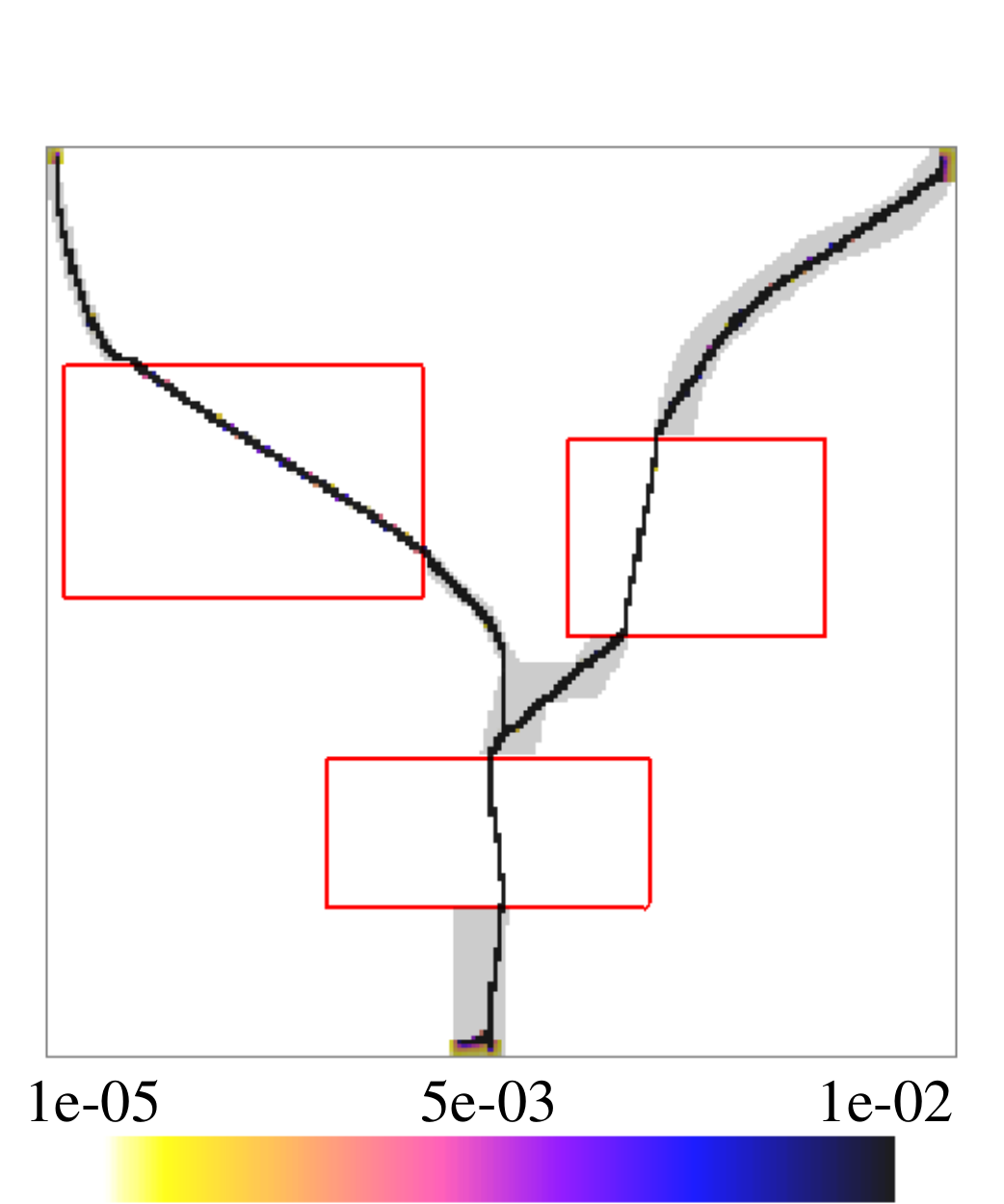}}}}
      &
      \tmpframe{\adjustbox{valign=m,vspace=0pt}{\includegraphics[trim={0.9cm 2.7cm 0.9cm 2.7cm},clip,width=\fraction\columnwidth,valign=t]{{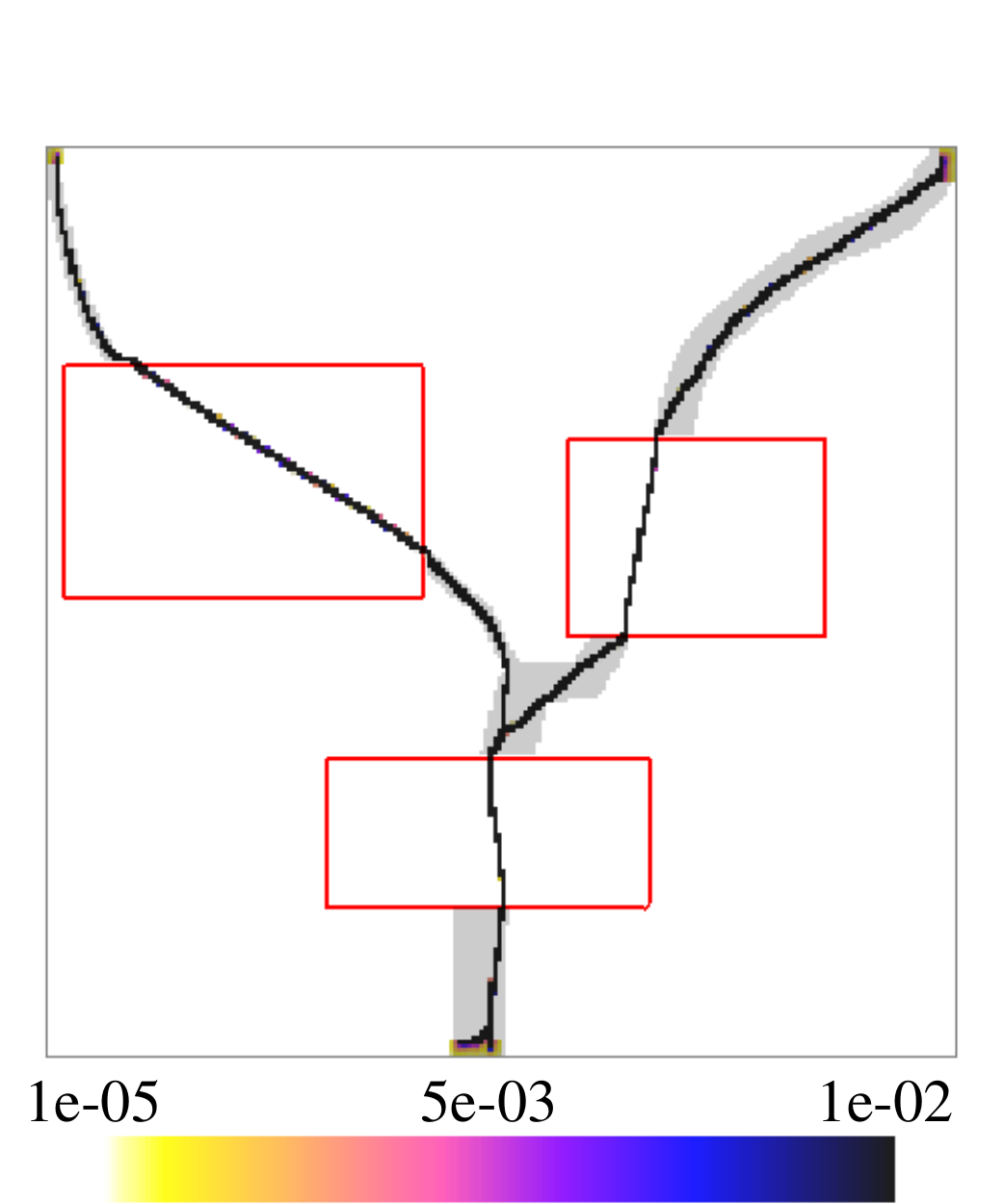}}}}
      &
      \tmpframe{\adjustbox{valign=m,vspace=0pt}{\includegraphics[trim={0.9cm 2.7cm 0.9cm 2.7cm},clip,width=\fraction\columnwidth,valign=t]{{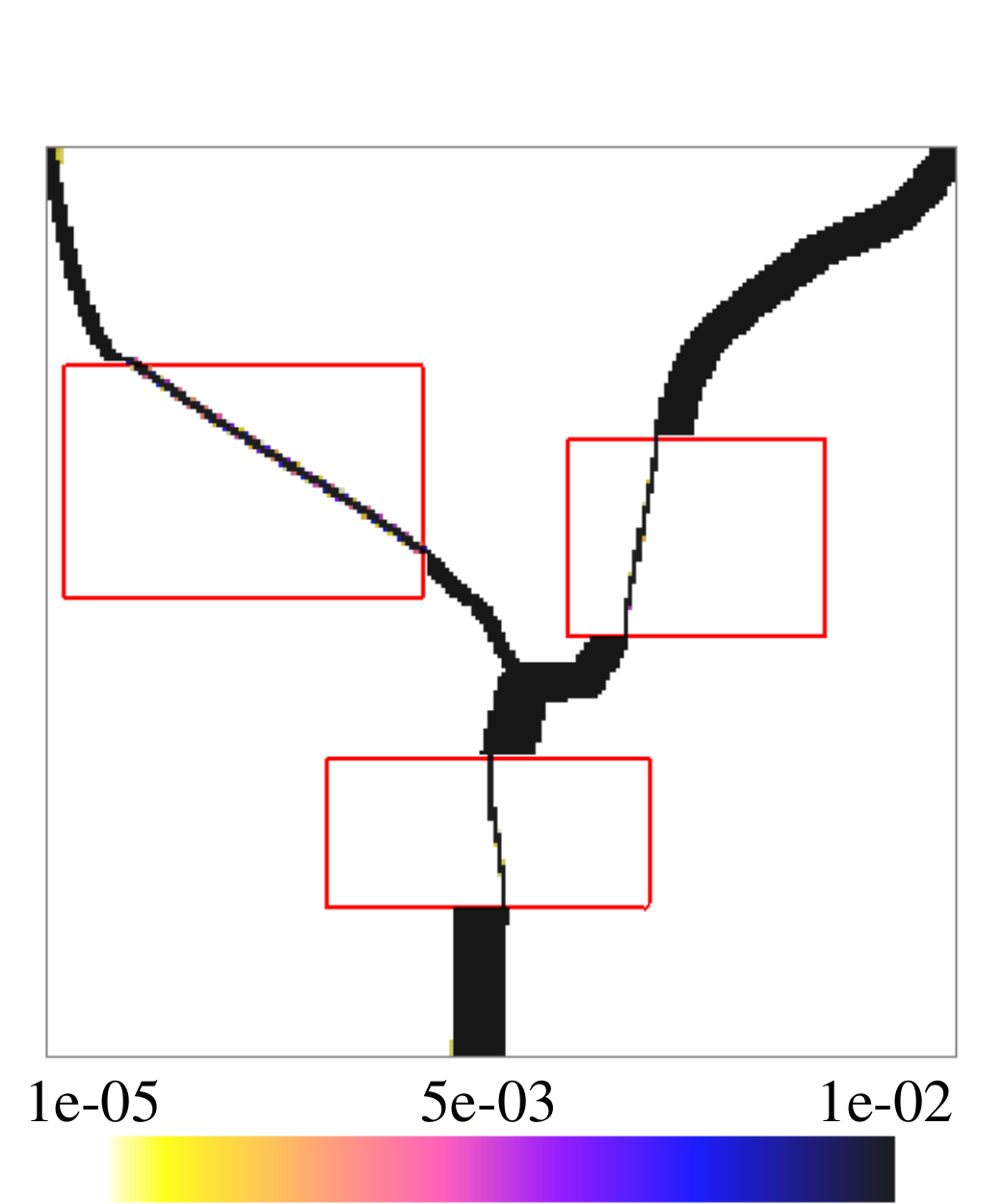}}}}
    \\
    \cline{2-5}
    &
    \raisebox{-.5\normalbaselineskip}[0pt][0pt]{\rotatebox[origin=c]{90}{\quad $\ImgTd(\Tdens)=50\PM(\Tdens)$}}
    &
    \tmpframe{\adjustbox{valign=m,vspace=0pt}{\includegraphics[trim={0.9cm 2.7cm 0.9cm 2.7cm},clip,width=\fraction\columnwidth,valign=t]{{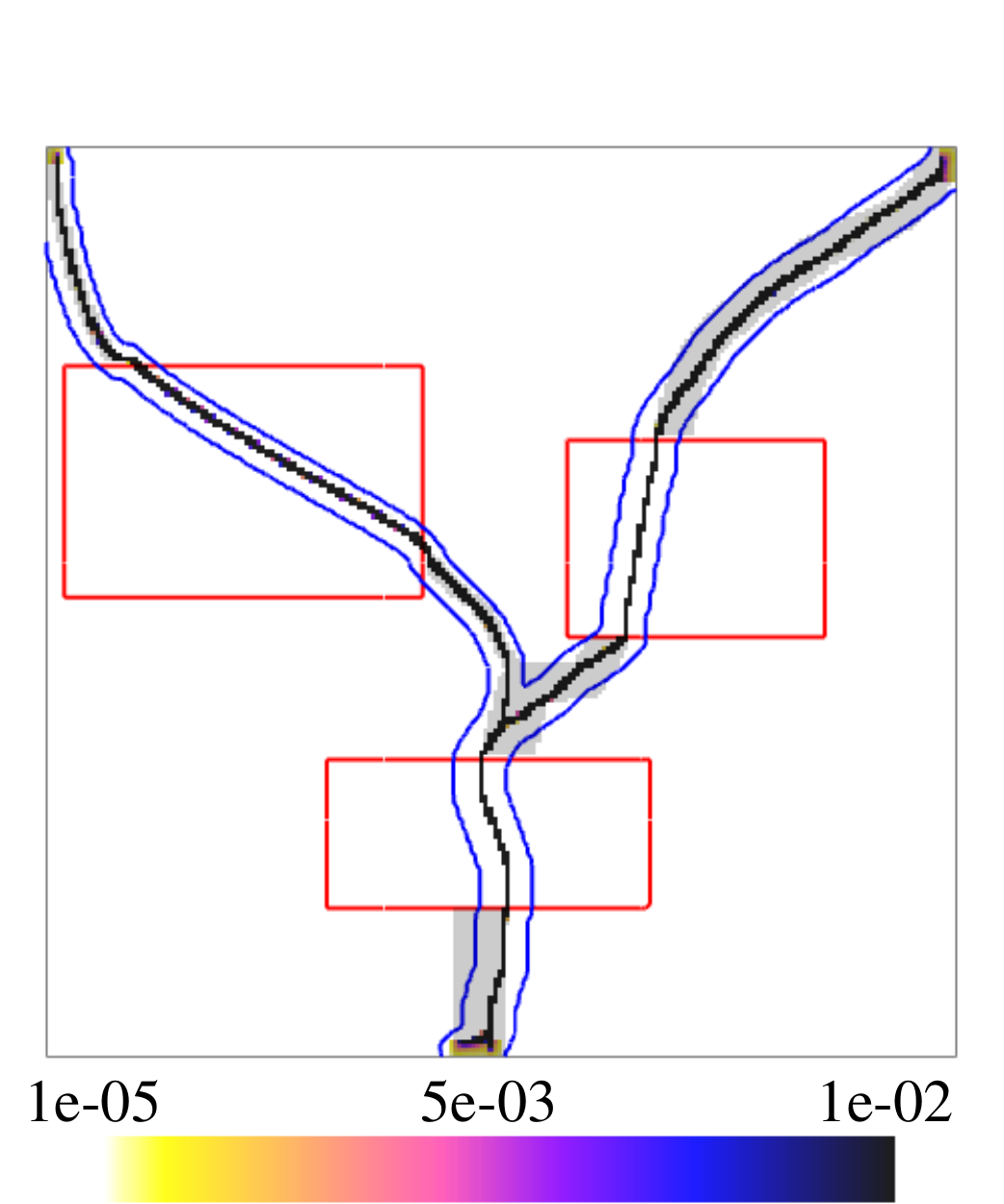}}}}
    &
    \tmpframe{\adjustbox{valign=m,vspace=0pt}{\includegraphics[trim={0.9cm 2.7cm 0.9cm 2.7cm},clip,width=\fraction\columnwidth,valign=t]{{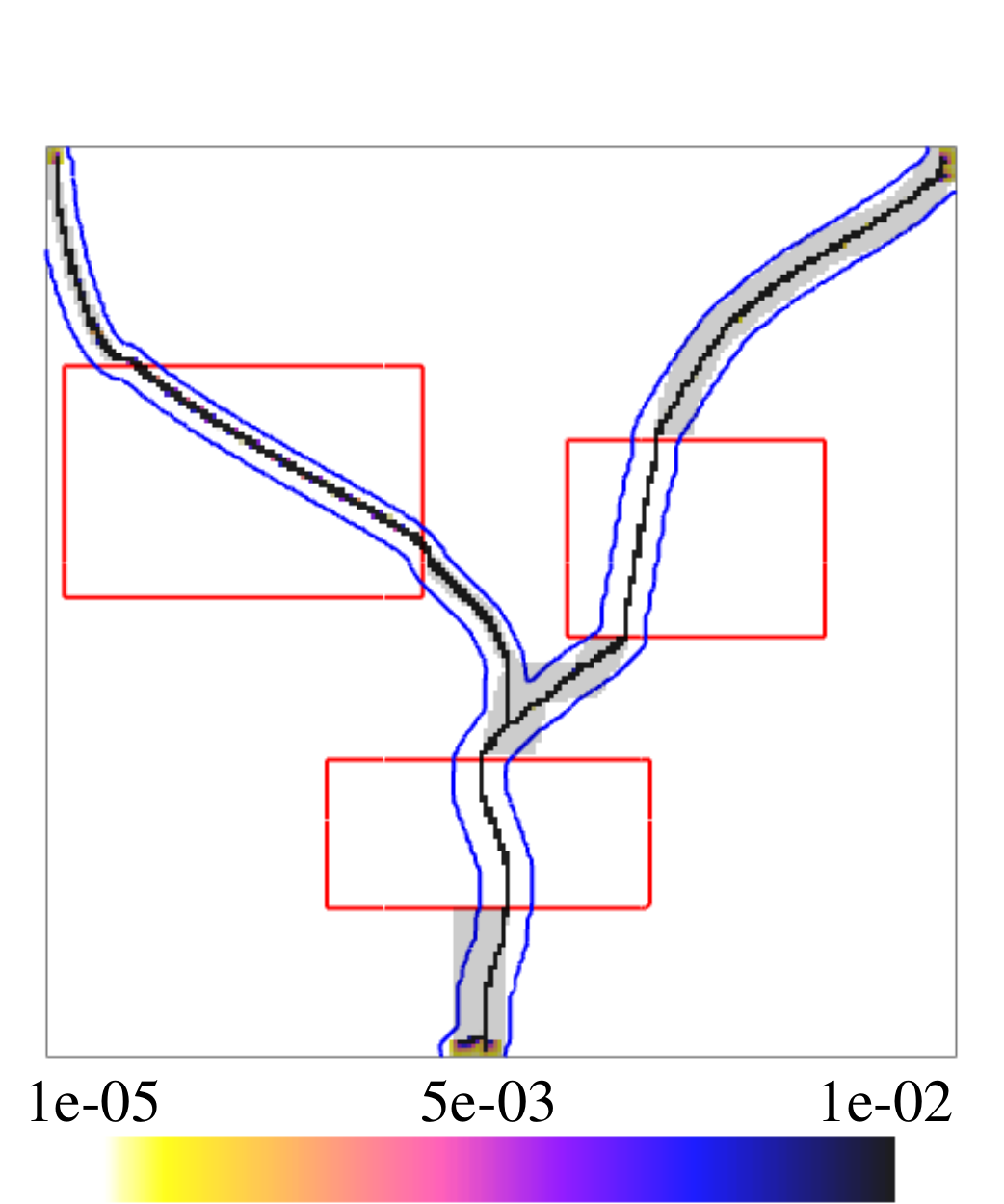}}}}
    &
    \tmpframe{\adjustbox{valign=m,vspace=0pt}{\includegraphics[trim={0.9cm 2.7cm 0.9cm 2.7cm},clip,width=\fraction\columnwidth,valign=t]{{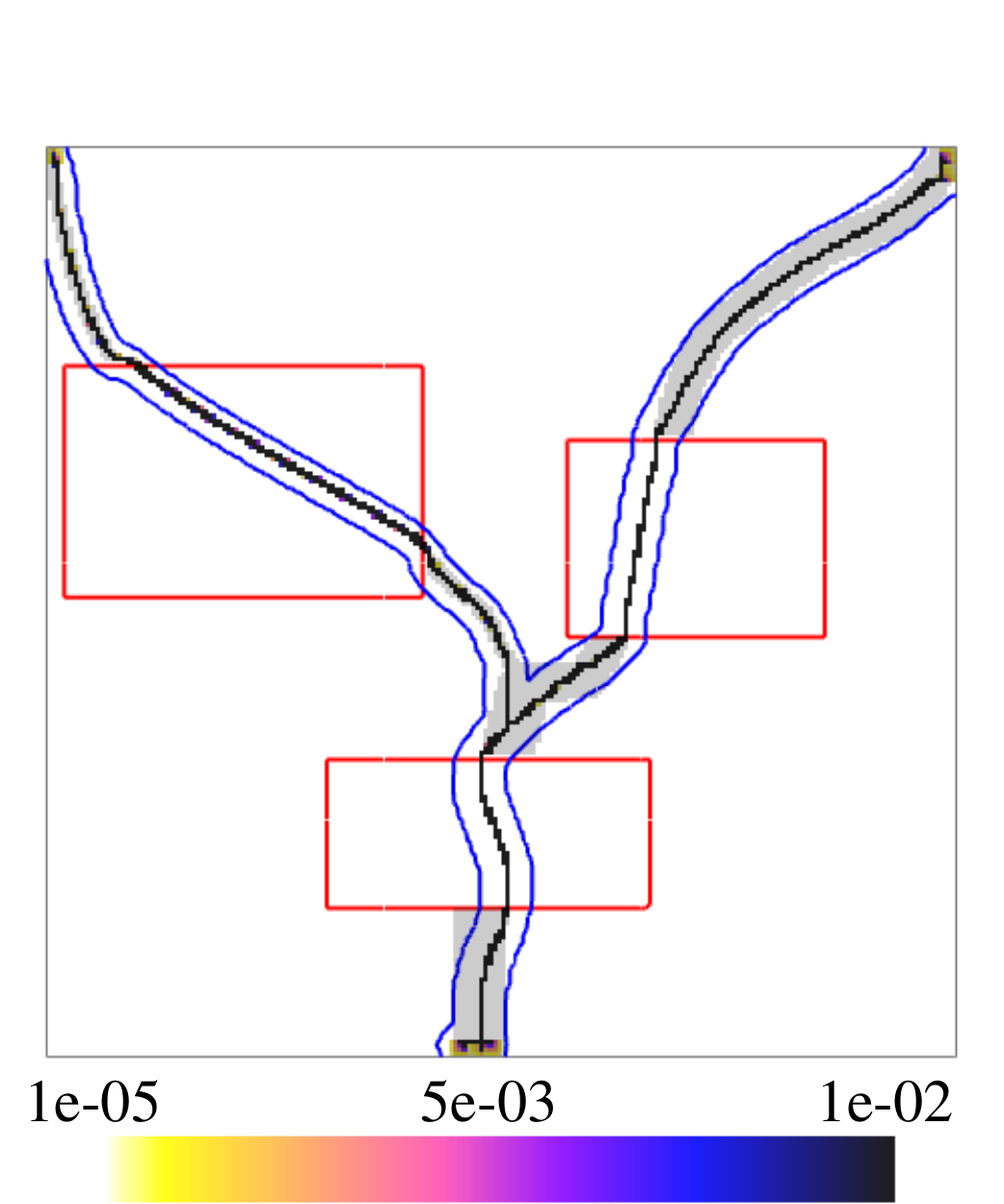}}}}
    &
    \tmpframe{\adjustbox{valign=m,vspace=0pt}{\includegraphics[trim={0.9cm 2.7cm 0.9cm 2.7cm},clip,width=\fraction\columnwidth,valign=t]{{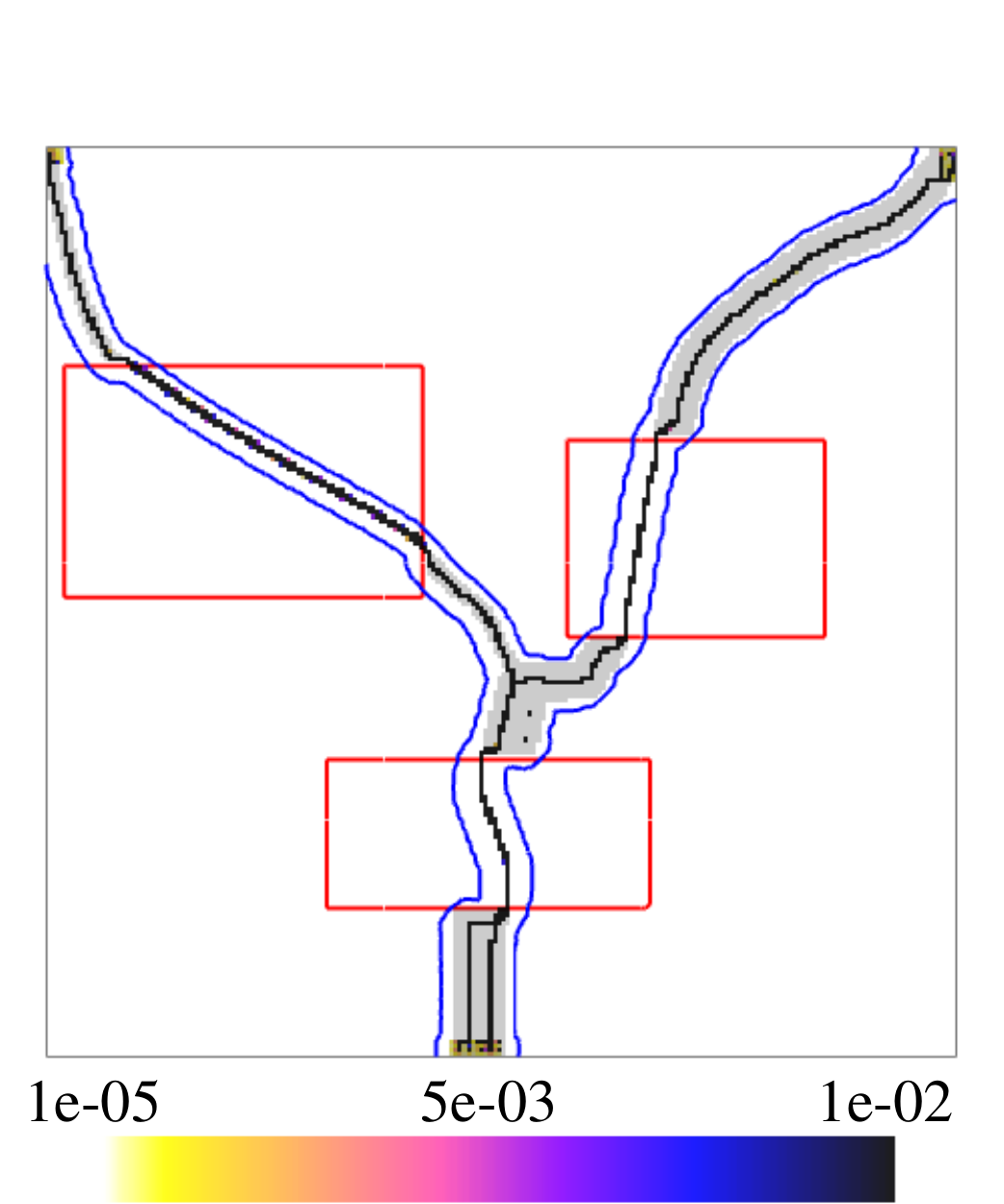}}}}
    \\
    \hline
  \end{tabular}
  \caption{Spatial distribution $\RecTdensH$ for different combinations of $\lambda$, maps, and confidence $\confidence$. When the \PM{} map is used, the blue line defines the boundary for which $\Rec{\Img}=\ImgTd(\RecTdensH)$ is $1e-3$. The red rectangle defines the mask region. The observed image $\Obs{\Img}$ is reported in light gray.}
  \label{fig:netnetwork-all}
  \end{figure}
}

We can see how the proposed approach effectively restores the connectivity of the network for all values of $\lambda$ and both maps. This result is the consequence of imposing a flux between the source term $\Source=\Dirac{O}$ and the sink term $\Sink=(\Dirac{P}+2\Dirac{Q})/3$. Having a disconnected network, which in our setting means having zero conductivity, would contradict the constraint given in~\cref{eq:poisson}. 

Since $\confidence=\confmask$, the discrepancy term has no influence within the rectangles that make up the mask, and the optimization of the regularization term $\Reg$ is the only factor driving the reconstruction. Thus, inside the mask, the reconstructed network is composed of (almost) straight channels or a union of them, reconnecting the ends of the original network cut by the mask. In this region, we can see some grid alignment artifacts.

Outside the mask, we see the influence of the discrepancy term. The network progressively fits to the observed data as $\lambda$ increases. In some cases we see some non-optimal structures such as parallel channels using the identity map with $\lambda=5e-2$ and the \PM{} map with $\lambda=1.0$. This is due to the fact that the fitting term tends to match the value of $\Obs{\Img}$, opposed to the concentrating effect coming from the \BT{} functional $\Reg$.

For both maps, the topology of the reconstructed network can change drastically even with a relatively small change of $\lambda$. For example, as $\lambda$ increases from $1e-2$ to $5e-2$, the upper right branched is fitted. Similarly, the bifurcation point in the middle of the domain is fitted by increasing $\lambda$ from $5e-2$ to $1e-1$. This sudden transition can be attributed to the fact that the functional $\Lyap$ is not convex and that we are using gradient-based optimization methods. This leads to diverse convergence paths within the optimization landscape. 

Both maps produced comparable results in terms of reconstructing the primary topological feature of the network, showing a similar change when $\lambda$ is altered.  However, the \PM{} map provides a smoother transition of the reconstructed network $\Rec{\Img}$ at the mask boundary. Moreover, it generates both skeletonization and a reconstruction of the network. 

\subsubsection{Images with lost data (case $\confidence=1$)}
\label{sec:netnetwork-confONE}
In this section, we propose a second experiment where we use the confidence term $\confidence\equiv 1$, more suitable for those problems where we are not aware of the location of the corruption. In the third and fourth rows in~\cref{fig:netnetwork-all} we report the spatial distribution of $\RecTdensH$ while changing the scaling parameter $\lambda$ of the discrepancy and the map $\ImgTd$.

The results do not differ significantly from those obtained using $\confidence=\confmask$. However, $\RecTdensH$ tends to have lower values and smaller support within the mask relative to the results obtained using $\confidence=\confmask$. For example, in the figure showing the combination of using $\confidence=\confmask$, the identity map, and $\lambda=5e-2$, the branching point within the lowest rectangle is moved upward. This effect can be attributed to the fact that the discrepancy term is now acting also within the mask $\Mask$, where the data $\Obs{\Img}$ fitted is equal to zero, penalizing the presence of the network. If $\RecTdensH$ is still present within the mask, this penalization effect results in a network with wider channels, as a wider section is required to carry the same flux. This is more evident when comparing the case with the identity map and $\lambda=1e0$.

\subsubsection{Influence of initial data $\Tdens_0$}
\label{sec:tdensini-obs}
In~\cref{sec:dmk} we already mentioned how the approximate solution $\OptTdensH$ of~\cref{prob:dmk} depends on the initial data. In this part, we illustrate the impact of this issue on the NIOT algorithm.
So far, our consideration has been limited to uniform data with $\Tdens_0=1$. However, in this instance, we employ initial data that is extrapolated from the observed data. In particular, we set
\begin{equation}
  \label{eq:tdensini-obs}
  \Tdens_0=\Obs{\Tdens}=\Scaling^{-1}(\Obs{\Img})+\TdensLift,
\end{equation}
were the lower bound $\TdensLift$ is required to avoid a strong penalization in the initial data where data are lost. In our experiment, we set it to $1e-5$ times the maximum value of $\Obs{\Img}$. The results are reported in the last two rows of~\cref{fig:netnetwork-all}. 

It is evident that the selection of initial data profoundly impacts the ultimate outcome. The magnitude of this influence is such that the selection of the map $\ImgTd$ and the scaling of the discrepancy $\lambda$ become almost insignificant. This effect can be beneficial when images are corrupted only by data loss. However, if this assumption is not valid, using such initial data may result in incorrect outcomes, as shown in the next section.

\subsubsection{Images with artifacts}
\label{sec:artifacts} 
In this section, we focus on a reconstruction problem in which false positive data are present in the image. Additional channels are included in~\cref{fig:network-inpainting-true} and then the mask shown in~\cref{fig:network-inpainting-true} is applied. The resulting $\Obs{\Img}$ can be seen in~\cref{fig:network-inpainting-obs}.

{
\def \fraction {0.235}
\def \extrah {1pt}
\begin{figure}
  \centering
  \begin{tabular}{|@{}c@{}|@{\hspace{\extrah}}c@{\hspace{\extrah}}|@{\hspace{\extrah}}c@{\hspace{\extrah}}|@{\hspace{\extrah}}c@{\hspace{\extrah}}|@{\hspace{\extrah}}c@{\hspace{\extrah}}|}
    \hline
    & $\lambda=1e-3$ & $\lambda=5e-2$  &  $\lambda=1e-1$ & $\lambda=1e0$ 
    \\
    \hline
    \raisebox{-.5\normalbaselineskip}[0pt][0pt]{\rotatebox[origin=c]{90}{$\Tdens_0=1$}}
    &
    \tmpframe{\adjustbox{valign=m,vspace=0pt}{\includegraphics[trim={0.9cm 2.7cm 0.9cm 2.7cm},clip,width=\fraction\columnwidth,valign=t]{{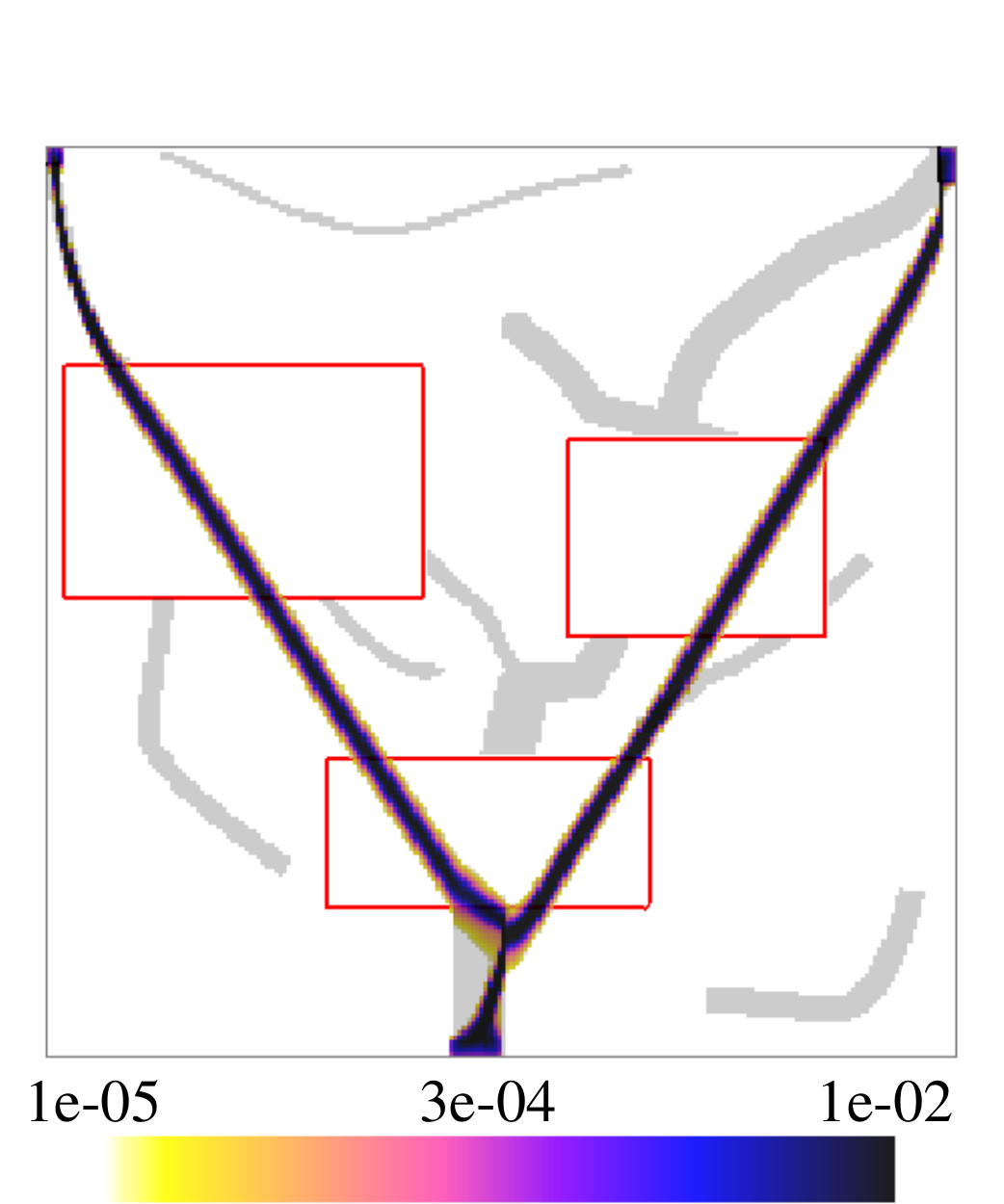}}}}
    & 
    \tmpframe{\adjustbox{valign=m,vspace=0pt}{\includegraphics[trim={0.9cm 2.7cm 0.9cm 2.7cm},clip,width=\fraction\columnwidth,valign=t]{{imgs/y_net_nref2/mask_medium/matrix_nref0_femDG0DG0_gamma5.0e-01_wd5.0e-02_wr0.0e+00_netnetwork_artifacts_ini0.0e+00_confONE_mu2iidentity_scaling1.0e+01_methodte_tdens.pdf}}}}
    & 
    \tmpframe{\adjustbox{valign=m,vspace=0pt}{\includegraphics[trim={0.9cm 2.7cm 0.9cm 2.7cm},clip,width=\fraction\columnwidth,valign=t]{{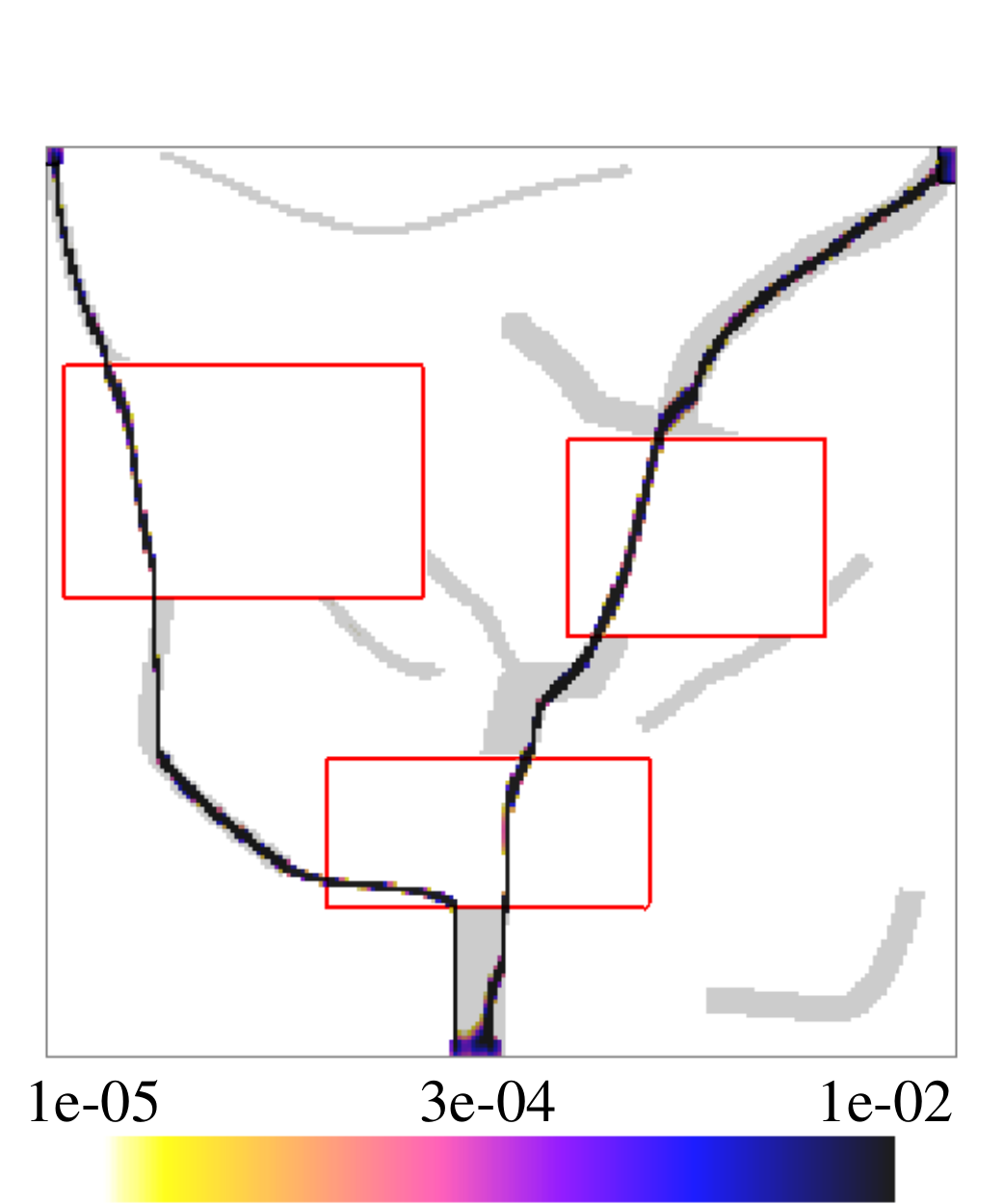}}}}
    & 
    \tmpframe{\adjustbox{valign=m,vspace=0pt}{\includegraphics[trim={0.9cm 2.7cm 0.9cm 2.7cm},clip,width=\fraction\columnwidth,valign=t]{{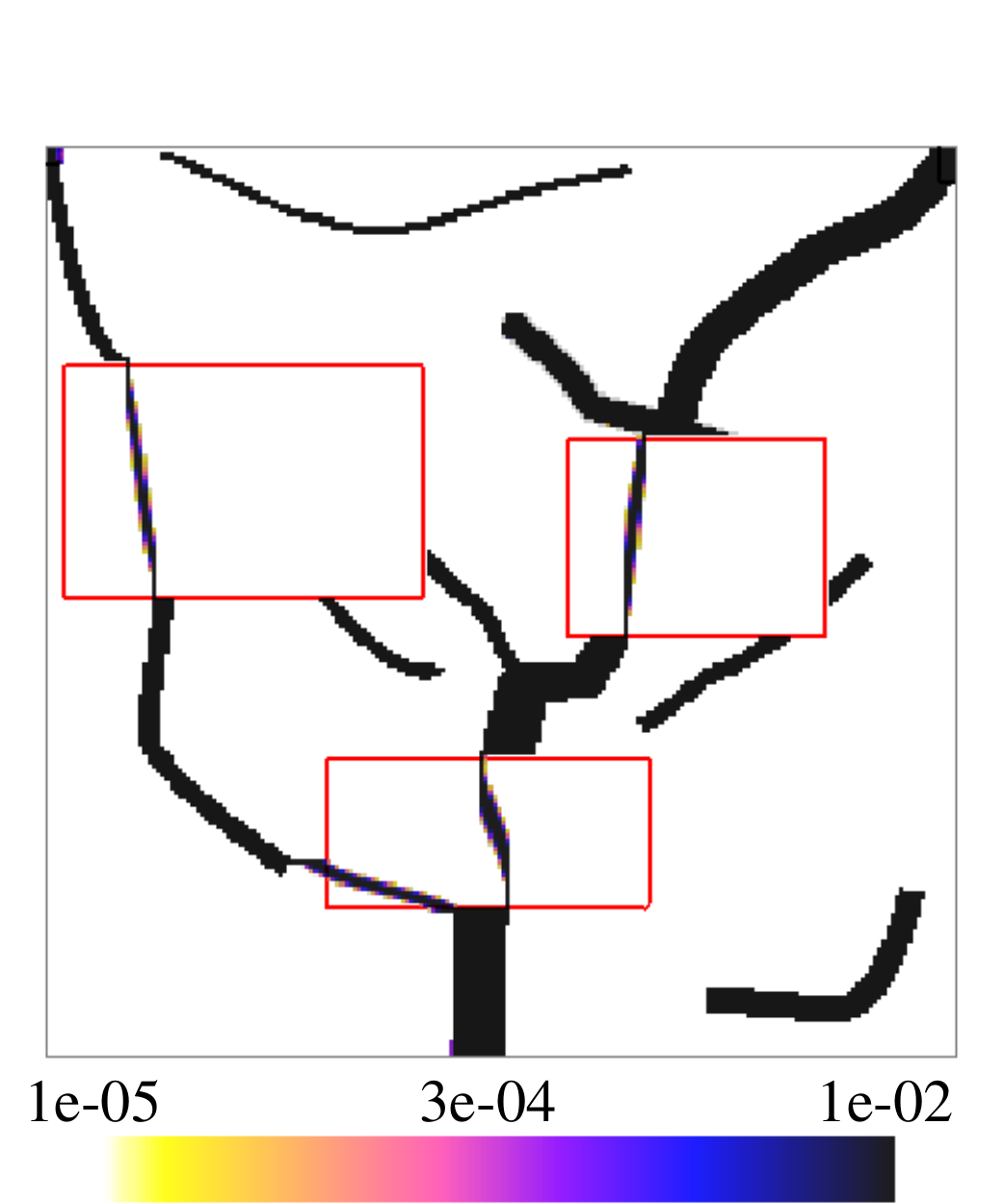}}}}
    \\
    \hline
    \raisebox{-.5\normalbaselineskip}[0pt][0pt]{\rotatebox[origin=c]{90}{$\Tdens_0=\Obs{\Tdens}$}}
    &
    \tmpframe{\adjustbox{valign=m,vspace=0pt}{\includegraphics[trim={0.9cm 2.7cm 0.9cm 2.7cm},clip,width=\fraction\columnwidth,valign=t]{{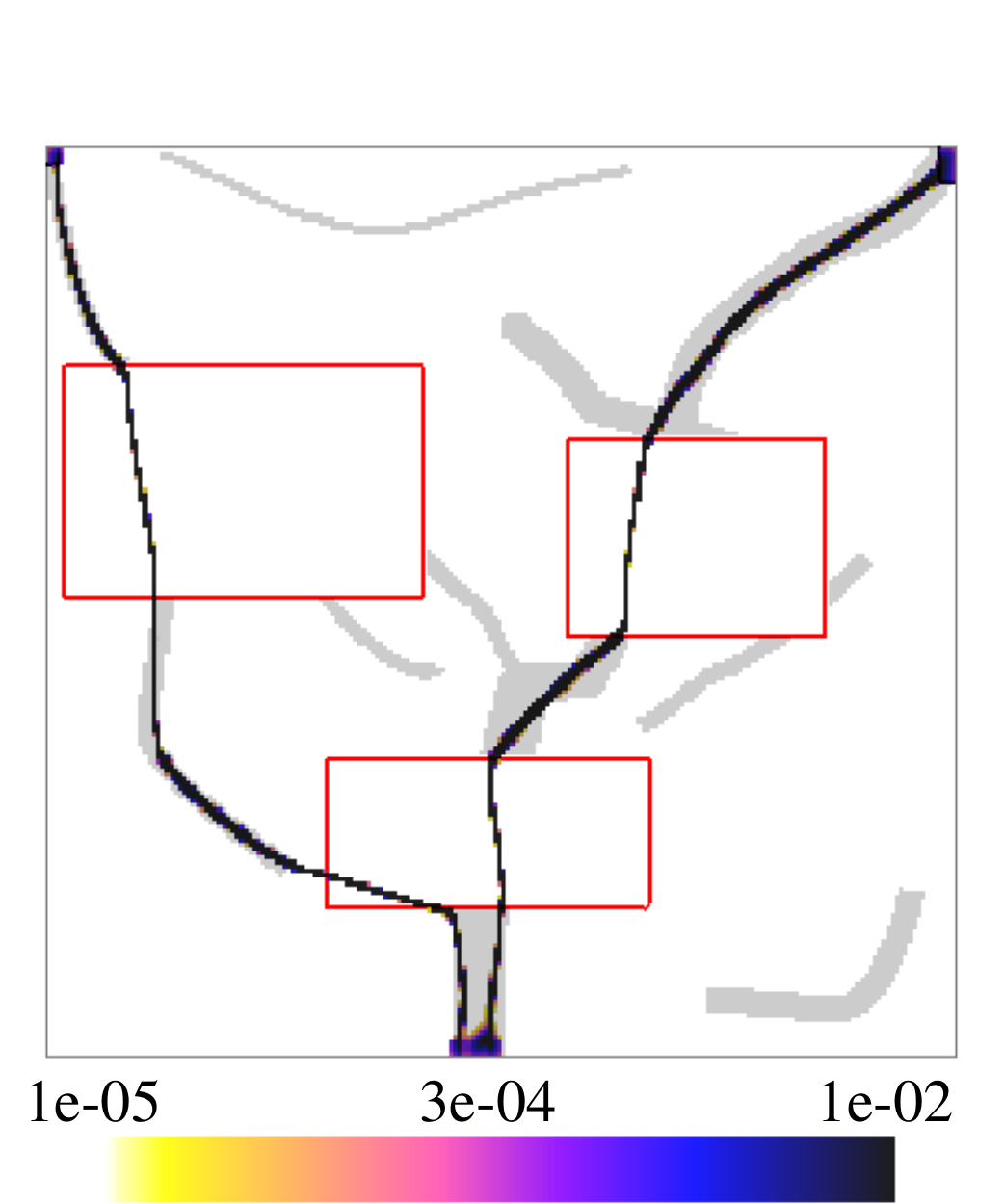}}}}
    & 
    \tmpframe{\adjustbox{valign=m,vspace=0pt}{\includegraphics[trim={0.9cm 2.7cm 0.9cm 2.7cm},clip,width=\fraction\columnwidth,valign=t]{{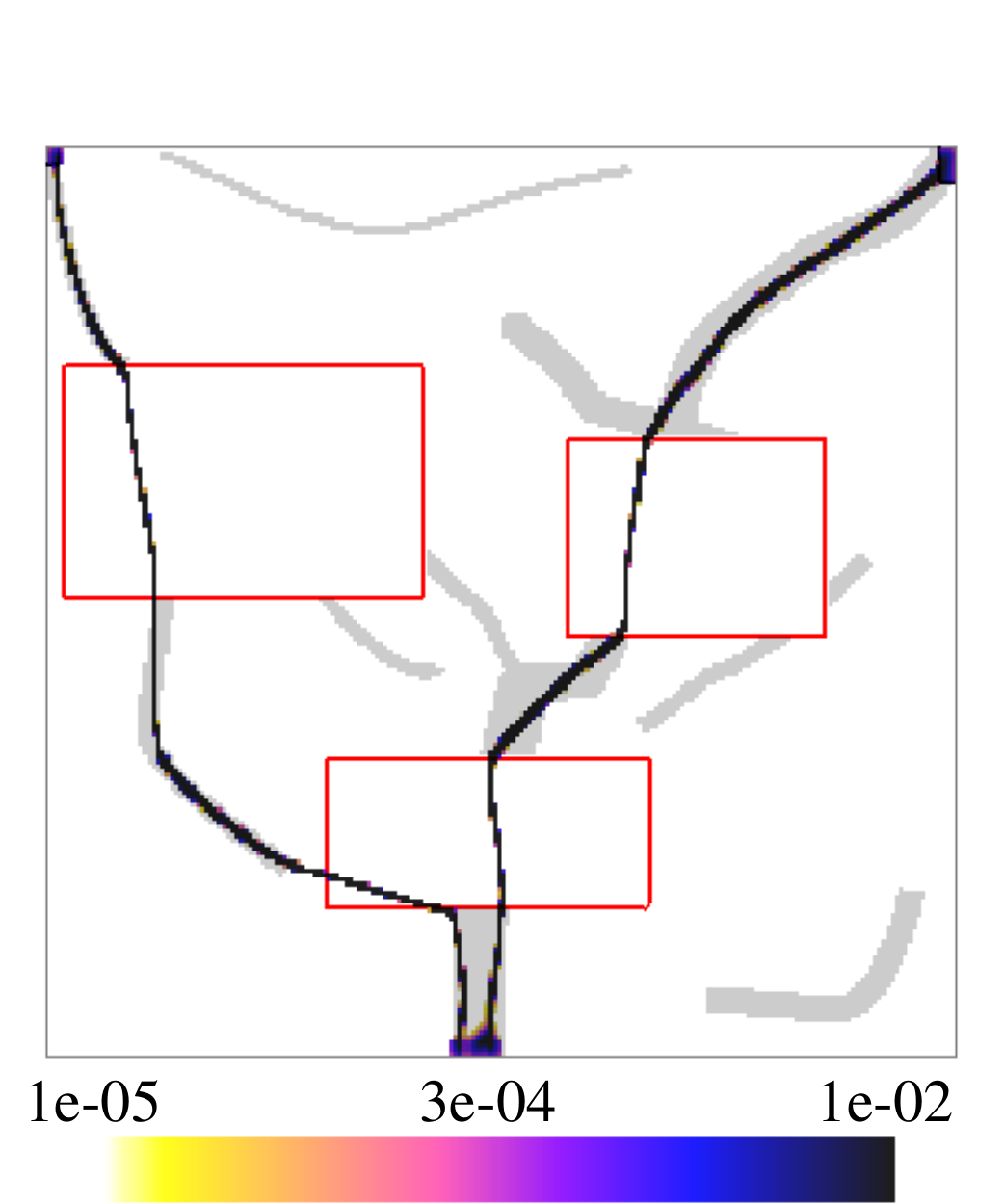}}}}
    & 
    \tmpframe{\adjustbox{valign=m,vspace=0pt}{\includegraphics[trim={0.9cm 2.7cm 0.9cm 2.7cm},clip,width=\fraction\columnwidth,valign=t]{{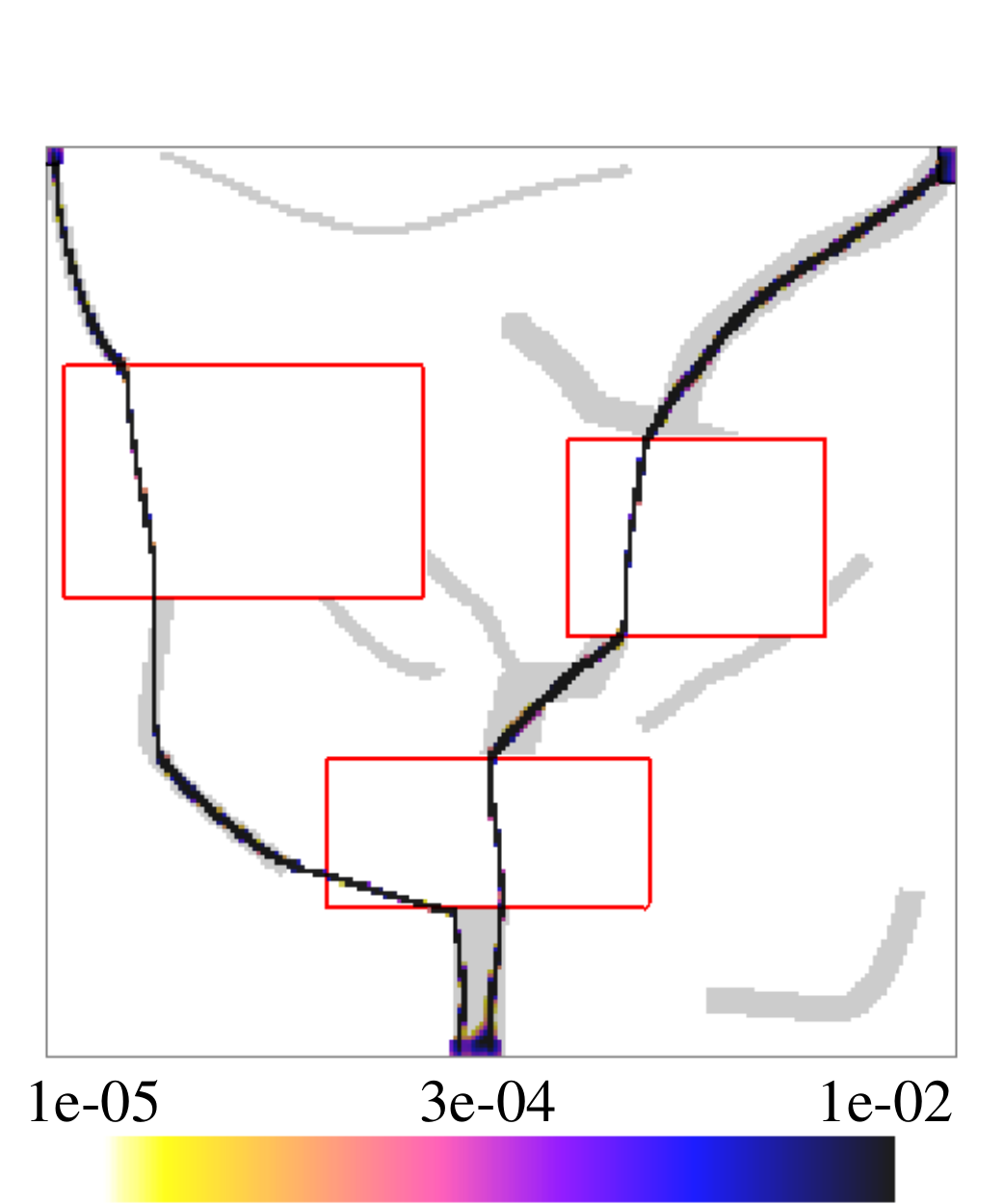}}}}
    & 
    \tmpframe{\adjustbox{valign=m,vspace=0pt}{\includegraphics[trim={0.9cm 2.7cm 0.9cm 2.7cm},clip,width=\fraction\columnwidth,valign=t]{{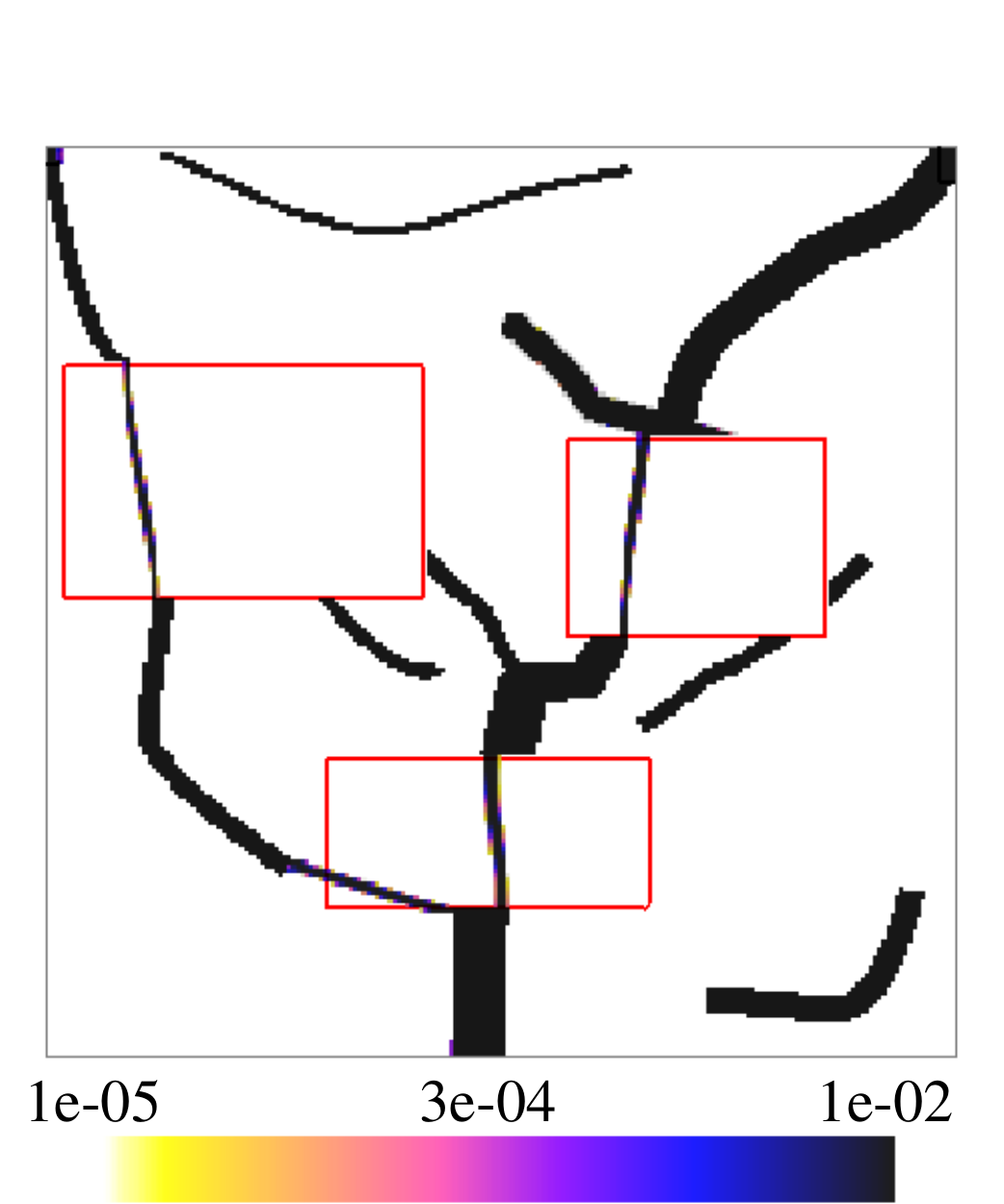}}}}\\
    \hline
  \end{tabular}
  \caption{Spatial distribution $\RecTdensH$ for $\lambda=1e-3,5e-2,1e-1,1e0$ (left to right) and $\Tdens_0=1,\Obs{\Tdens}$ (top to bottom). Confidence $\confidence=1$ and identity map. The observed data $\Obs{\Img}$ is obtained by adding additional artifacts to the original network in~\cref{fig:network-inpainting-true}, and then removing them within the mask (red rectangles). The result $\Obs{\Img}$ is reported in light gray.
}
\label{fig:artifacts}
\end{figure}
}

In this case, we set the confidence measure $\confidence=1$, since we typically do not know the location of false positive data. We use the identity map with $\alpha=1e1$. In~\cref{fig:artifacts}, we report the reconstructed networks using different $\lambda>0$ and initial data $\Tdens_0=1,\Obs{\Tdens}$.
 
The results show how large values of $\lambda$ (rightmost panels) lead to overfitting of the data, including the artifacts that we would like to remove. On the other hand, a network more similar to the original Y-shaped network is obtained with $\lambda=1e-3,5e-2$ and unbiased initial data $\Tdens_0=1$. This suggests that the usage of smaller values of $\lambda$ and uniform initial data should be preferred if false positive data is present. 

These experiments show an additional feature of the proposed method, its capability to remove artifacts from the corrupted data. This procedure is referred to in the literature on image processing as pruning \cite{russ2006image}. In our model, this selection is the result of an optimization process, where the presence of the functional $\Lyap$ removes those channels that do not carry enough flux. Similar ideas can be found in~\cite{baptista2020network} in the contest of graph pruning.

\begin{remark}
  We considered only two types of initial data, uniform data $\Tdens_0=1$ and $\Tdens_0=\Obs{\Tdens}$, with the latter strongly biased by the observed data. However, the choice is not necessarily binary. In fact, it is easy to image a transition between the two configurations by progressively blurring $\Obs{\Tdens}$ with a Gaussian filter. With such a procedure, the initial data is less biased by the observed data but still contains some information about the network that may be useful in the reconstruction process.
\end{remark}

\subsection{Frog tongue}
\label{sec:frog_tongue}
As an example of a real biological geometry with a reasonable 2D spatial representation, we turn to the classical work of Julius Cohnheim from 1872 \cite{cohnheim1872investigations}. In his work, Cohnheim describes the vasculature of a frog tongue and includes a detailed colored drawing. As this geometry is 2D in nature, we can represent results without accounting for out-of-plane connections. The original image is shown in~\cref{fig:cohnheim}. From this drawing, we segment the arterial network (shown in red in~\cref{fig:frog-network-cohnheim}) and the tongue itself. Segmentations are performed using a semi-automatic procedure used in~\cite{hodneland2019}, followed by a manual inspection to ensure spatial connectivity of the arterial vessel tree.

\def \fraction {0.4}
\def \subfraction {0.99}
\begin{figure}
  \centering
  \begin{subfigure}[t]{\fraction\textwidth}
  \centerline{
     \includegraphics[width=\subfraction\columnwidth,valign=t]{{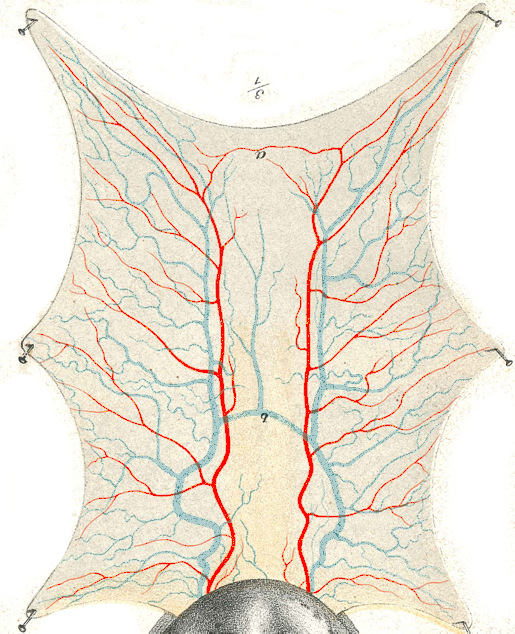}}
  }
  \caption{Original drawing by Cohnheim with the arterial (red) and venal (blue) blood vessel network. The binary map of the arterial network is $\Obs{\Img}$.}
  \label{fig:frog-network-cohnheim}
  \end{subfigure}
\qquad
  \begin{subfigure}[t]{\fraction\textwidth}
  \centerline{
     \includegraphics[trim={.4cm 0cm 0.2cm 2.6cm},clip, width=\subfraction\columnwidth,valign=t]{{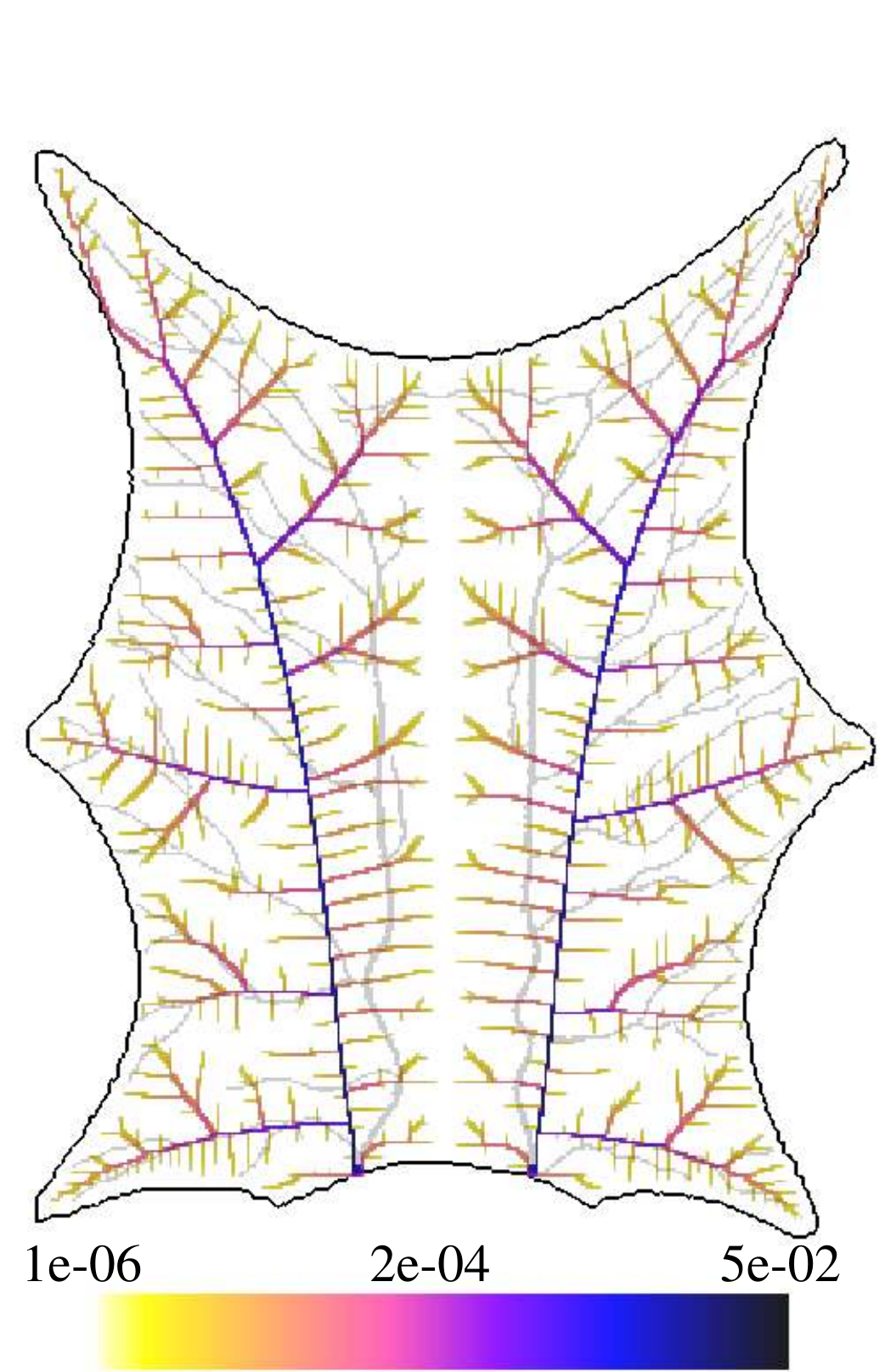}}
  }
  \caption{Spatial distribution of $\OptTdensH$
  solving for $\gamma=0.5$ and $\Tdens_0=1$.}
  \label{fig:frog-network-opttdens}
  \end{subfigure}
  \begin{subfigure}[t]{\fraction\textwidth}
  \centerline{
    \includegraphics[trim={.4cm 0cm 0.2cm 2.6cm},clip, width=\subfraction\columnwidth,valign=t]{{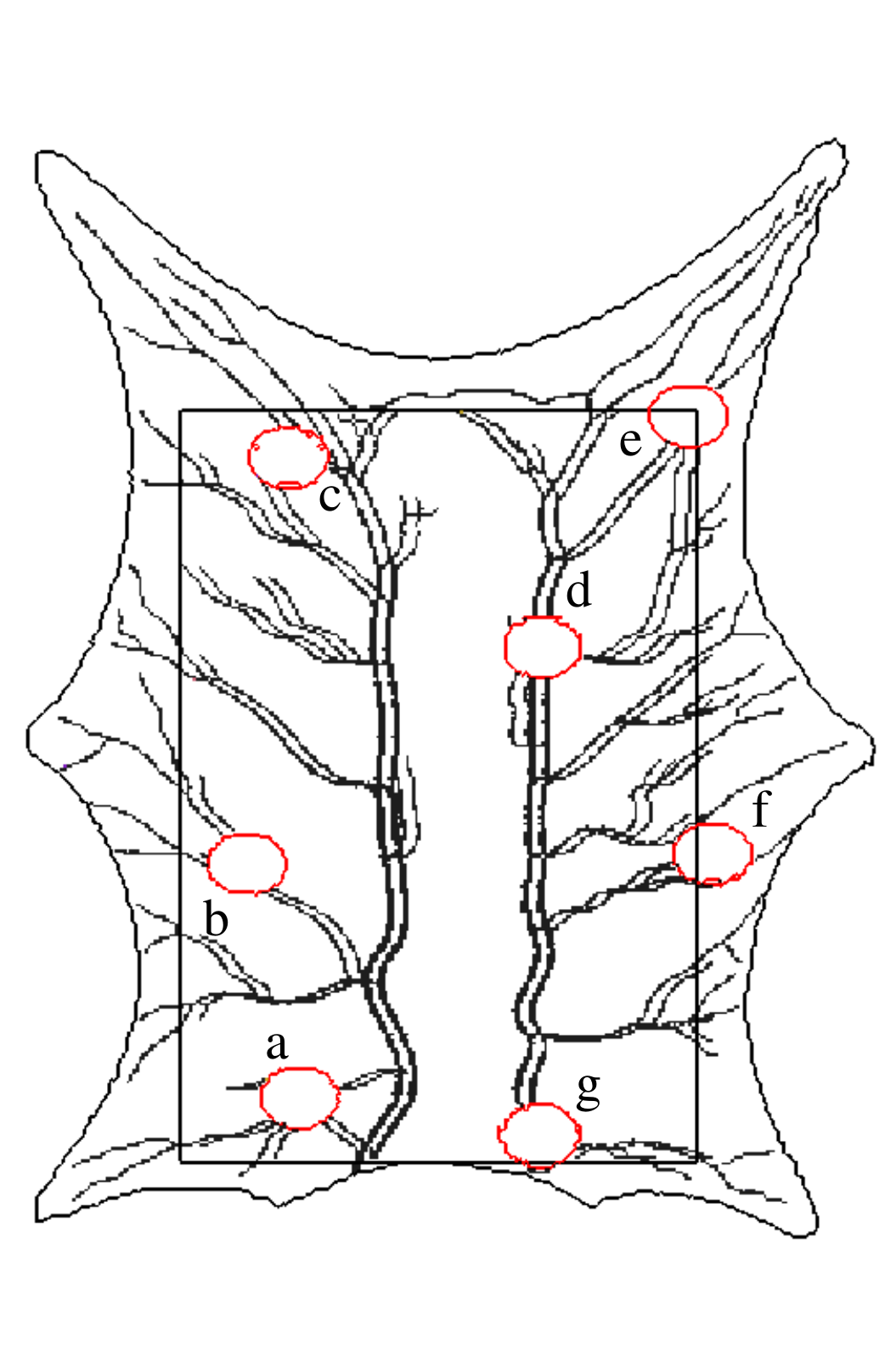}}
  }
  \caption{$\Obs{\Img}$ with artifacts and missing data used in~\cref{sec:frog-shifted}. The red circles mark the location of the mask used.}
  \label{fig:frog-network-shifted}
  \end{subfigure}
  \qquad
  \begin{subfigure}[t]{\fraction\textwidth}
    \centerline{
    \includegraphics[trim={.4cm 0cm 0.2cm 2.6cm},clip, width=\subfraction\columnwidth,valign=t]{{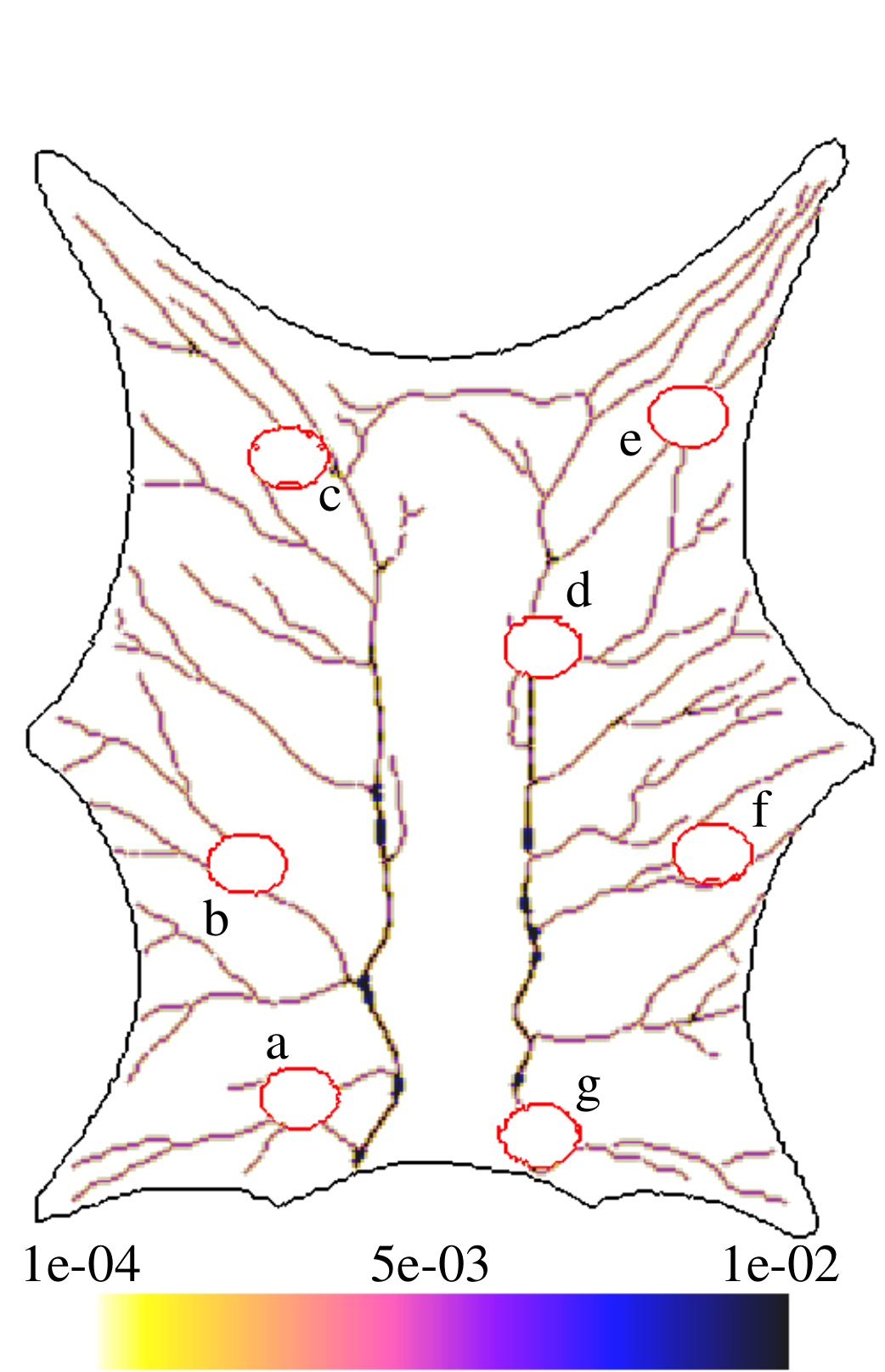}}
  }
  \caption{Spatial distribution of $\Obs{\Tdens}$ used in~\cref{sec:frog-mu}.}
  \label{fig:frog-network-mupou}
  \end{subfigure}
  \caption{Frog tongue data. In~\cref{fig:frog-network-cohnheim}, the original drawing from~\cite{cohnheim1872investigations}. In~\cref{fig:frog-network-opttdens}, the spatial distribution of $\OptTdensH$ solving~\cref{prob:dmk} 
  for $\Sink=1$ within the contour area and $\Source$ being the sum of two piecewise constant functions located at the root of the two main channels. In~\cref{fig:frog-network-shifted,fig:frog-network-mupou}, the observed data $\Obs{\Img}$ used in~\cref{sec:frog-shifted,sec:frog-mu}, respectively.
  The circles indices the area where the mask is applied, the mask also used in~\cref{sec:frog-network}.
  }
  \label{fig:cohnheim}
\end{figure}

We deduce the terms $\Source$ and $\Sink$ used in~\cref{prob:dmk} directly from~\cref{fig:cohnheim}. We set $\Source$ as the sum of two piecewise constant functions located at the root of the two main channels, having the same mass injection rate and support. The sink term $\Sink$ is set equal to one on the support of the frog tongue, modeling a uniform absorption of nutrients/oxygen. The two terms are balanced to ensure that the total mass injected is equal to the total mass absorbed. This configuration, which has few concentrated sources and a distributed sink (or viceversa), can be recognized in several natural networks, such as blood flow or water flow in a river basin.

Three experiments on the frog tongue data are presented. The first two experiments resemble those presented in~\cref{sec:y-net}. One experiment simulates the case where the source of corruption is loss of data that causes network disconnection, while the other simulates the presence of artifacts in the observed data. The observed network $\Obs{\Img}$ used in the second experiment is reported in~\cref{fig:frog-network-shifted}.

In the third experiment, we present an alternative approach to reconstruct the arterial network of the frog tongue. The idea is to estimate a conductivity $\Obs{\Tdens}$ (shown in~\cref{fig:frog-network-mupou}) from the binary image $\Obs{\Img}$ by using the Hagen-Poiseuille law in~\cref{eq:pouseille-law}. We use these enhanced data within~\cref{eq:inpainting-niot-intro} with the identity map, where the NIOT approach becomes a model-driven inverse problem written in terms of the conductivity $\Tdens$.

\subsubsection{Binary images with lost data}
\label{sec:frog-network}
The first network inpainting experiment consists of reconstructing the binary map that describes the support of the arterial network. This image is corrupted using a mask composed of seven circles that removes certain portions of the network. Their location is reported in~\cref{fig:cohnheim}.

We first calibrate our model using $\lambda=0$. The corresponding optimal $\OptTdensH$, reported in~\cref{fig:frog-network-opttdens}, is composed of a high-conductivity skeleton and a series of small channels that fill the domain of the sink term $\Sink$. The maximum value of $\OptTdensH$ is $2e-2$. So we set $\Scaling=5e1$ to match the value of a binary representation of the arterial network.

In~\cref{fig:frog-network} we report the results obtained by fitting the observed data using different values of $\lambda$ and $\confidence=\confmask$ or $\confidence=1$. The initial data $\Tdens_0=1$.

{
\def \fraction {0.32}
\def \extrah {0.02em}
  \begin{figure}
    \centering
    \begin{tabular}{|@{}c@{}|@{\hspace{\extrah}}c@{\hspace{\extrah}}|@{\hspace{\extrah}}c@{\hspace{\extrah}}|@{\hspace{\extrah}}c@{\hspace{\extrah}}|}
      \hline
      &$\lambda=1e-3$ & $\lambda=1e-1$ & $\lambda=1e0$ 
      \\
      \hline
      \raisebox{-.5\normalbaselineskip}[0pt][0pt]{\rotatebox[origin=c]{90}{$\confidence=\confmask$}}
      &
      \adjustbox{valign=m,vspace=0pt}{\includegraphics[trim={.4cm 2.6cm 0.2cm 2.6cm},clip, width=\fraction\columnwidth,valign=t]{{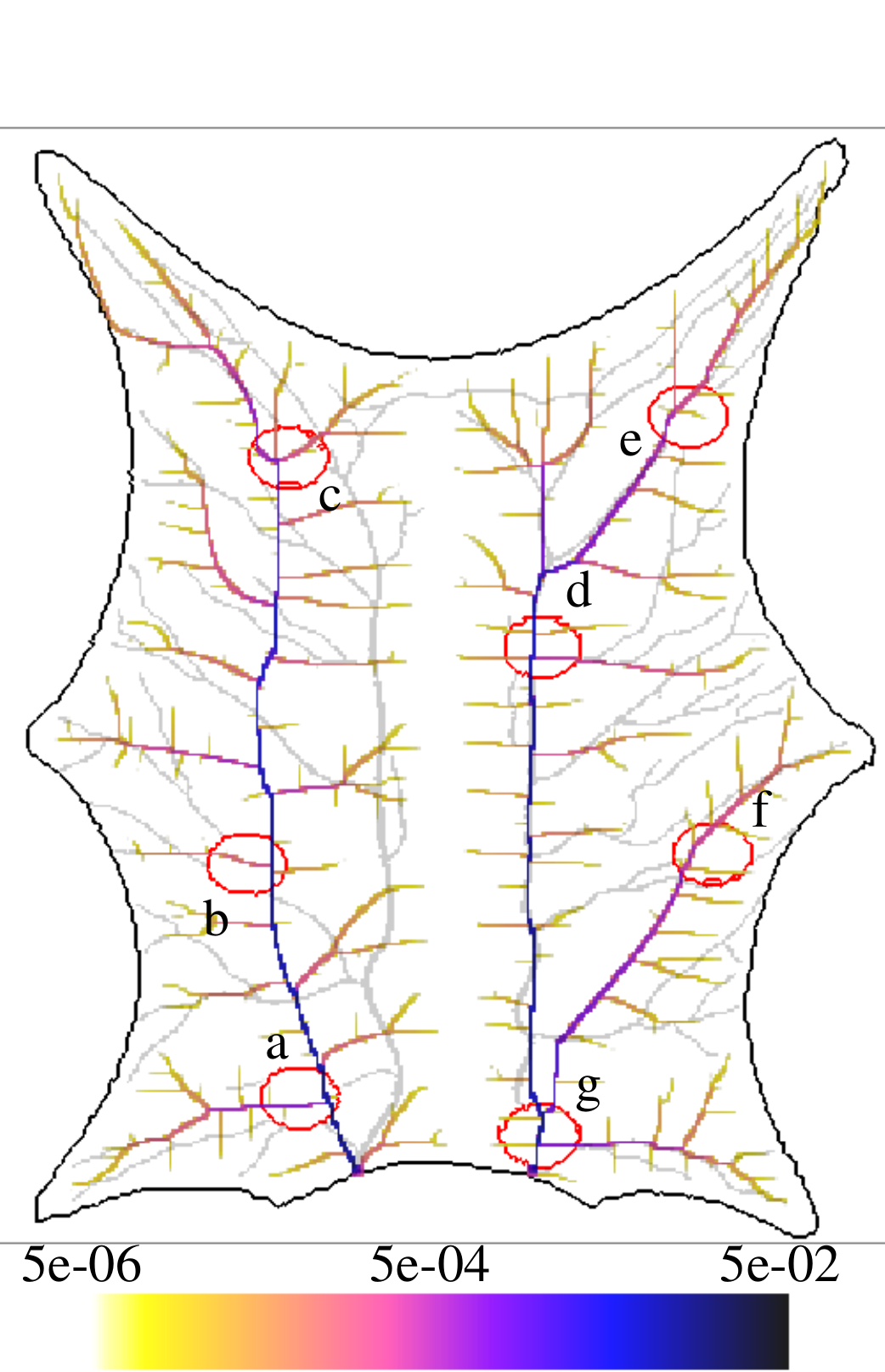}}}
      &
      \adjustbox{valign=m,vspace=0pt}{\includegraphics[trim={.4cm 2.6cm 0.2cm 2.6cm},clip, width=\fraction\columnwidth,valign=t]{{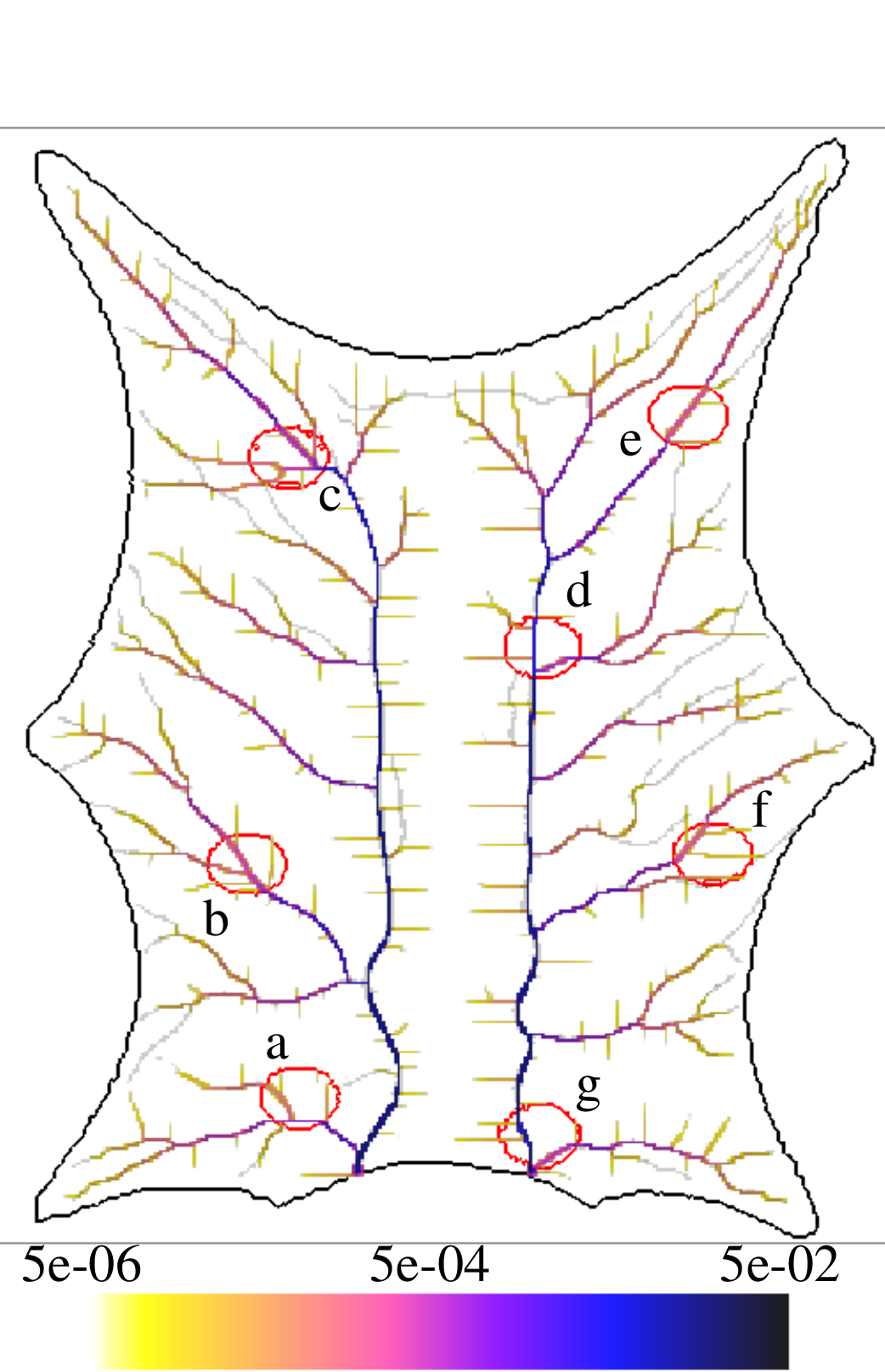}}}
      &
      \adjustbox{valign=m,vspace=0pt}{\includegraphics[trim={.4cm 2.6cm 0.2cm 2.6cm},clip, width=\fraction\columnwidth,valign=t]{{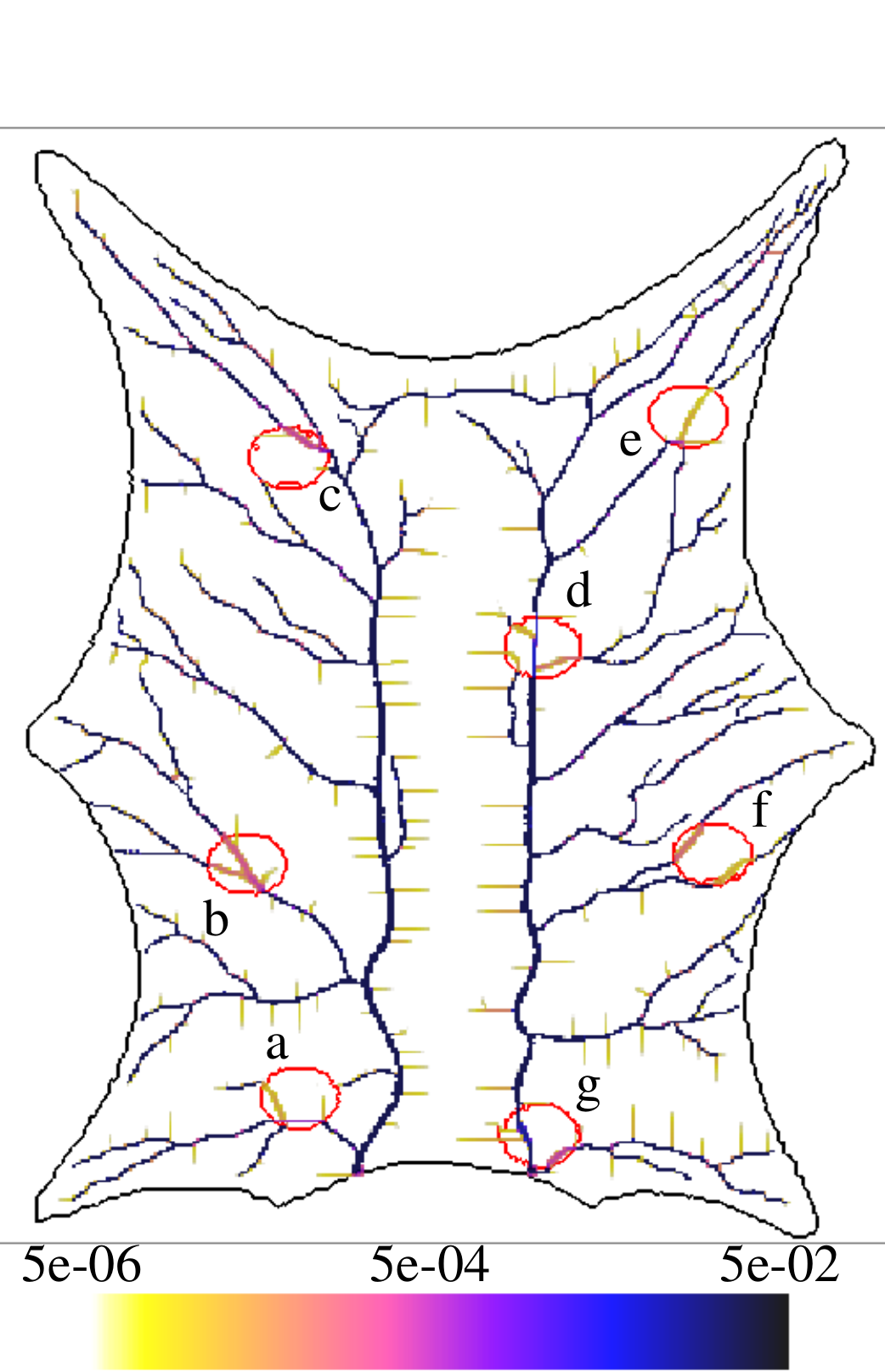}}}
      \\
      \hline
      \raisebox{-.5\normalbaselineskip}[0pt][0pt]{\rotatebox[origin=c]{90}{$\confidence=1$}}
      &
      \adjustbox{valign=m,vspace=0pt}{\includegraphics[trim={.4cm 2.6cm 0.2cm 2.6cm},clip, width=\fraction\columnwidth,valign=t]{{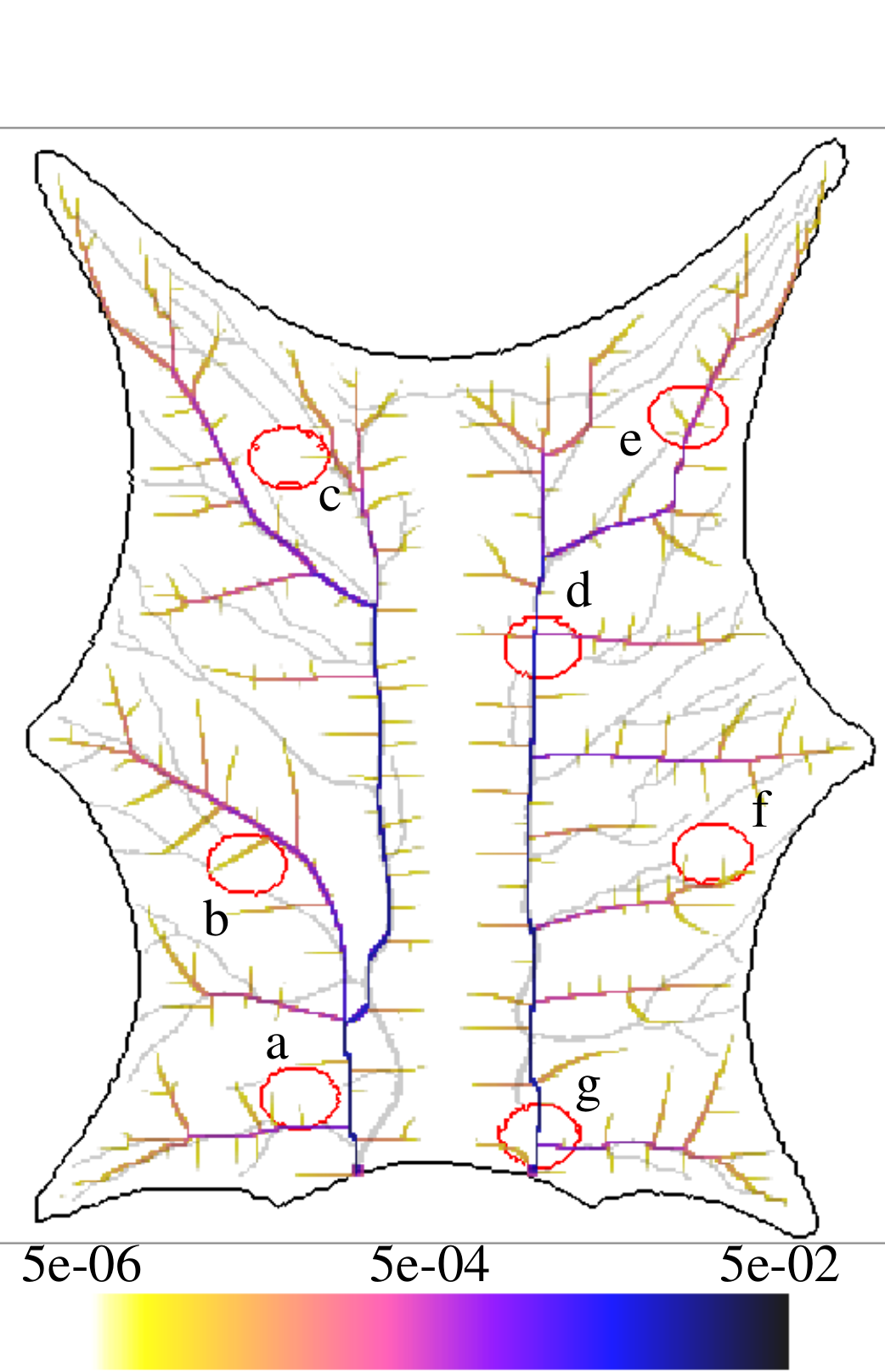}}}
      &
      \adjustbox{valign=m,vspace=0pt}{\includegraphics[trim={.4cm 2.6cm 0.2cm 2.6cm},clip, width=\fraction\columnwidth,valign=t]{{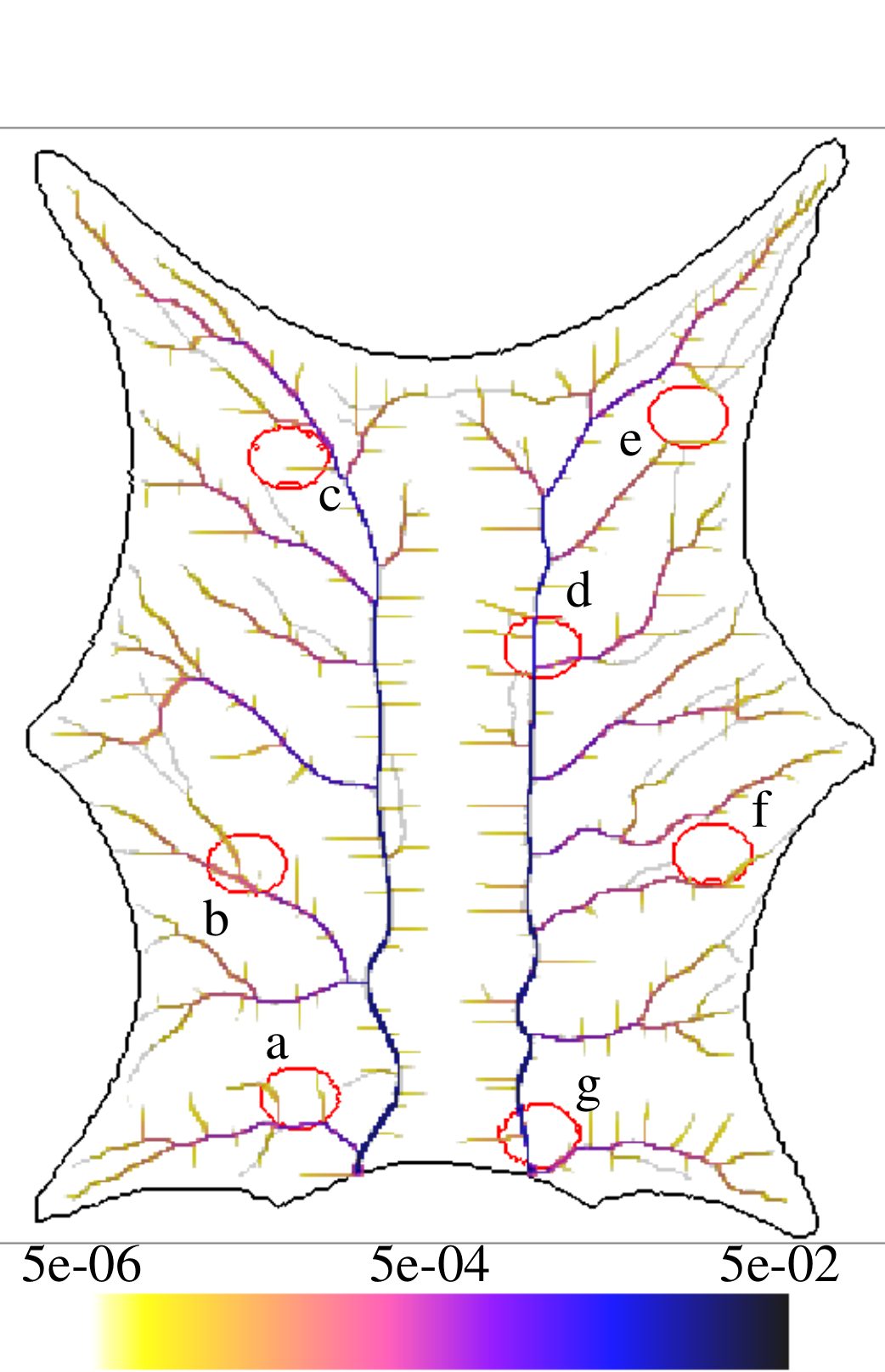}}}
      &
      \adjustbox{valign=m,vspace=0pt}{\includegraphics[trim={.4cm 2.6cm 0.2cm 2.6cm},clip, width=\fraction\columnwidth,valign=t]{{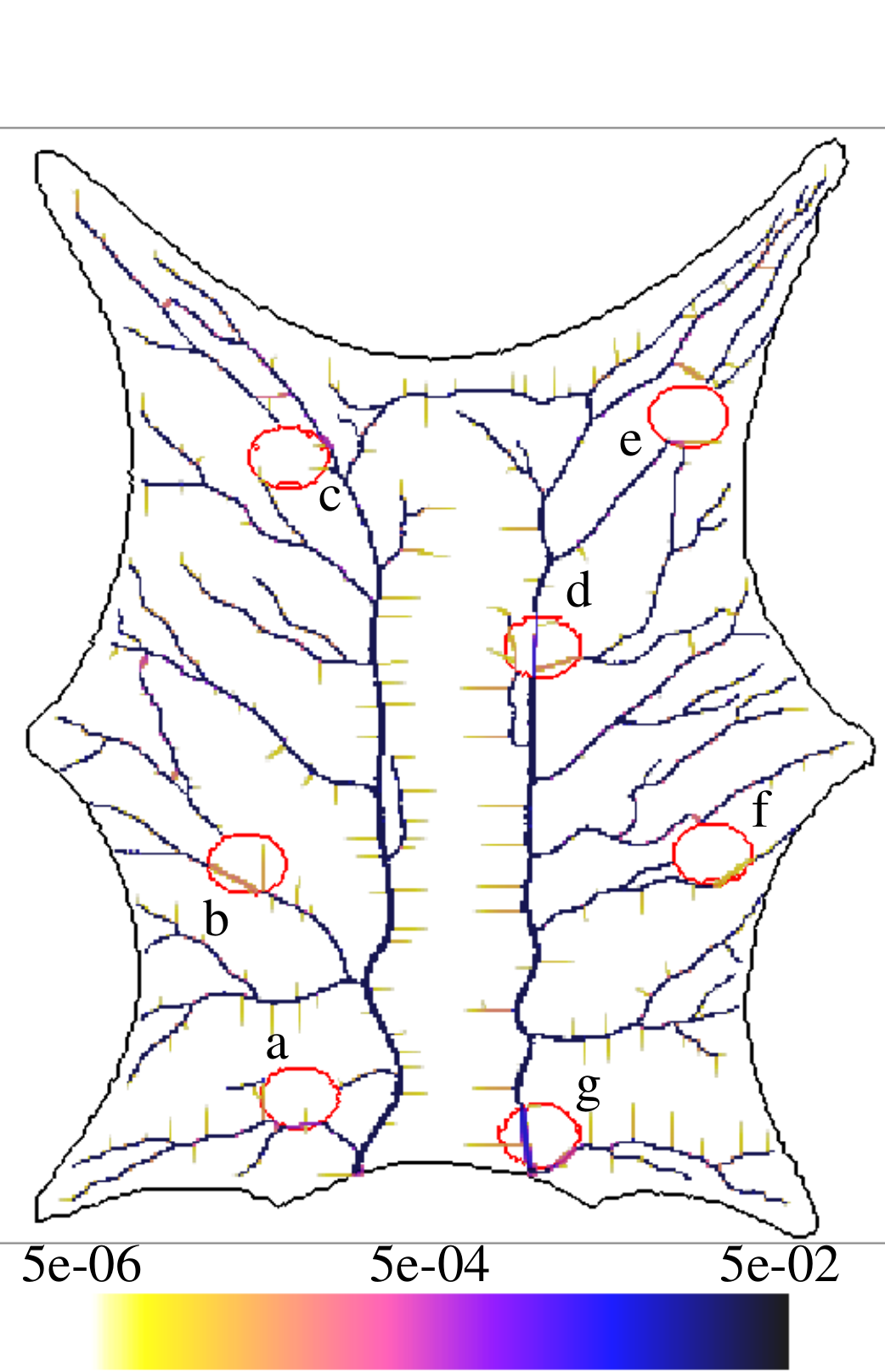}}}
    \\
    \hline
    \multicolumn{4}{c}{
    \quad
    \includegraphics[trim={0.cm 0cm 0.cm 24.83cm},clip,width=0.33\columnwidth,valign=t]{{imgs/medium/mask02/matrix_nref0_femDG0DG0_gamma5.0e-01_wd1.0e+00_wr0.0e+00_netnetwork_ini0.0e+00_confONE_mu2iidentity_scaling5.0e+01_methodte_tdens.pdf}}
    }
  \end{tabular}
  \caption{Spatial distribution $\OptTdensH$ for $\lambda=1e-3,1e-1,1e0$ 
  and $\confidence=\confmask,1$. We use the identity map and $\Tdens_0=1$. The observed data $\Obs{\Img}$ (reported in light gray) is the binary map of the arterial network removed within the mask.}
  \label{fig:frog-network}
  \end{figure}
}

In all experiments, the curved portion of the network within circle (d) and the parallel channels within circle (a) are consistently not recovered. This can be attributed to the fact that within the mask, these configurations are not optimal for the $\Reg$ functional. 

For $\lambda=1e-3$, the discrepancy term already affects the resulting network. The main right branch is approximately fitted for both $\confidence=\confmask$ and $\confidence=1$. However, despite the fact that these two confidence measures are different only within the mask, the left branch resulting networks are completely different. This can again be attributed to the non-convexity of the functional $\Reg$ and the usage of gradient-based optimization methods. 

Within circle (g), note that a good approximation of the bifurcation point is created in all cases, even in the case $\lambda=0$ shown in~\cref{fig:frog-network}. In fact, without any channel within the local mask, the flux cannot be carried from the (right) source to the sink.

At $\lambda=1e-1$, the Y-shaped network within circle (b) is recovered for both $\confidence=\confmask$ and $\confidence=1$. In the first case, the reconstructed network is not totally accurate, but the main structure of the network is recovered within all circles, excluding the inaccuracies already mentioned within circles (a) and (d). In the second case, the penalization imposed by having $\confidence=1$ within the mask removes the channels within circles (c), (e) and (f). However, outside the circular masks, the overall structure of the network is recovered with good quality.

For $\lambda=1e0$ and $\confidence=\confmask$, the reconstructed network fits the observed data accurately outside the mask at the cost of having some errors within the circle (e). 

For $\lambda=1e0$ and $\confidence\equiv 1$, an inconsistent network is obtained within circle (b). No network is created here, and the tail of corrupted branched is ``fed'' by another sub-branch. This reconstructing error is not particularly satisfactory from a modeling point of view, because this would lead to flux flowing in the wrong direction. For $\lambda=1e0$ and $\confidence=\confmask$, the reconstructed network fits the observed data very well outside the mask at the cost of having some errors within circle (e).

We attribute these phenomena to the combination of different factors:
\begin{itemize} 
  \item Within the mask we have that $\confidence\equiv1$ and $\Obs{\Img}=0$, hence  the discrepancy term is penalizing any presence of the network, similarly to what happened in the previous section. This effect can also be seen with circles (c), (e), and (f), where channels outside the mask are selected.
  \item The approach is currently blind to the direction of the flux, since this information is not encoded in our inpainting approach, neither in the data nor in the model.
  \item The fitted data is a binary map. It contains only information about the support of the network, but not its hierarchical structure, from the root to the leaves.
\end{itemize}
In~\cref{sec:frog-mu}, we will try to mitigate the effect of the last factor by extracting more information from the data. 

\subsubsection{Binary images with artifacts}
\label{sec:frog-shifted}
We now present the experiment in which the data are corrupted by the presence of artifacts, where $\Obs{\Img}$ is shown in~\cref{fig:frog-network-shifted}. This image is obtained by shifting the original network within the rectangular area and merging it with the original network. The data are also corrupted within the circular masks shown in~\cref{fig:frog-network}. 

In this false positive regime, the reconstruction of the network is also required to prune redundant channels and obtain the main structure of the network. We use $\Tdens_0=1$ and $\confidence=1$, and we show the results in~\cref{fig:shifted} the results obtained for different $\lambda$.

{
\def \fraction {0.32}
\def \extrah {0.02em}
\begin{figure}
    \centering
    \begin{tabular}{|@{\hspace{\extrah}}c@{\hspace{\extrah}}|@{\hspace{\extrah}}c@{\hspace{\extrah}}|@{\hspace{\extrah}}c@{\hspace{\extrah}}|}
      \hline
        $\lambda=1e-2$ & $\lambda=1e-1$ & $\lambda=1e0$
      \\
      \hline
      \adjustbox{valign=m,vspace=0pt}{ \includegraphics[trim={0.cm 2.65cm 0.cm 2.65cm},clip,width=\fraction\columnwidth]{{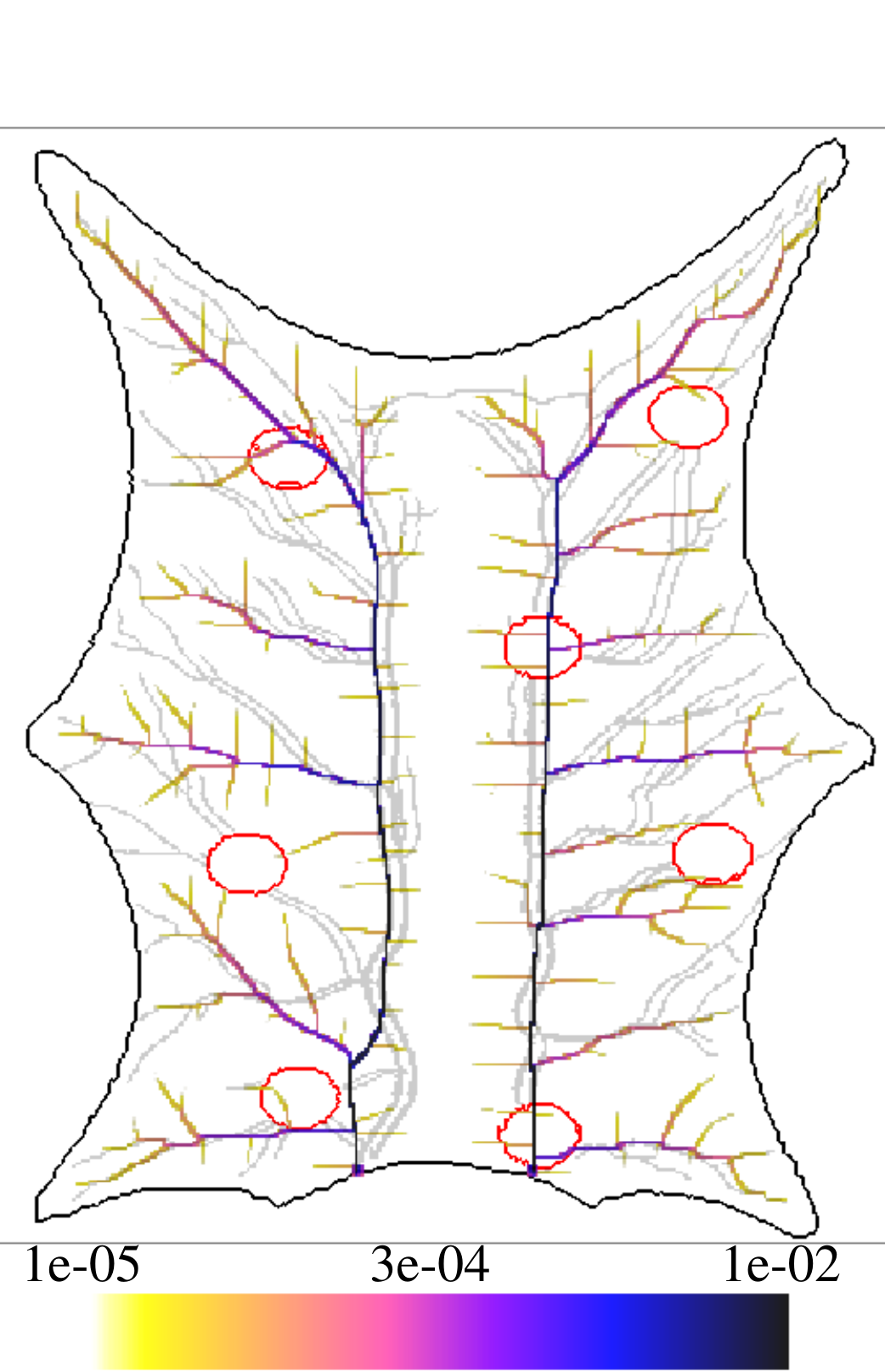}}}
      &
      \adjustbox{valign=m,vspace=0pt}{ \includegraphics[trim={0.cm 2.65cm 0.cm 2.65cm},clip,width=\fraction\columnwidth]{{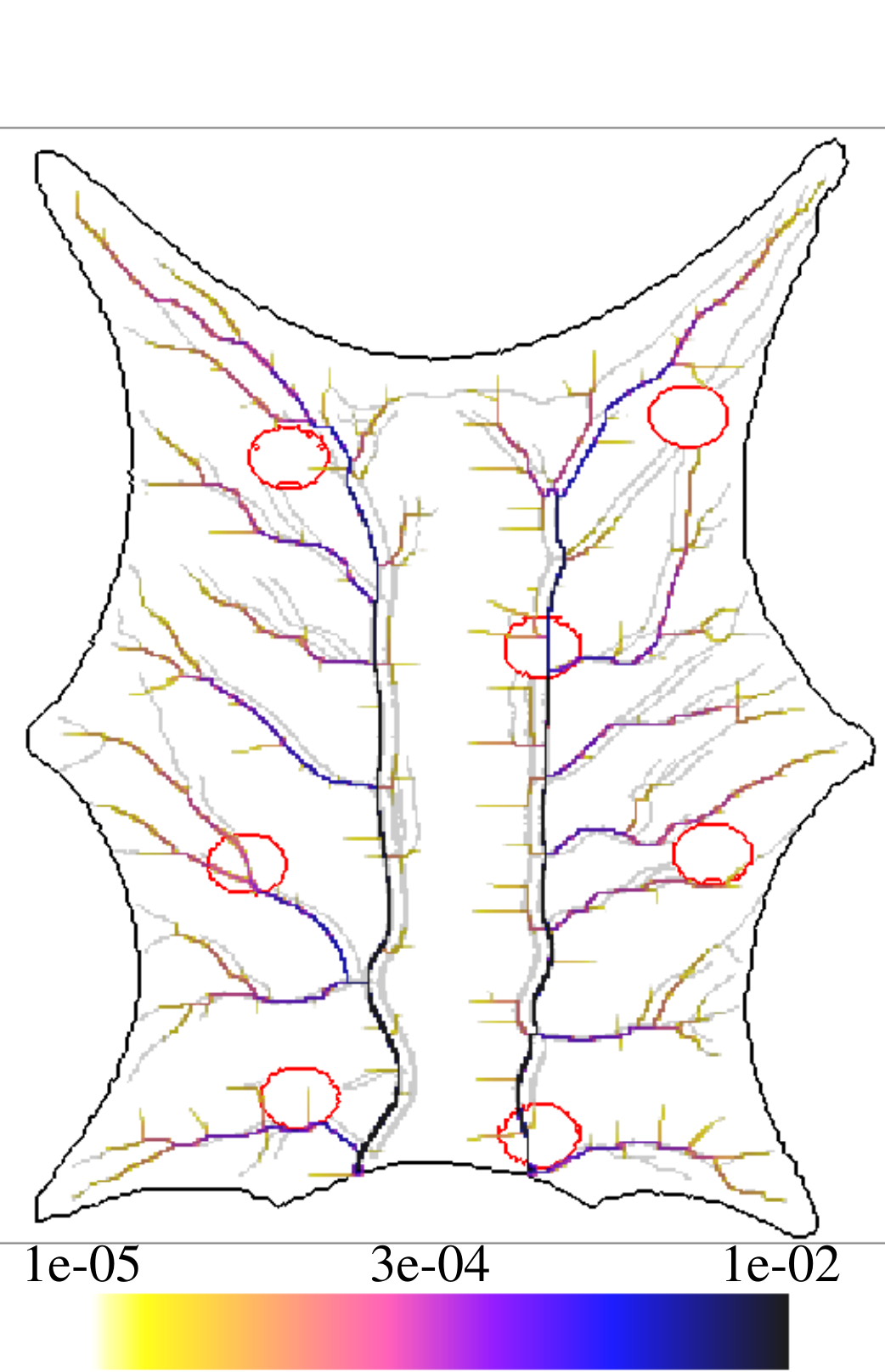}}}
      &
      \adjustbox{valign=m,vspace=0pt}{ \includegraphics[trim={0.cm 2.65cm 0.cm 2.65cm},clip,width=\fraction\columnwidth]{{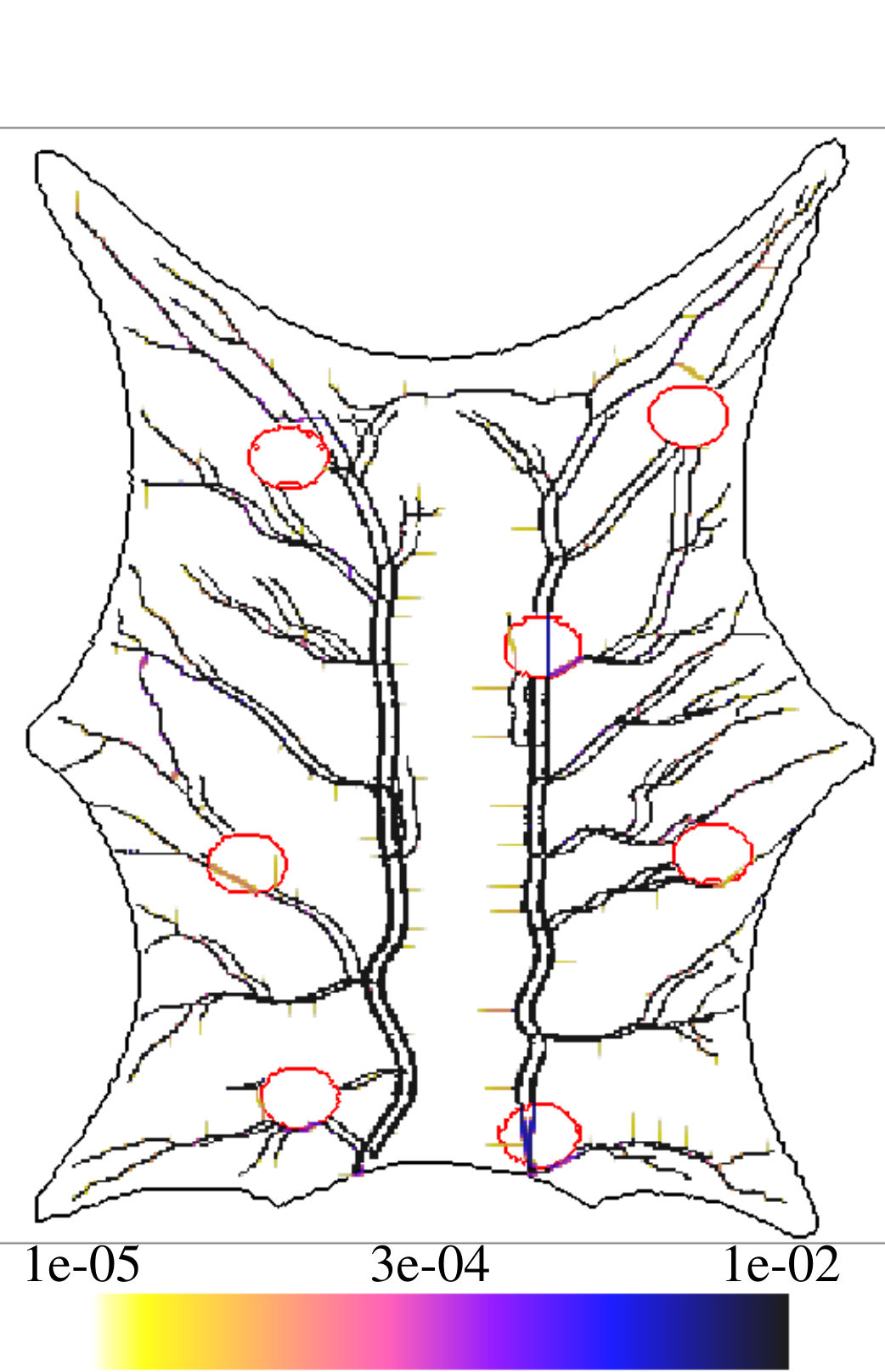}}}
      \\
      \hline
    \multicolumn{3}{c}{
    \includegraphics[trim={0.cm 0cm 0.cm 24.83cm},clip,width=0.33\columnwidth,valign=t]{{imgs/medium/mask02/matrix_nref0_femDG0DG0_gamma5.0e-01_wd1.0e-02_wr0.0e+00_netnetwork_shifted_inside2_ini0.0e+00_confONE_mu2iidentity_scaling5.0e+01_methodte_tdens.pdf}}
    }
  \end{tabular} 
  \caption{
    Spatial distribution of $\Obs{\Img}$ (top left panel) and of $\RecTdensH$ for different $\lambda=5e-2,1e-1,1e0$ ($\confidence=1$, $\Tdens_0$, and identity map). The observed image $\Obs{\Img}$ is reported in light gray.
  }
  \label{fig:shifted}
\end{figure}

At $\lambda=1e-2$ we already recover a good approximation of the original network, with the main structure of the network being recovered but with a strong pruning of the redundant channels. At $\lambda=1e-1$, the network we obtain even better results, like the reconstruction of the Y-shaped network in the circle (b), but we start seeing some overlay of small channels with the shifted network, in particular along the two main channels. For $\lambda=1.0$ we fall into the overfitting regime.

In the context of pruning, let us also note how reducing $\lambda$ cuts the vessel at the top of the image that connects the left and right basins of the frog arterial structure. In fact, this redundant structure cannot be taken into account in the current \BT{} penalization.

\subsubsection{Enhanced image data}
\label{sec:frog-mu}
In this section, we present an alternative approach to reconstructing the arterial network of the frog tongue. The idea is to use as observed data $\Obs{\Img}$ an estimate of the conductivity of the network. We can extract the conductivity from the image using the local thickness (see \cite{hildebrand1997new} for a precise definition of this quantity) of the network and the expression $\PouTdens = \PouseilleConstant r^p$ derived from the Hagen-Poiseuille law.

The observed data $\Obs{\Img}$ is obtained through the following steps (we used ImageJ library~\cite{schindelin2012fiji} to perform the first two).
\begin{enumerate}
  \item We compute the local thickness $t$ of the binary representation $\True{\Img}$ of arterial network.
  \item We get a skeleton $s$ of $\True{\Img}$, which is a binary map inside the support of $\True{\Img}$ that is only 1 pixel wide.
  \item We define an approximate Dirac measure of the 1d skeleton $s$ as $S=s/h$, where $h$ is the width of 1 pixel.
  \item We compute the conductivity $\PouTdens=\PouseilleConstant r^p S$ with $\PouseilleConstant=5.0e+02$ and $p=3$.  
  \item Finally, we apply a \PM{} map to $\PouTdens$, obtaining the data $\Obs{\Img}$ to be fitted. The exponent $m$ and the time $t^*$ that define the \PM{} map are those given in~\cref{eq:mt}.
  The value $\PouseilleConstant=5.0e+02$ was experimentally found so that $\Obs{\Img}$ and the optimal $\OptTdensH$ reported in the right panel of~\cref{fig:cohnheim} have the same order of magnitude. 
\end{enumerate}
The last preprocessing step provides a smoother representation of the conductivity having the same support of the arterial network. 

The corruption step involved in the application of the mask shown in~\cref{fig:cohnheim} and the resulting conductivity $\Obs{\Tdens}$ is reported in~\cref{fig:frog-network-mupou}. Then we use these data within the optimization problem in~\cref{eq:inpainting-niot-intro} using $\Obs{\Img}=\Obs{\Tdens}$ and the identity map with $\alpha=1$. In~\cref{fig:frog-mu}, we report the spatial distribution of the conductivity $\RecTdensH$ using increasing values of $\lambda$ and $\confidence=\confmask,1$. In this case, the values of $\lambda$ are higher than those presented in~\cref{sec:frog-network} as the fitted data have lower values (the maximum value is $1e-2$).

The results are consistent with the previous experiments, the network is progressively reconstructed as $\lambda$ increases. In this case, however, we see that the Y-shaped network with circle (b) is reconstructed even for $\lambda=1e0$ and $\confidence=1$. This example shows that the process of retrieving the conductivity from the image can be beneficial for the inpainting process proposed in this paper.

{
\def \fraction {0.3}
\def \extrah {0.02em}
\begin{figure}
  \begin{tabular}{|@{}c@{}|@{\hspace{\extrah}}c@{\hspace{\extrah}}|@{\hspace{\extrah}}c@{\hspace{\extrah}}|@{\hspace{\extrah}}c@{\hspace{\extrah}}|}      
  \hline
  &$\lambda=1e1$ & $\lambda=1e3$ & $\lambda=1e4$ 
      \\ 
      \hline
      \raisebox{-.5\normalbaselineskip}[0pt][0pt]{\rotatebox[origin=c]{90}{$\confidence=\confmask$}}
      &
      \adjustbox{valign=m,vspace=0pt}{\includegraphics[trim={.4cm 2.6cm 0.2cm 2.6cm},clip, width=\fraction\columnwidth,valign=t]{{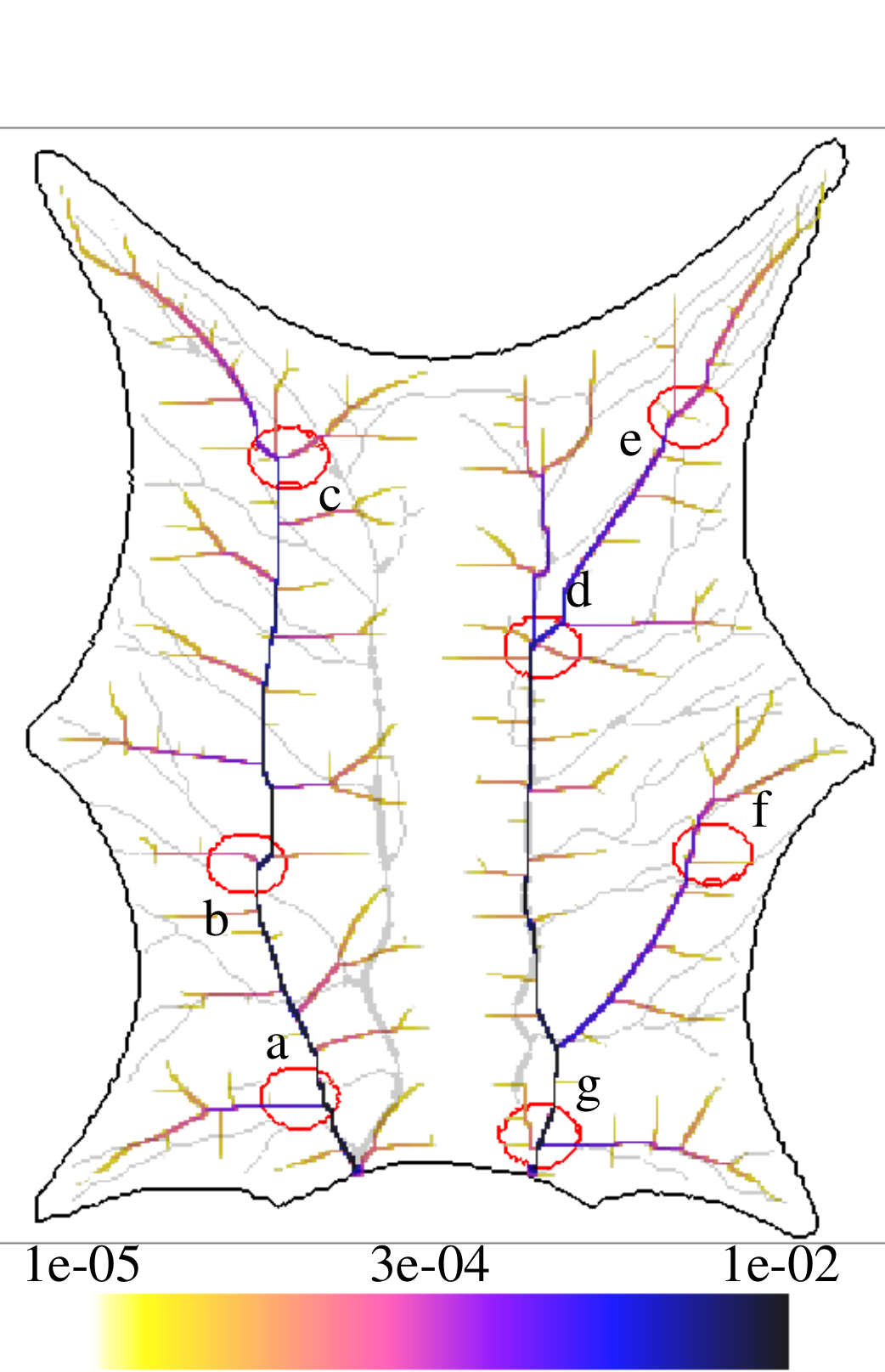}}}
      &
      \adjustbox{valign=m,vspace=0pt}{\includegraphics[trim={.4cm 2.6cm 0.2cm 2.6cm},clip, width=\fraction\columnwidth,valign=t]{{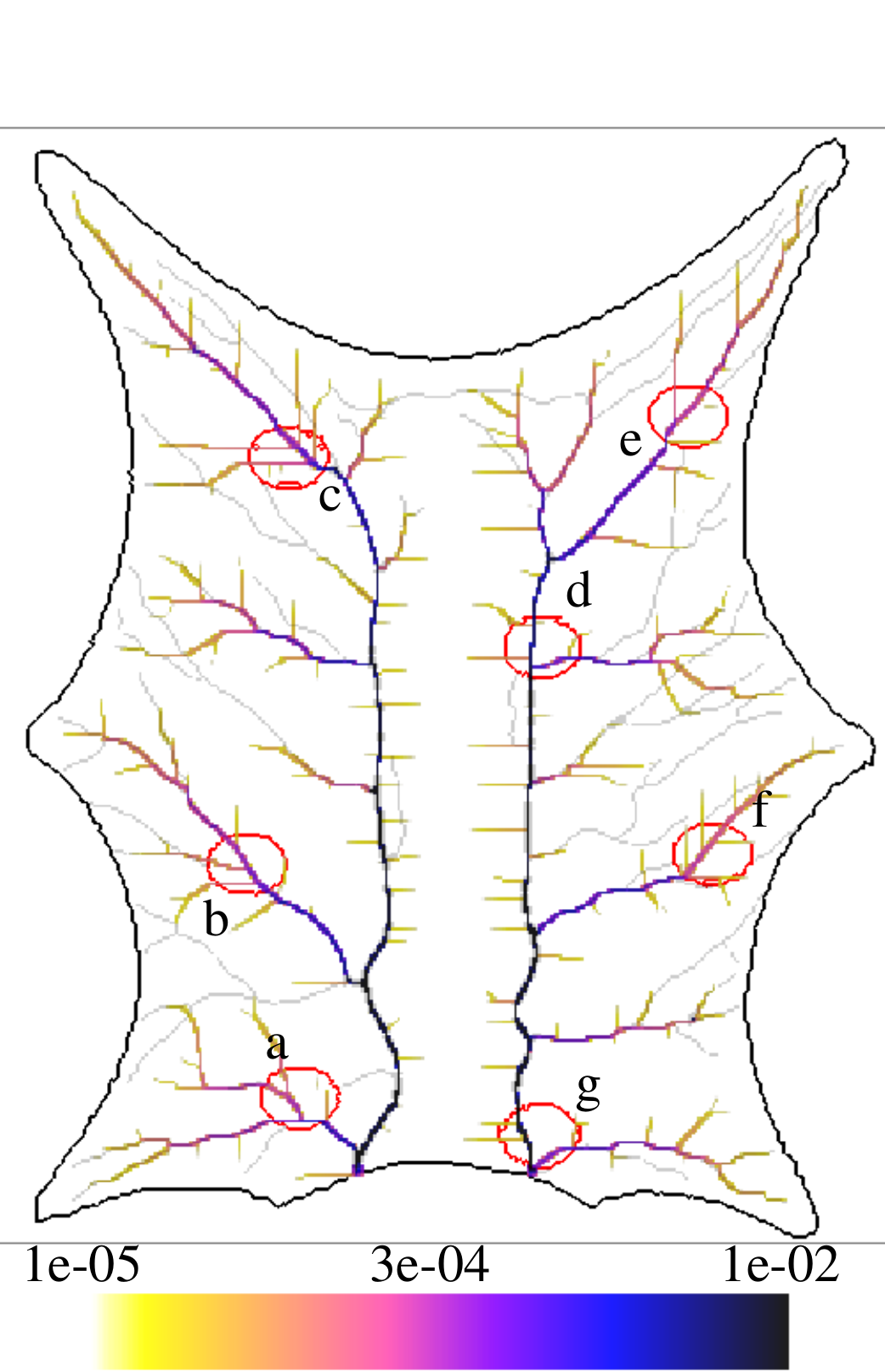}}}  
      &
      \adjustbox{valign=m,vspace=0pt}{\includegraphics[trim={.4cm 2.6cm 0.2cm 2.6cm},clip, width=\fraction\columnwidth,valign=t]{{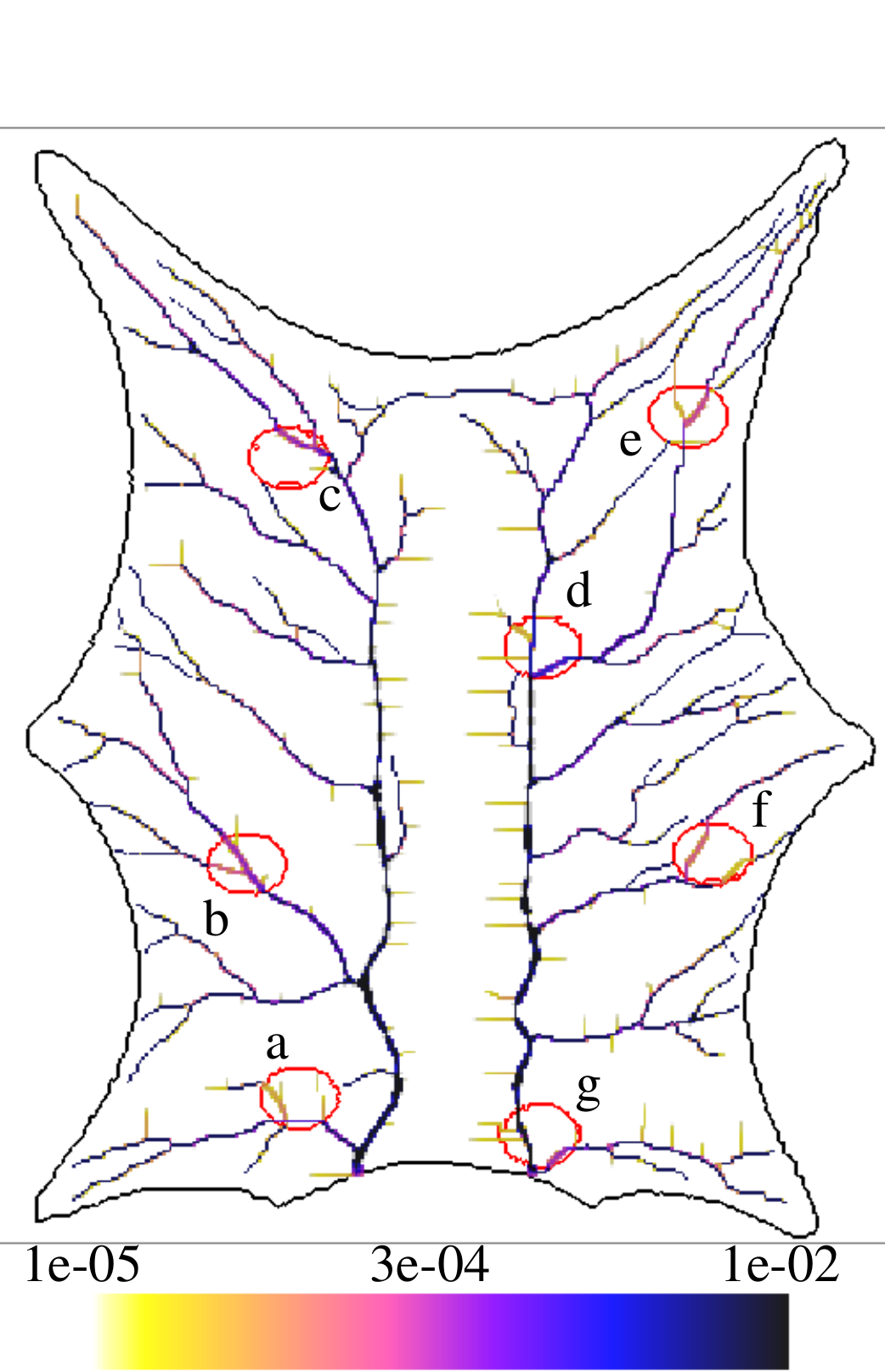}}}
      \\
      \hline
      \raisebox{-.5\normalbaselineskip}[0pt][0pt]{\rotatebox[origin=c]{90}{$\confidence=1$}}
      &
      \adjustbox{valign=m,vspace=0pt}{\includegraphics[trim={.4cm 2.6cm 0.2cm 2.6cm},clip, width=\fraction\columnwidth,valign=t]{{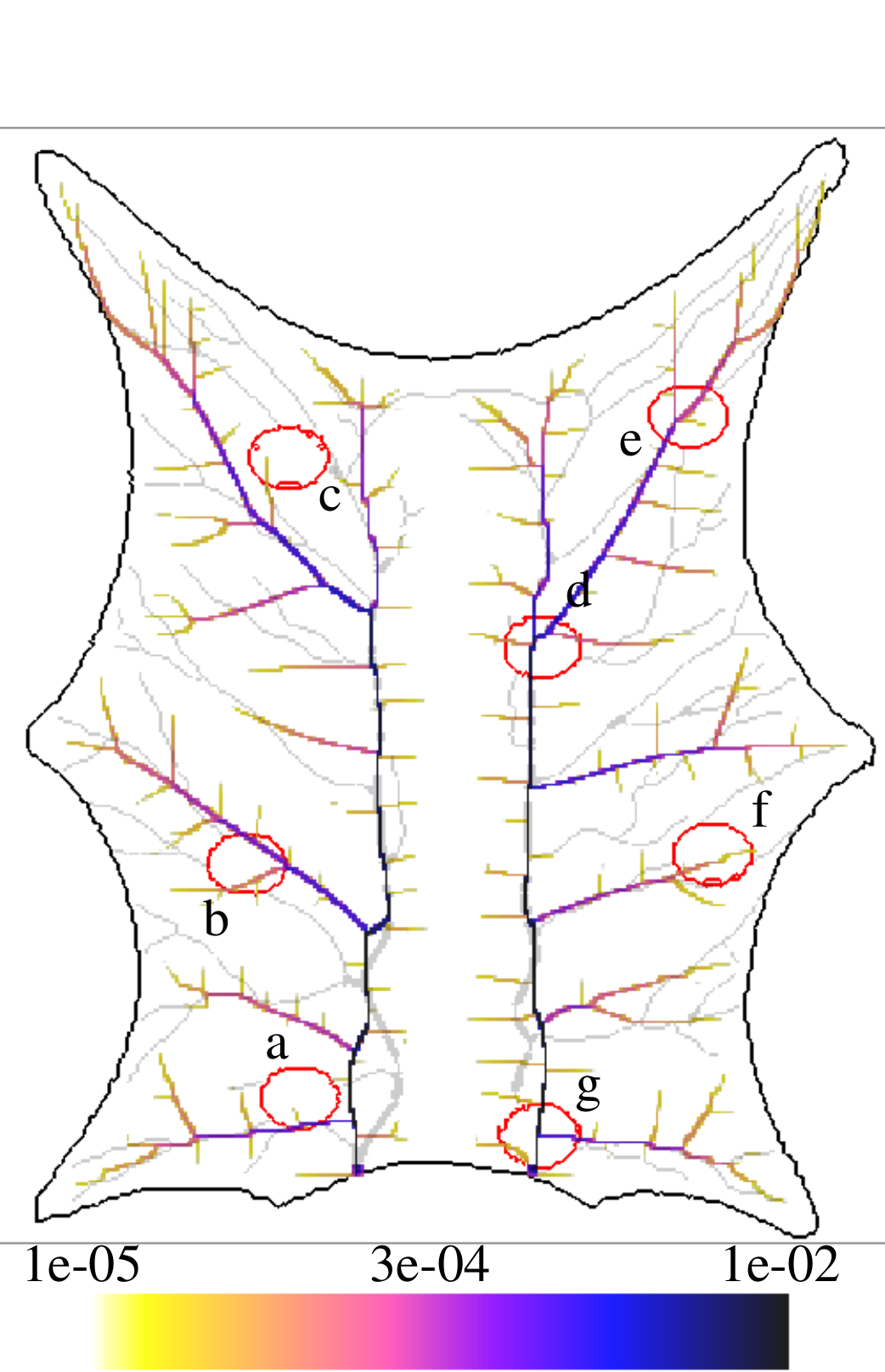}}}
      &
      \adjustbox{valign=m,vspace=0pt}{\includegraphics[trim={.4cm 2.6cm 0.2cm 2.6cm},clip, width=\fraction\columnwidth,valign=t]{{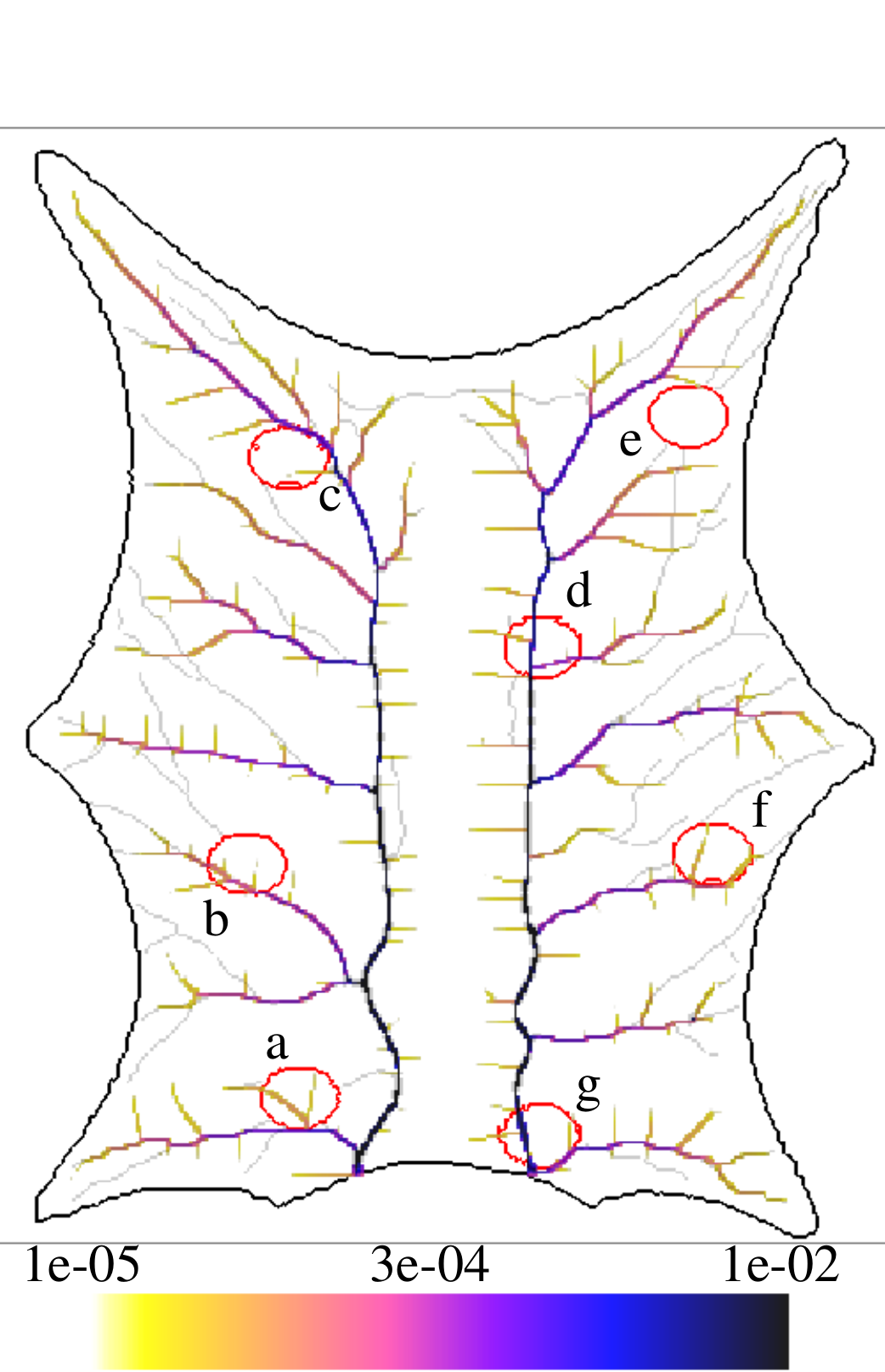}}}
        &
      \adjustbox{valign=m,vspace=0pt}{\includegraphics[trim={.4cm 2.6cm 0.2cm 2.6cm},clip, width=\fraction\columnwidth,valign=t]{{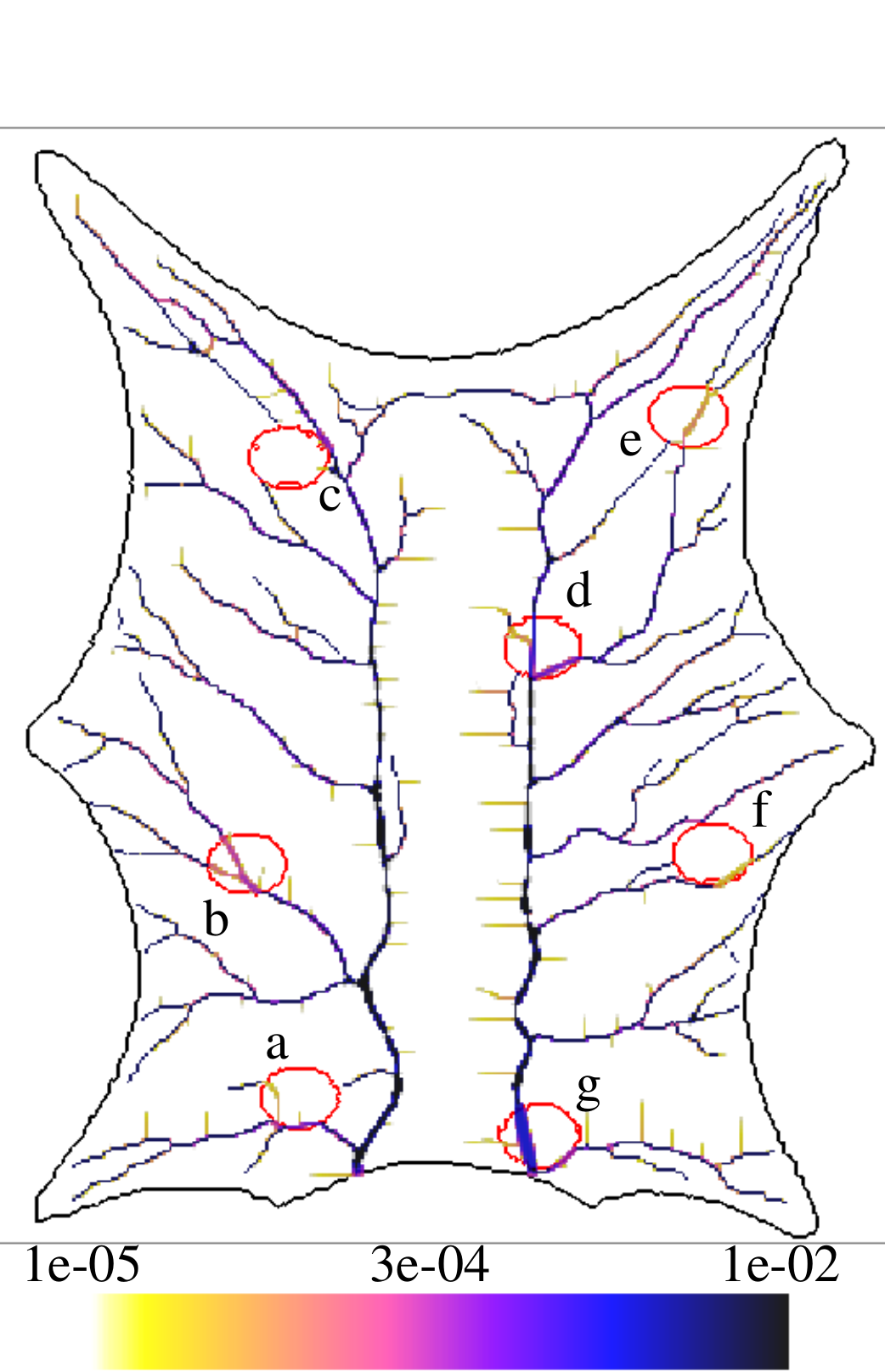}}}
      \\
      \hline
      \multicolumn{4}{c}{
      \qquad
      \includegraphics[trim={0.cm 0cm 0.cm 24.83cm},clip,width=0.33\columnwidth,valign=t]{{imgs/medium/mask02/matrix_nref0_femDG0DG0_gamma5.0e-01_wd1.0e+04_wr0.0e+00_netmup3.0e+00zero5.0e+02_ini0.0e+00_confONE_mu2iidentity_scaling1.0e+00_methodte_tdens.pdf}}
      }
    \end{tabular}
    \caption{Spatial distribution $\RecTdensH$ for $\lambda=1e1,1e3,1e4$ ($\confidence=1$, $\Tdens_0$), and identity map. The observed data $\Obs{\Img}$ is the conductivity $\Obs{\Tdens}$ is shown~\cref{fig:frog-network-mupou}. 
    Its support is reported in light gray.}
    \label{fig:frog-mu}
    \end{figure}
}   

\section{Conclusions}
\label{sec:conclusions}
We showed how using the branch-inducing functional $\Reg$ as the regularization functional in variational inpainting methods offers a natural way to encode physical information in the network reconstruction problem.

The proposed approach is particularly effective in restoring connectivity in corrupted networks. This is obtained by imposing a flux passing through the network, and it is possible only under the assumption that a source term $\Source$ and a sink term $\Sink$ can be identified. Even rough estimations presented in the examples led to a good reconstruction of the network.

Characterizing a network as a support of conductivity $\Tdens$ allows for easy merging of the modeling information and observed data. These two components of the model are linked by the map $\ImgTd$. In this paper, we presented only two maps, the simple identity map and the \PM{} map. The latter is a useful tool in the case of the network for which the Hagen-Pouseille law given in~\cref{eq:pouseille-law} holds.

Unfortunately, the approach requires the adjustment of several parameters, summarized in~\cref{tab:parameters}. In most of the examples presented, these parameters were estimated experimentally, but in general their tuning should be driven by the problem at hand.

The experiments presented in the papers suggest that this tuning should also be driven by the way data are corrupted. If good data are available and the corruption is only due to data loss, high values of $\lambda$ can be used. In this setting, it may be beneficial to choose the initial data $\Tdens_0$ derived from the observed data and exploit the ``memory'' effect described in~\cref{sec:tdensini-obs}.
On the other hand, if the corruption is due to the presence of false positive data, small values of $\lambda$, $\confidence=1$, and $\Tdens_0=1$ should be used.
In this setting, the proposed approach has also shown interesting capabilities in removing artifacts from highly corrupted data, performing this selection according to an optimization principle, as shown in~\cref{sec:artifacts,sec:frog-shifted}.

The proposed model can be extended in several directions, enriching both its fitting and modeling component. In the first component, we can include further terms, like the flux or the pressure, when this information can be extracted from the data. Different discrepancy metrics can be used in place of the weighted $L^2$-norm considered here. Natural candidates are \OT{}-based distances, similar to the unbalanced metrics described in~\cite{todeschi2024unbalanced}, since these distances can deal with densities concentrated on support with a small width-to-length ratio. We did not explore this possibility in this paper, but it is a natural extension of the work presented here.

The second component of the NIOT approach can be enriched, for example, by including a localized penalization of the term $\int_{\Domain} \Tdens^{\gamma}\dx$ in~\cref{eq:lyap}, similarly to what was done in~\cite{Facca-et-al:2018}. Networks having loops can be included using fluctuating source and sink terms as in~\cite{PhysRevE.107.024302}.

%% file: appendix.tex
\section{Appendix}
\label{sec:appendix}
Here we present some well-known results on the porous media equation and how to set the exponent $m$ and the time $t^*$ to obtain the desired properties of the \PM{} map introduced in~\cref{sec:map}.

\subsection{Porous media equation and the Barenblatt solution}
\label{sec:pm}
\newcommand{\Dimp}{n}
When the initial data the porous media equation in~\cref{eq:porous-media} is given as $\pmu_0(x)=M\delta_{0}(x)$ with $M>0$, then there exists a closed-form solution of the problem, known as the Barenblatt solution~\cite{barenblatt1972}. It is given by
\begin{equation}
    \label{eq:barenblatt}
    \pmu(t,x) = t^{-\alpha}(A-B|x|^2t^{-2\beta})_{+}^{\frac{1}{m-1}} 
\end{equation}
where $\alpha$, $\beta$, $A$, and $B$ are positive constants that can be expressed in terms of 
the dimension $\Dimp$, the exponent $m$, and the initial mass $M$. Following~\cite{fasano2006problems}, they can be written as
\begin{equation*}
  \alpha = \beta \Dimp, \quad 
  \beta = \frac{1}{2+(m-1)\Dimp}, \quad 
  A=\left(\frac{B^{\Dimp/2}M}{w_\Dimp z_{m,\Dim}}\right)^{2\beta(m-1)} \quad B=\beta \frac{m-1}{2m},
  \end{equation*}
where $w_\Dimp$ is the volume of the $\Dimp$ sphere and $z_{m,\Dimp}$ is given by
\[
  z_{m,\Dimp}=\int_{0}^{\pi/2}[\cos (\theta)]^{\frac{m+1}{m-1}}[\sin(\theta)]^{\Dimp-1}d \theta.
\]

The support of the Barenblatt solution $\pmu(t,x)$ is a sphere expending in time with radius
\begin{equation}
    \label{eq:radius-support}
    r(t)=t^{\beta}B^{-1/2}A^{1/2}.
\end{equation}

\subsection{Exponent $m$ and time $t^*$ for the \PM{} map}
We want to find the exponent $m$ and the time $t^*$ such that, given a singular measure of mass $M$ centered at $0$ at time $t=0$, the evolution of this mass according to~\cref{eq:porous-media} 
\[
M = \PouseilleConstant r(t^*)^{p}
\]
for any $M$. This is equivalent to saying such that the following equality holds
\begin{equation*}
  r(t^*)=(t^*)^{\beta} B^{-1/2} A^{1/2} =
  (t^*)^{\beta}  B^{-1/2}  \frac{B^{n\beta(m-1)/2} M^{\beta(m-1)}}{(w_n,z_{m,n})^{\beta(m-1)} } 
  =\PouseilleConstant^{-1/p}M^{1/p}=r
\end{equation*}
for all $M$. This can be done by matching the exponents for $M$ and proportionality constants, getting the following expression for the exponent $m$ and the time $t$ given by
\begin{equation}
  \label{eq:mt}
  m=\frac{2+p-\Dimp}{p-\Dimp}
  \quad
  t^* = \left(\PouseilleConstant^{-1/p} \frac{(w_n,z_{m,n})^{\beta(m-1)}}{B^\beta}\right)^{1/\beta}.
\end{equation}
Note that in our inpainting problem we need to set $\Dimp=\Dim-1$.